\documentclass[10pt,a4paper]{article}

\RequirePackage[OT1]{fontenc}
\RequirePackage{amsthm,amsmath}
\RequirePackage[numbers,comma,sort&compress]{natbib}
\RequirePackage[colorlinks=true,backref=page]{hyperref}
\hypersetup{
    linkcolor = {blue},
    citecolor = {blue},
    urlcolor = {blue}
    }
\usepackage[utf8]{inputenc}
\usepackage[T1]{fontenc}
\usepackage{times}
\usepackage{verbatim}
\usepackage{fullpage}
\usepackage{pdfpages}
\usepackage{amssymb,latexsym,amsfonts}
\usepackage{dsfont, bm}
\usepackage{color}
\usepackage{graphicx}
\usepackage{multicol}
\usepackage{mathtools}
\usepackage{xspace}
\usepackage{bold-extra}
\usepackage{microtype}
\usepackage{mathptmx} 
\usepackage{stmaryrd}
\usepackage[inline,shortlabels]{enumitem}
\usepackage{mathrsfs}
\usepackage{authblk}
\usepackage{todonotes}
\usepackage{float}
\usepackage{subfigure}
\usepackage{mathtools}

\definecolor{Green}{RGB}{0,150,60} 
\definecolor{violet}{RGB}{200,0,200}
\definecolor{vert}{RGB}{12,120,20}
\definecolor{rouge}{RGB}{200,0,0}
\newcommand{\av}[1]{{{\displaystyle{\langle} {#1}\displaystyle{\rangle}}}}

\numberwithin{equation}{section}
\newtheorem{theorem}{Theorem}[section]
\newtheorem{corollary}[theorem]{Corollary}
\newtheorem{remark}[theorem]{Remark}
\newtheorem{lemma}[theorem]{Lemma}
\newtheorem{proposition}[theorem]{Proposition}
\newtheorem{hyp}{Hypothesis}

\newtheorem{exmp}{Example}[section]
\newtagform{cyan}{\color{blue}(}{)}
\usetagform{cyan}

\newcommand{\law}[1]{\mathscr{L}\!\textit{aw}\left(#1 \right)}
\newcommand{\tlaw}{\mathscr{L}\!\textit{aw}}

\newcommand{\rpm}{\widetilde{N}(ds,dz)}
\newcommand{\rpc}{\nu_s (dz) ds}
\newcommand{\lm}{\nu_s (dz)}
\newcommand{\leqpt}{\preccurlyeq}
\newcommand{\indic}[1]{\mathds{1}_{\{#1\}}}

\newcommand{\Blue}[1]{{\color{blue}{#1}\color{black}\xspace}}
\newcommand{\EE}{\mathbb{E}}
\newcommand{\RR}{\mathbb{R}}
\newcommand{\Var}{\mathbb{V}\!{\rm{ar\,}}}
\newcommand{\oeX}{\overline{X}^{\varepsilon}}
\newcommand{\oeeX}{{\widehat{X}}^{\varepsilon}}
\newcommand{\ueX}{\underline{X}^{\varepsilon}}
\newcommand{\wtX}{\widetilde{X}}

\newcommand{\Sint}{\int_{-\infty}^\infty}
\newcommand{\Eint}{\int_{B(\varepsilon)}}
\newcommand{\WEint}{\int_{\mathbb{R}\backslash B(\varepsilon)}}

\newcommand{\normP}[1]{\left\| #1 \right\|_{L^p(\Omega)}}
\newcommand{\normPt}[1]{\| #1 \|_{L^p(\Omega)}}

\newcommand{\cpbdg}{C_p^{\tiny{\textit{BDG}}}}

\newcommand{\ckunpT}{C^{\tiny{\textit{Kun}}}_{p,T}}
\newcommand{\ckunp}{C^{\tiny{\textit{Kun}}}_{p}}
\newcommand{\ie}{\textit{i}.\textit{e}.\xspace }
\providecommand{\keywords}[1]{\textbf{\textit{Keywords---}} #1}

\begin{document}
\title{On the $\varepsilon$--Euler--Maruyama scheme for time-inhomogeneous jump-driven SDEs}
\author[1]{Mireille Bossy\footnote{mireille.bossy@inria.fr}}
\author[1]{Paul Maurer\footnote{paul.maurer@inria.fr}}
\affil[1]{Universit{\'e} C{\^o}te d'Azur, Inria, CNRS, France}
\date{August, 2025}
\maketitle

\begin{abstract}
    We consider a class of general SDEs with a jump integral term  driven by a time-inhomogeneous  Poisson random measure. We propose a two-parameters Euler-type scheme for this SDE class and prove an optimal rate for the strong convergence with respect to the $L^p(\Omega)$-norm and for the weak convergence, {considering integration over $n$ uniform time-steps.} 
    One of the primary issues to address in this context is the approximation of the noise structure when it can no longer be expressed as the increment of random variables. We extend the Asmussen--Rosiński approach to the case of a fully dependent jump coefficient and time-dependent Poisson compensation, handling contribution of jumps smaller than $\varepsilon$ with an appropriate Gaussian substitute and exact simulation for the large  jumps contribution. 
    For any $p \geq 2$, under hypotheses required to control the $L^p$-moments of the process, we obtain a strong convergence rate of order $1/p$. Under standard regularity hypotheses on the coefficients, we obtain a weak convergence rate of $1/n+\varepsilon^{3-\beta}$, where $\beta$ is the Blumenthal--Getoor index of the underlying Lévy measure. We compare this scheme with the Rubenthaler's approach where the jumps smaller than $\varepsilon$ are neglected, providing strong and weak rates of convergence in that case too. The theoretical rates are confirmed by numerical experiments afterwards. We apply this model class for some anomalous diffusion model related to the dynamics of rigid fibres in turbulence.\bigskip\newline
\keywords{stochastic differential equations with jumps; time-inhomogeneous {P}oisson random measures;  {E}uler--{M}aruyama scheme; strong rate of convergence; weak rate of convergence; anomalous diffusion in turbulence}
\end{abstract}

{\tableofcontents}

\section{Introduction}

We consider general stochastic differential equations (SDEs) with jumps, of the form:
\begin{align}\label{eq:intro_SDEs}
  X_t & = X_0+ \int_0^t a (s, X_s) d s + \int_0^t b (s, X_s) d W_s + \int_0^t
  \int_{- \infty}^{+ \infty} c ({s}, X_{s^-}, z)  \ \big(N (d s, d z) - \nu_s (d z) \ ds\big),
\end{align}
where $(W_t)_{t\geq 0}$ is a standard Brownian motion and $N (dt ,d z)$ is a  {Poisson random} measure with associated  time-dependent compensator $\nu_t (d z) dt$. In particular, we are interested into the construction of a simulation algorithm for such SDEs, as well as providing $L^p$-strong convergence and weak convergence rates for this numerical scheme.

In this paper, we focus on the case where the family of compensator measures has infinite activity, for all $t$ in $[0,T]$ (or in any identified time-subinterval that we assume extended to $[0,T]$ for simplicity): 
\begin{align*}
    \int_{\mathbb{R}} \nu_t(dz) = +\infty.
\end{align*}
Other cases lead to SDEs driven by compound Poisson process, that are nonetheless covered by our  convergence results.  

Numerical schemes for diffusion processes such as 
\begin{align}\label{eq:brownian_SDEs}
  X_t & = X_0+ \int_0^t a (s, X_s) d s + \int_0^t b (s, X_s) d W_s,
\end{align} 
has been extensively studied. In that case, the simulation of the noise increments $(W_{t_i+1} - W_{t_i}, t_i, i=0,\ldots, n)$ does not lead to any difficulty other than that of performance. 
The literature on the numerical analysis of approximation  schemes for \eqref{eq:brownian_SDEs} therefore develops  and analyses convergence of Euler-type methods under ever-weaker hypothesis on the smoothness and boundedness  of the coefficients.
This is in stark contrast to the situation we are studying, since jump SDEs can lead to a wide variety of noise structures.  
We refer to Section \ref{sec:review} for a detailed review on existing approximation results for Lévy-driven SDEs. This covers from the first extension of \eqref{eq:brownian_SDEs} -- when the Brownian motion  $(W_t)_{t\geq 0}$ is replaced by a Lévy process $(L_t)_{t\geq 0}$, which corresponds to a coefficient in the jumps integral term of the form $c(s,X_{s^-},z) = \sigma(s,X_{s^-}) z$ and $\nu_s(dz) = \nu(dz)$ in \eqref{eq:intro_SDEs} -- to SDEs driven by a {Poisson random} measure $N(ds,dz)$, when the jump coefficient $c(s,X_{s^-},z)$ has a direct impact on the noise structure, since it modifies the size $z$ of the jumps. In this situation, the notion of increment itself is meaningless, and one need to rely on other techniques to construct a numerical scheme.

In this paper, we  restrict our SDE setting to Lipschitz-in-space coefficients, to focus on the construction and the convergence analysis for a simulation algorithm, approximating jointly the noise structure and the equation. The regularity of the jump part with respect to the time, however, might be quite low without significantly affecting the convergence rate. As a result, we can simply rely on an time-integrability condition for $c$ and $\nu_t$.

\subsection{Motivations from modelling} 

This study is mainly motivated by the simulation of stochastic models, with a focus on investigating Lévy processes, and particularly $\alpha$-stable processes, arising as the limit distribution of the generalised Central Limit Theorem for independent random variables with infinite variance. One notable example of such random variables is given, in statistical physics for the description of turbulence in fluids,  by the probability distribution of statistically independent and fractal v{o}rtices in 2D turbulent flow, as explored in \cite{farge}. 

\paragraph{Modelling with finite variance dynamics.  } Nevertheless, even in such singular phenomena, natural cutoffs may enable the presence of a second moment, albeit potentially very high. For such reason, a wide class of models  introduce tempered stable processes, as discussed in \cite{dubrulle}. The Lévy measure for such processes is given by $\nu(dz) = \varphi(z) / |z|^{1+\alpha} {dz}$, where $\varphi(z)$ represents a cutoff function for the jump size. A simple example arises when $\varphi(z)$ is a truncation function (i.e., $\varphi(z) = \indic{{|z|} \leq B}$), resulting in a truncated stable process, where the jumps are spectrally truncated by an upper bound $B \geq 0$. This construction allows for the finite moments of all orders while maintaining power-law tails in the probability distribution function at a fixed time. Detailed mathematical analysis of tempered stable processes has been conducted by Rosiński \cite{rosinski}, who notably demonstrates their long time convergence to Gaussian processes under the Skorohod topology. 

From a numerical point of view, it should be acknowledged that an exact simulation algorithm for the increment of some classes of tempered stable process has been developed by Dassios, Lim and Qu in \cite{dassios} in the finite variation case (i.e., $\alpha \in (0,1))$. However, as of now, there is no existing method capable of exactly simulating the increments of tempered stable processes with infinite variation. 

\paragraph{From diffusive to anomalous regimes. }Lévy processes with finite variance evolve in a diffusive regime, which means that their variance grows linearly with time. This characteristic fails to capture certain physical phenomena exhibiting diverse variance regimes, such as anomalous diffusion (i.e., stochastic processes $X$ satisfying $\mathbb{E}[X_t^2] \simeq t^\beta$ with $\beta \neq 1$). 
In statistical physics, discrete non-Markovian stochastic processes are often employed as models for such phenomena (see, e.g., \cite{klafter}), under the name of Lévy Walks. Nonetheless, the rigorous mathematical analysis of Lévy Walks poses challenges, and their calibration remains difficult compared to continuous Markovian models. An intriguing alternative arises in the form of additive Markov processes, which are stochastic processes with independent increments that may not be stationary. These processes share similarities with Lévy processes, with the distinction lying in their time-dependent Lévy measure $\nu_t(dz)$.
A first motivation to introduce and simulate SDEs driven by this class of processes emerged from a prior study conducted by co-authors. In \cite{campana}, Campana, Bossy  and Bec  introduce and analyse a stochastic model  for the dynamics of  the orientation of small rods in  a 2D turbulent flow. The model is formulated as a standard Brownian SDE, which is shown to be suitable at the long-time equilibrium regime. However, for more complex statistics like rods tumbling variance,   this model significantly underestimates the impact of a transient regime and the contribution of a sharp local activity before reaching the long-time behaviour. We refer to \cite{campana} and to Section \ref{sec:application} for more details. 

\paragraph{Additive models in ecology and finance. } Stochastic models driven by general additive processes are also useful in  ecology as movement models. For instance, Burte, Cointe, {Perez, Mailleret and Calcagno} \cite{cointe} analyse  the population dynamics of parasitic insects. Depending on the initial population density, the mean squared displacement (i.e., the empirical variance) of insect positions exhibits various regimes, displaying either superdiffusive or subdiffusive behaviour. Furthermore, under specific conditions, the distribution of the individuals is not Gaussian and demonstrate heavy tails that are characteristic of stable-like processes.

Finally, since they offer more flexibility by allowing time-inhomogeneity, additive processes are able to capture financial option prices over a range of different maturities and strikes, while also providing a better fit to the observed marginal distributions of returns in financial asset empirical time series.  For instance, in \cite[{Part IV}]{tankov}, Tankov and Cont introduce an additive exponential model for risk-neutral option pricing, and conduct a comparison with some local volatility models.

\subsection{Approximation of Lévy-driven SDEs}\label{sec:review}
Over the past two decades, there has been a notable surge of interest in the simulation of SDEs driven by Lévy processes, leading to the exploration of various techniques.  The Brownian motion $(W_t)_{t\geq0}$ replaced by a
Lévy process $(L_t)_{t\geq0}$, still allows the construction of an Euler scheme, and in this case, $L^2$-strong convergence results has been proven by Kohatsu-Higa and  Protter \cite{kohatsu_euler}, while Talay and Protter  \cite{talay_protter} provided a bound for the weak convergence.  
However,  there is no systematic technique for the exact simulation of the law of the Lévy increment $L_{t_{i+1}} - L_{t_i}$, except for some specif cases, which arises the need to develop approximations method for the simulation of Lévy processes, and more generally for stochastic processes driven by Poisson random measures.

\paragraph{A first approach,  neglecting the small jumps. }
When the driving Lévy process has finite activity, meaning that it undergoes only a finite number of jumps during any given finite time interval, simulation becomes straightforward using the exact compound Poisson process representation given by $L_t = \sum_{i=1}^{N_t} Z_i$, where $(N_t)_{t\geq0}$ is a Poisson process and the $Z_i$ are independent and identically distributed random variables (assuming the existence of an algorithm to generate $Z_1$). However, challenges arise when dealing with Lévy processes possessing infinite activity, as they can entail an infinite number of arbitrarily small jumps within a finite time span. To tackle this obstacle, one approach involves neglecting jumps smaller than a fixed threshold $\varepsilon > 0$, thereby simulating only the finite activity portion of the process. In \cite{rubenthaler}, Rubenthaler conducted an analysis of a numerical scheme following this technique, providing  asymptotic behaviour of a normalised error process $\rho_n (\overline{X}^n - X)$, with a sequence $(\rho_n)_{n\geq0}$ as an appropriate rate going to infinity.  Another approach proposed by Rubenthaler and Wiktorsson \cite{rubenthaler2003improved} involves a series representation of the jumps of Lévy processes, offering a better rate of convergence, albeit applicable in specific cases only.

\paragraph{The Asmussen--Rosiński approximation. }
Asmussen and Rosiński  \cite{asmussen} introduced a method to enhance the precision of approximations in handling jumps smaller than $\varepsilon$ in Lévy-driven processes. Instead of discarding these small jumps, they proposed replacing them with an appropriate Brownian motion and employed a Berry--Essen technique to estimate the error. This method, widely explored in the literature, has become known as the Asmussen--Rosiński approximation. In \cite{rio},  Rio presented a Berry--Essen theorem  which offers an estimation of the convergence speed of the Central Limit Theorem in terms of transport distances. Specifically, Rio demonstrated that if $\mu_n$ represents the law of the renormalised sum of independent centred random variables $X_1,\ldots,X_n$, the Wasserstein distance $\mathcal{W}_p(\mu_n,\nu)$ between $\mu_n$ and the standard normal distribution $\nu$ can be bounded based on the Lyapunov coefficients $L_{p+2,n} = \sum_{i=1}^{n} \mathbb{E}[|X_i|^{p+2}]$, where $p\in(1,2]$. In  \cite{fournier},  Fournier showed that the usual $L^2$-strong convergence rate (i.e., $\mathbb{E}[\sup_{0 \leq i \leq n} |\overline{X}^n_{t_i} - X_{t_i}|^2]^\frac{1}{2} \leq C n^{-1/2}$) holds for the discretisation of the SDE $d X_t = \sigma(X_t) d L_t$, using a numerical scheme based on the Asmussen--Rosiński approximation. Fournier utilised Rio's theorem with $p=2$ to estimate the discrepancy between the small jumps and their approximations with respect to the $L^2$-norm. In \cite{dereich}, for the  convergence analysis of the same class of SDEs, Dereich 
uses  the Komlós--Major--Tusnády coupling (more precisely its generalisation by Zaitsev)  between Lévy  process, having its Lévy measure supported in a small ball $B(0,\varepsilon)$, and  a Wiener process (see the references therein). This technique can be used directly in dimension $d>1$. However, it seems to overestimate the speed of convergence by a log factor in contrast with Fournier's approach. 
 In the same vein, we would like to mention the weak error bound  derived in \cite{jacod_meleard}, where the authors demonstrated that $|\mathbb{E}[f(\overline{X}^n_T)] - \mathbb{E}[f(X_T)]| \leq C n^{-1}$, for a {\it modified} Euler scheme,  recovering the results of  Protter and Talay for the genuine (but not always constructive) Euler one. In this case,  the modification assumes some apriori weak convergence hypothesis on the approximation of the increment of the underlying Lévy process.

Although the Asmussen--Rosiński technique might not be optimal in some cases when it comes to weak approximation (see \cite{kohatsu2014optimal}, Remark 24), it has still been employed in several types of schemes. Notably, it has been employed in jump-adapted discretisations, such as those discussed in \cite{kohatsu_tankov} and more recently in \cite{tretyakov}, with both papers providing weak error bounds. One may also note the recent work of Bally and Qin \cite{bally}, who extended this method to an SDE closely related to \eqref{eq:intro_SDEs}, using a space-time Gaussian measure to approximate the small jump integral, and provide a bound for the total variation distance.

\medskip

In this paper, we have chosen to perform an extension of the Asmussen--Rosiński approach for the SDE \eqref{eq:intro_SDEs}. The first reason of this choice is due to the flexibility of the Gaussian approximation, allowing to deal with a complex jump coefficient. This enable{s} us to provide a simulation tool for a very wide range of stochastic models. A second reason lies in the existence of appropriate tools to conduct a strong convergence analysis. In fact, the theorem of Rio mentioned above has been extended by Bobkov \cite{bobkov} to all Wasserstein distance $\mathcal{W}_p$ with $p \geq 1$, which provide a powerful tool to obtain estimates with respect to any $L^p$-norm.  Concerning weak error analysis, we take advantage of the explicit scheme construction to quantify a two parameters rate of convergence,   involving the size of the mesh and the size of the small jump cutoff $\varepsilon$.

\subsection{Plan of the paper and  main techniques used}
In the present work, we aim to study the convergence rate of a numerical scheme for \eqref{eq:intro_SDEs}, where the jumps of the system smaller than some cutoff value $\varepsilon$ are either ignored or  approximated by appropriate Gaussian variables. Our analysis,  in $L^p$-norm and weak convergence, leads to optimisation choices between $\varepsilon$ and the number $n$ of {uniform} time steps. 

In section \ref{sec:hypo_prelim}, we present a general technique to represent and simulate stochastic integrals with respect to a {Poisson random} measure, which allows us to introduce the $\varepsilon$-Euler-Maruyama scheme $\oeX$, for which we discuss about its practical implementation. We also mention a version $\ueX$ of the scheme where the small jumps are only neglected.

In Section \ref{sec:strong_convergence}, we prove the strong convergence of this scheme in the $L^p$-norm by leveraging the Euler-Peano scheme (i.e., SDE \eqref{eq:intro_SDEs} with frozen coefficients) as a pivot term. The convergence of the Euler-Peano scheme is derived through a Gronwall argument. The time-inhomogeneity of the Lévy measure introduce new terms compared to the standard case, that can be estimated as long a sufficient time-integrability property is satisfied for the jump term. Subsequently, we compare the Euler-Peano scheme and the $\varepsilon$-Euler-Maruyama scheme in terms of the $L^p$-norm. This comparison is made possible by an estimation for the transport distance between the small jumps and their Gaussian approximation, which we obtain by the means of Bobkov's result \cite{bobkov}. In order to effectively utilise this result into the error analysis, we use a parametric optimal transport technique to construct a suitable version of the scheme, measurable with respect to the filtration of the process and satisfying a $\mathcal{W}^p$-optimal coupling property. 

In Section \ref{sec:weak_convergence}, we prove a weak convergence rate based on flow regularity of the solution process that ha{s} been proved recently in \cite{breton}.  

In Section \ref{sec:numerics}, we present numerical experiments, demonstrating the optimality of our theoretical bounds. 

Finally, we introduce in Section \ref{sec:application} an application in physics, by the proposition of a stochastic additive model to describe the anomalous diffusion tumbling of small rods in a 2D turbulent flow. 

\subsection{Summary of main results}
We refer to the corresponding theorems for detailed hypotheses on the coefficients of the SDE \eqref{eq:intro_SDEs} (in particular to simplify the discussion here we assume at least an 1/2-H{ö}lder regularity in time of the coefficient). 

The proposed schemes (with or without compensation for small jumps) are two-approximation parameters schemes, with parameters $n$ and $\varepsilon$.    We set out to determine the behaviour of the error as a function of $n$ and $\varepsilon$. 
We consider first strong $L^p(\Omega)$-error, with $p\geq 2$ in Theorem \ref{theoreme_maurer}  and \ref{thm:noapprox}, obtaining convergence by linking optimally $\varepsilon$ as a function of $n$. We get  the following convergence rates (see Corollaries \ref{cor:mainresult} and \ref{cor:scheme_cut_BG_index}): for a given number $n$ of homogeneous discretisation steps $(t_i)_{0 \leq i \leq n}$ of the interval $[0,T]$ and the optimal choice $\varepsilon = n^{-\left(\frac{1}{2}+\frac{1}{p} \right)}$, for the compensated scheme $\oeX$ 
\begin{align}\label{eq:cv_intro}
\Big\| \sup_{0 \leq i \leq n} \left| X_{t_i} - \oeX_{t_i} \right| \Big\|_{L^p(\Omega)} \leq C n^{ -\frac{1}{p} },
\end{align}
where  $C$ is a constant independent of $n$.

For the non-compensated scheme $\ueX$, the optimal choice is  $\varepsilon = n^{-\left(\frac{1}{p} \frac{2}{2 - \beta} \right)}$, to achieve the same rate 
\begin{align}
\Big\| \sup_{0 \leq i \leq n} \left| X_{t_i} - \ueX_{t_i} \right| \Big\|_{L^p(\Omega)} \leq C n^{ -\frac{1}{p} },
\end{align}
where $\beta$ is the Blumenthal--Getoor index of the compensator $(\nu_t)_{t \in [0,T]}$ (we refer to Corollary \ref{cor:scheme_cut_BG_index}  for a precise definition), and $C$ is a constant independent of $n$. 
For the typically standard RMS error (i.e. $p=2$), this means that  the compensated scheme performs better when the Blumenthal--Getoor index is greater than $1$, as it was observed  in \cite{dereich} for SDEs driven by Lévy noise.  It is, moreover, for these $\beta$ values greater than one that the numerical simulation of the jump integral becomes particularly cumbersome. It is therefore very cost-effective to work with bigger cutoff $\varepsilon$ while keeping the same precision in this case. In our analysis, note that the Blumenthal--Getoor index  threshold giving the advantage to the compensate scheme  is increasing with $p$ (for $p=4$, the threshold is $\beta=4/3$). We refer to the end of Section \ref{sec:main_strong} for a detailed behaviour of this threshold. 
\medskip

The situation is radically different when we consider weak convergence analysis, that gives the advantage to the compensated scheme for any $\beta$  with, for a smooth test function $\Phi$, 
\begin{equation*}
|\EE[\Phi(X_T)] - {\EE}[\Phi(\oeX_T)]| \leq \frac{C}{n} + C\varepsilon^{3 - \beta}, 
\end{equation*}
whereas
\begin{equation*}
|\EE[\Phi(X_T)] - {\EE}[\Phi(\ueX_T)]| \leq \frac{C}{n} +  C\varepsilon^{2 - \beta}. 
\end{equation*}
As detailed in the next section, the compensated scheme could require some extra computation (in worse cases, typically some numerical integration at each time step) but this price is more than offset by the gain in $\varepsilon$ while simulating the jumps stochastic integral, as  illustrated in Section \ref{sec:numerics}.

\subsection*{Notations  }

\noindent$\small{\bullet}$ We note $\mathbb{N}^\ast=\{1,2,\ldots\}$ the set of all positive integers and $\mathbb{R}_{+}= [0,+\infty)$ the set of positive or null real numbers.

\noindent$\small{\bullet}$ For a random variable $X$ on a probability space $(\Omega, \mathcal{F}, \mathbb{P})$, we note $\law{X}$ the law of $X$, and $\Var{(X)}$ be the variance of $X$.\smallskip

\noindent$\small{\bullet}$ We denote by $\mathcal{B}(\RR)$ the Borel $\sigma$-algebra on $\RR$.

\noindent$\small{\bullet}$ The notation $\mathcal{N}(m, \sigma^2)$ represents the normal distribution with a mean of $m$ and variance of $\sigma^2$.\smallskip

\noindent$\small{\bullet}$ For a real $T>0$, let $n\in\mathbb{N}^\ast$ be the number of discretisation steps. We consider {a} homogeneous  grid on $[0,T]$, namely $t_i = i\frac{T}{n}$ for $i\in\{0,...,n\}$. We denote the {càdlàg version of} the freezing time on the grid $t\mapsto \eta(t)= \lfloor \frac{t n}{T}\rfloor \frac{T}{n}$, and the corresponding index $t\mapsto \rho(t)= \lfloor \frac{t n}{T}\rfloor$.  \smallskip

\noindent$\small{\bullet}$ We denote by $B(\varepsilon)=\{ z \in \mathbb{R} ;  |z| < \varepsilon \}$ the open ball centred in zero with radius $\varepsilon>0$.

\noindent$\small{\bullet}$ For two real numbers $x$ and $y$, we write
\begin{equation*}
x \leqpt y, 
\end{equation*}
if there exists a constant $C>0$ that does not depends on $n$ and $\varepsilon$, such that 
\begin{equation*}
 x \leq C y.
\end{equation*}
This notation will be essentially used in the results statements and in some proofs to lighten the presentation, although it is always possible to specify an upper-bound for the constant $C$.  
\smallskip

\noindent$\small{\bullet}$  We will make use of the $L^q$-Wasserstein distance between two distributions $\mathcal{L}_1$ and $\mathcal{L}_2$, for $q \geq 1$,  defined as 
\begin{align*}
    \mathcal{W}_q(\mathcal{L}_1,\mathcal{L}_2) \coloneqq \inf_{(X_1,X_2) \in \pi(\mathcal{L}_1,\mathcal{L}_2)} \mathbb{E}[|X_1-X_2|^q]^{1/q},
\end{align*}
where $(X_1,X_2) \in \pi(\mathcal{L}_1,\mathcal{L}_2)$, are two random variables $X_1,X_2$ on $( \Omega,\mathcal{F},\mathbb{P})$ verifying $\law{X_i} = \mathcal{L}_i$.

\section{Preliminaries and schemes}\label{sec:hypo_prelim}

Let $T > 0$. We consider a probability space $(\Omega, \mathcal{F}, \mathbb{P})$, and a $\sigma$-finite kernel $\nu_{\cdot}(\cdot):[0,T] \times \mathcal{B}(\RR) \rightarrow [0,+\infty]$, satisfying 
\begin{align*}
    \int_{\mathbb{R} \backslash B(\varepsilon)} \nu_t(dz) < +\infty,  \quad \text{for all $t\in[0,T]$ and $\varepsilon >0$.}
\end{align*}
We consider a {Poisson random} measure $(N ({[0,t]} ,A), \ t \in [0,T], \ A \in \mathcal{B} (\mathbb{R}))$ on $(\Omega, \mathcal{F}, \mathbb{P})$ whose intensity is the $\sigma$-finite measure $\textit{Leb}_{[0,T]} \otimes \nu$, where $\textit{Leb}_{[0,T]}$ designates the Lebesgue measure on the interval $[0,T]$. We refer to \cite{last2017lectures} for the precise definition of a $\sigma$-finite kernel and the detailed construction of such a {Poisson random}  measure. 
We assume that the probability space $(\Omega, \mathcal{F}, \mathbb{P})$ is rich enough to contain a  standard $\mathbb{R}$-valued Brownian motion $(W_t)_{t\in[0,T]}$ independent from $N$.

The compensated Poisson measure of $N$ will be denoted as $\rpm \coloneqq N(ds,dz) - \nu_s(dz) ds$.
We also introduce the filtration $(\mathcal{F}_t)_{t \in [0,T]}$ generated by $W$ and $N$ on $[0,T]$.

\subsection{{$L^p$}-wellposedness and moments bound for the solution of \eqref{eq:intro_SDEs} }

The existence in $L^p$ and the uniqueness of strong solutions for \eqref{eq:intro_SDEs} with Lipschitz coefficients can be obtained by standard arguments, such as in Theorem 3.1 in Kunita \cite{kunita}.  

In the class of SDEs \eqref{eq:intro_SDEs}, we paid attention to the time-behaviour of the coefficients, in particular the jump coefficient $c$.  For that reason, we choose to work with a framework similar to Breton and Privault \cite{breton}, that requires weaker integrability hypothesis than in \cite{kunita}. 

Let  $a:[0, T] \times \mathbb{R} \rightarrow \mathbb{R}$, $b:[0,
T] \times \mathbb{R} \rightarrow \mathbb{R}$ and $c:[0, T] \times \mathbb{R}
\times \mathbb{R} \rightarrow \mathbb{R}$ be measurable deterministic
coefficients. 

\begin{hyp}{\hspace{-0.15cm}\Blue{\bf  \textrm Lipschitz condition. (\ref{hyp:lipschitz}).}}
\makeatletter\def\@currentlabel{ {\bf\textrm H}$_{\mbox{\scriptsize\bf\textrm{Moment$\bm{(p)}$}}}$}\makeatother
\label{hyp:lipschitz}
The compensator measure family $(\nu_t)_{t\in[0,T]}$ satisfies for all $t\in [0,T]$, $\int_{\mathbb{R} \backslash B(\varepsilon)} \nu_t(dz) < + \infty$,  {for all $\varepsilon >0$.} 
There exists a measurable function $L_{a,b}: [0,T] \rightarrow \mathbb{R}_{+}$, and a measurable function $L_c:[0,T] \times \mathbb{R} \rightarrow\mathbb{R}_{+}$ such that, for all $(x,y)$ 
in $\RR^2$, 
\begin{equation*}
\begin{aligned}
| a (t, x) - a (t, y)|  + | b (t, x) - b (t, y) | & \leq L_{a,b}(t) | x - y |, & \quad t \in [0,T],\\
| c (t, x, z) - c (t, y, z) | &  \leq L_c (t,z)  \  | x - y |,  & \quad  t \in [0,T], \ z\in \mathbb{R}.
\end{aligned}
\end{equation*}
We define
\begin{align*}
\overline{L}_c(t,z)  = L_c(t,z) \vee |c(t,0,z)|. 
\end{align*}
There exists $p\geq 2$ such that $X_0 \in L^p(\Omega)$, {$a(\cdot,0)$, $b(\cdot,0)$} and $L_{a,b}(\cdot)$ are in $L^p([0,T])$. Moreover, defining the map $t\in [0,T] \mapsto \psi_p(t)$ by 
\begin{equation} \label{eq:def_psi_p}
    \psi_p(t) = \left( \int_{-\infty}^{+\infty} \left| \overline{L}_c(t,z)  \right|^2 \nu_t(dz) \right) ^ {p/2} + \int_{-\infty}^{+\infty} \left|  \overline{L}_c(t,z)  \right|^p \nu_t(dz), 
\end{equation}
we assume that  $\psi_p \in L^1([0,T])$. 
\end{hyp}

Under \ref{hyp:lipschitz} with $p\geq 2$, according to Theorem 3.1 in \cite{breton}, there exists a unique strong solution $(X_t)_{t\in[0,T]}$ of \eqref{eq:intro_SDEs} in $L^p(\Omega)$ and, with $m (p, T)$ depending  on the coefficients $a,b$ and $c$,  we have the moment bound  
\begin{equation} \label{eq:moment_bound}
    \mathbb{E} \Big[ \underset{t \in [0, T]}{\sup} | X_t |^p \Big] =m
    (p, T) \  < \ +\infty.
\end{equation}

\subsection{Representation and simulation of a {stochastic} Poisson integral} \label{sub:simulation_poisson_integral}

Let $F:\mathbb{R}_+ \times \mathbb{R} \times \Omega \rightarrow \mathbb{R}$ be a{n} $(\mathcal{F}_t)$-predictable process such that $\int_0^T \int_{-\infty}^{+\infty} |F(s,z)|^2 \rpc < +\infty$  $\mathbb{P}$-a.s.

Using a threshold  $\varepsilon >0$, the stochastic {Poisson} integral 
\begin{equation*}
I(F) \coloneqq \int_0^T \int_{-\infty}^{+\infty} F({s},z) \rpm 
\end{equation*}
can be separated into its "large jumps" part $I_l^{\varepsilon}(F)$ and its "small jumps" part $I_s^{\varepsilon}(F)$, that is 
\begin{equation*}
    I_l^{\varepsilon}(F)  = \int_0^T \int_{\mathbb{R} \backslash B(\varepsilon)} F({s},z) \rpm,  \qquad \qquad
    I_s^{\varepsilon}(F)  = \int_0^T \int_{B(\varepsilon)} F({s},z) \rpm.
\end{equation*}

\paragraph{The small jumps integral.  } Generally,  exact simulation of  the "small jumps" integral $I_s^{\varepsilon}(F)$ is not possible, unless we are in the case of finite activity, (i.e. $\int_{\mathbb{R} } \nu_t(dz) < +\infty$).  As discussed in the introduction, we adopt the method of the Asmussen--Rosiński (AR) approximation.  This approach substitutes the stochastic integral $I_s^{\varepsilon}(F)$ with a Gaussian random variable possessing an equivalent variance.
\begin{equation*}
    \law{ I_s^{\varepsilon}(F)}  \simeq \left( \int_0^T \int_{B(\varepsilon)} \mathbb{E} [|F(s,z)|^2] \rpc \right)^{\frac{1}{2}}  \mathcal{N}(0,1).
\end{equation*}
While this approximation only requires to compute the second moment of the small jumps, it still allows to get satisfying error bounds in the $L^{p}$-Wasserstein distance for any $p \geq 1$. More precisely, {we establish in Section \ref{sub:proof_wasserstein}} the following Proposition \ref{prop:wasserstein} in the case where $F$ is  deterministic. 
For this purpose, we add the following assumption regarding the integrability condition for $F$. 

\begin{proposition} \label{prop:wasserstein}
 Let $F:\mathbb{R}_+ \times \mathbb{R} \rightarrow \mathbb{R}$ be a deterministic measurable map, and $q\geq1$. There exists $\varepsilon_F > 0$ such that $F(t,\cdot)$ is non identically zero on $B(\varepsilon_F)$ for any $t\in[0,T]$ and
\begin{align}\label{eq:hypo_propo:wasserstein}
  \int_0^T \int_{B(\varepsilon_F)} | F(s,z) |^{q+2} \nu_s (dz) d s  \ + \ \int_0^T \left( \int_{B(\varepsilon_F)} | F(s,z) |^2 \nu_s (dz) \right) ^{q/2+1} d s  \ < +\infty.
\end{align}
Then there exists a constant $\mathcal{A}(q)$,  only depending on $q$,  such that for every $0 < \varepsilon < \varepsilon_F$, the following inequality holds for any $t_0, t \in [0,T]$:
\begin{align*}
   &  \mathcal{W}_{q} \left( \law{\int_{t_0}^t \int_{B(\varepsilon)} F({s},z) \rpm}, \ \mathcal{N} \left( 0 , \int_{t_0}^t \int_{B(\varepsilon)} | F(s,z) |^2 \nu_s (d z) d s \right) \right)\\ 
    & \leq \mathcal{A}(q) \left(  \frac{\int_{t_0}^t \int_{B(\varepsilon)} |
    F(s,z) |^{q + 2} \nu_s (d z) d s}{\int_{t_0}^t \int_{B(\varepsilon)}  | F(s,z) |^2 \nu_s (d z) d s}\right) ^{1/q}.
\end{align*}
\end{proposition}

\paragraph{The large jumps integral.  } Represented as a finite random sum, the "large jumps" integral $I_l^{\varepsilon}(F)$ can be  simulated exactly in many application cases.  Indeed, separating the Poisson measure from the compensation, 
\begin{align}
    I_l^{\varepsilon}(F) & = \int_0^T \int_{\mathbb{R} \backslash B(\varepsilon)} F({s},z) N(ds,dz) \  - \int_0^T \int_{\mathbb{R} \backslash B(\varepsilon)} F(s,z) \rpc,
\end{align}
and since $\mathbb{R} \backslash B(\varepsilon)$ is bounded below in the sense of \cite{apple} (i.e its closure does not contain zero), one can rewrite 
\begin{equation} \label{randomsum}
    \int_0^T \int_{\mathbb{R} \backslash B(\varepsilon)} F({s},z) N(ds,dz) = \sum_{0\leq u\leq T} F(u,\Delta P(u)) \indic{\mathbb{R} \backslash B(\varepsilon)}(\Delta P(u)),
\end{equation}
where the process $(P(t),t\geq 0)$ is defined by $P(t) = \int_{\mathbb{R} \backslash B(\varepsilon)} z\, N({[0,t]},dz)$. 

Let $N^{\varepsilon}(t) = N([0,t],\mathbb{R} \backslash B(\varepsilon))$. Then $N^{\varepsilon}$ is a time-inhomogeneous Poisson process, with intensity $\lambda^{\varepsilon}(t) = \int_{\mathbb{R} \backslash B(\varepsilon)} \nu_t (dz)$, and jump times $T^{\varepsilon}(j) = \inf\{t \in [0,T], N^{\varepsilon}(t)= j \}$.

Finally, setting $Z^{\varepsilon}(n) = \int_{\mathbb{R} \backslash B(\varepsilon)} z \,N({[0,T^{\varepsilon}(n)]},dz) - \int_{\mathbb{R} \backslash B(\varepsilon)} z\, N({[0,T^{\varepsilon}(n-1)]},dz)$, a straightforward extension of Theorem 2.3.9 in \cite{apple} to the time-inhomogeneous case shows that $P$ is nothing but the (time-inhomogeneous) compound Poisson process $$P(t) = \sum_{j=1}^{N^{\varepsilon}(t)} Z^{\varepsilon}(j).$$
In particular, one has $\Delta P(u) = Z^{\varepsilon}(j)$ for $j$ such that $u = T^{\varepsilon}(j)$, and \eqref{randomsum} becomes 
\begin{equation}\label{eq:rep_large_jumps_int}
    \int_0^T \int_{\mathbb{R} \backslash B(\varepsilon)} F({s},z) N(ds,dz) = \sum_{j=1}^{N^{\varepsilon}(T)} F(T^{\varepsilon}(j),Z^{\varepsilon}(j)) \ \indic{ \mathbb{R} \backslash B(\varepsilon) }(Z^{\varepsilon}(j)).
\end{equation}

\subsection{The $\varepsilon$-Euler-Maruyama scheme and its implementation}\label{sec:the_EM_scheme}

Let $N^\varepsilon(\cdot)$ be the  Poisson process with intensity function $\lambda^\varepsilon(t) = \int_{\mathbb{R}\backslash B(\varepsilon)} \nu_t(dz)$  considered above. We introduce its   jump times $T^\varepsilon(j)$, and  set the corresponding jump sizes to be 
{
$$Z^{\varepsilon}(j) = \int_{\mathbb{R} \backslash B(\varepsilon)} z \ N([0,T^{\varepsilon}(j)],dz) - \int_{\mathbb{R} \backslash B(\varepsilon)} z \ N([0,T^{\varepsilon}(j-1)],dz),$$}
 with the notation used in \eqref{eq:rep_large_jumps_int}. 
\smallskip

For $n \in \mathbb{N}^{\ast}$, we define $0 = t_0 < \cdots < t_n = T$, a discretisation of the interval $[0, T]$ with {uniform} steps, i.e $t_i=i\frac{T}{n}$. Let $(\xi_i)_{1 \leq i \leq n}$ be a sequences of i.i.d standard Gaussian random variables, independent from $W$ and $N$. 
The $\varepsilon$-Euler-Maruyama scheme for the simulation of the solution of \eqref{eq:intro_SDEs} is then defined by $\oeX_{t_0} = X_0$ and for $1\leq i\leq n$,
\begin{equation}\label{eq:numerical_scheme}
    \begin{aligned} 
    \oeX_{t_i} = & \oeX_{t_{i - 1}} + \left( a (t_{i - 1}, \oeX_{t_{i -
   1}}) - \int_{t_{i-1}}^{t_i} \int_{\mathbb{R} \backslash B(\varepsilon)} c(s,\oeX_{t_{i - 1}},z) \rpc \right) \frac{T}{n} +   b (t_{i - 1}, \oeX_{t_{i - 1}}) \ (W_{t_i} - W_{t_{i-1}}) \\ & 
   + \left( \int_{t_{i-1}}^{t_i} \int_{B(\varepsilon)} c^2(s,\oeX_{t_{i - 1}},z) \rpc \right)^\frac{1}{2} \xi_{i} +  \sum_{j=N^\varepsilon(t_{i-1})+1}^{N^\varepsilon(t_{i})} c(T^\varepsilon(j),\oeX_{t_{i -
   1}},Z^\varepsilon(j)).
    \end{aligned}
\end{equation}

In contract with Euler schemes for diffusion, effectively implementing this scheme requires: 
\begin{enumerate}[label=(\Alph*)]
    \item to generate the Poisson process $N^\varepsilon(\cdot)$ on $[0,T]$; 
    \item to generate the jump sizes random variables $Z^\varepsilon(j)$ for any $j \leq N^\varepsilon(T)$;
    \item to compute the integrals $\int_{t_{i-1}}^{t_i} \int_{\mathbb{R} \backslash B(\varepsilon)} c(s,\oeX_{t_{i - 1}},z) \rpc$ and $\int_{t_{i-1}}^{t_i} \int_{B(\varepsilon)} c^2(s,\oeX_{t_{i - 1}},z) \rpc$ at each step $0 \leq i \leq n$. \label{point_c} 
\end{enumerate}

\paragraph{(A) } In the literature, two primary classical approaches are commonly used for simulating the Poisson process.
\begin{itemize}
    \item When the intensity function $\lambda^\varepsilon (\cdot)$ can be explicitly inverted, one may perform a time-scale transformation of an homogeneous Poisson process with rate 1 to generate $N^\varepsilon$, as it is described in \cite{cinlar} (pp. 96-97). This method will however be less efficient when the inversion of $\lambda^\varepsilon(\cdot)$ has to be done numerically. 
    
    \item If the intensity function verifies $\lambda^\varepsilon(t) \leq \lambda^\ast$ for any $t \in [0,T]$, then one can rely on the thinning method (see \cite{lewis}), which is based on a acceptance--rejection technique, and does not require to invert the intensity function. The thinning method will be less efficient if $\lambda^\ast$ is too large, leading to high rejection probability. Note also that thinning cannot be used in the case of singular compensator measure such as $\nu_t(dz) = \frac{1}{\sqrt{t}} \nu(dz)$, whose noise remains admissible in our class of SDE, as long as $t\mapsto c(t,x,z)$ is not too singular in its turn (see the hypothesis \ref{hyp:peano} below). 
\end{itemize}

\paragraph{(B) }  To generate the random variables $Z^{\varepsilon}(j)$ for $j \leq N^{\varepsilon}(T)$, one can derive from the work of Orsingher, Ricciuti and Toaldo  \cite{randomsums} (see part 2.2) that, {providing that $\nu_t$ is absolutely continuous with respect to the Lebesgue measure for every $t \in [0,T]$,} the conditional distribution of $Z^{\varepsilon}(j)$ given the jump times $T^{\varepsilon}(j)$ is {itself} absolutely continuous with respect to the Lebesgue measure, and satisfies
\begin{equation} \label{jump_size_law}
  \forall B \in \mathcal{B}(\mathbb{R}) \qquad \mathbb{P} \left( Z^{\varepsilon}(j) \in B \text{ } | \text{ } T^{\varepsilon}(j) = t \right)
  = \frac{\nu_t (B \cap \mathbb{R} \backslash B(\varepsilon))}{\nu_t(\mathbb{R} \backslash B(\varepsilon))}.
\end{equation}
In simple cases, this distribution can be simulated exactly by inversion or acceptance--rejection method. Note in particular that if the time-dependence of the compensator is multiplicative, i.e if one has $\nu_t(dz) = \phi(t) \nu(dz)$, then the distribution \eqref{jump_size_law} is homogeneous in time, and the jump sizes $Z^{\varepsilon}(j)$ become i.i.d random variables.

\begin{exmp}\label{example:nu}
    Let $\rpc = \phi(s) |z|^{-(1+\alpha)} \indic{|z| \leq b} {dz \, ds}$ be the compensator of a time-inhomogeneous $b$-truncated $\alpha$-stable process, with $\alpha \in (0,2)$ and $b > \varepsilon$.  The random variables $Z^{\varepsilon}(j)$ are thus i.i.d, with the same law as $Z^{\varepsilon}(1)$. We derive from \eqref{jump_size_law} that the cumulative density function $G^{\varepsilon}$ of $Z^{\varepsilon}(1)$ is explicitly given by 

    \begin{align*}
    G^{\varepsilon} (x) & =  \frac{1}{2}  \ 
    \left\{\begin{array}{rl}
    0, & \mbox{if } x < - b\\
     \dfrac{(- x)^{- \alpha} - b^{- \alpha}}{\varepsilon^{- \alpha}
    - b^{- \alpha}}, & \mbox{if }  x \in [- b, - \varepsilon]\\
    1, & \mbox{if }  x \in [- \varepsilon, \varepsilon]\\
    1 + \dfrac{\varepsilon^{- \alpha} - x^{-
    \alpha}}{\varepsilon^{- \alpha} - b^{- \alpha}},& \mbox{if }  x
    \in [\varepsilon, b]\\
    2, & \mbox{if } x \geq b.
    \end{array}\right. 
    \end{align*}
    This allows to compute the quantile function of $G^{\varepsilon}$, namely $Q^{\varepsilon} (y) = \inf \left\{ x \in \mathbb{R};~  {y \leq G^{\varepsilon} (x)} \right\}$:
    \begin{align*}
        \mbox{for all }  y \in] 0, 1], \quad Q_{\varepsilon} (y) & = 
  \left\{\begin{array}{rl}
    - \{ 2 y (\varepsilon^{- \alpha} - b^{- \alpha}) + b^{- \alpha} \}^{-
    \frac{1}{\alpha}}, & \mbox{if }  y \in (0, \frac{1}{2}],\\
    \{ (1 - 2 y)  (\varepsilon^{- \alpha} - b^{- \alpha}) + \varepsilon^{-
    \alpha} \}^{- \frac{1}{\alpha}}, & \mbox{if } y \in (\frac{1}{2}, 1].
  \end{array}\right.
    \end{align*}
Then, to simulate the random variable $Z^{\varepsilon}(1)$, one can rely on a straightforward inversion method: starting with a uniform variable $U \sim \mathcal{U} (] 0, 1])$, it is sufficient to compute $Q^{\varepsilon} (U) \sim Z^{\varepsilon}(1)$.
\end{exmp}

\medskip

\paragraph{(C) } Finally, the two integrals $\int_{t_{i-1}}^{t_i} \int_{\RR  \backslash B(\varepsilon)} c(s,\oeX_{t_{i - 1}},z) \rpc$ and $\int_{t_{i-1}}^{t_i} \int_{B(\varepsilon)} c^2(s,\oeX_{t_{i - 1}},z) \rpc$ can be pre-computed if the dependence of $c(\cdot,\cdot,\cdot)$ in the space variable $x$ is multiplicative, i.e if $c(s,x,z) = \sigma(x) w(s,z)$. In the general case, one may rely on a Taylor expansion of $c(s,x,z)$ when $z \rightarrow 0$ to compute the integral on $B(\varepsilon)$ when $\varepsilon$ is small enough, or use numerical integration techniques (see for example Appendix \ref{sec:details_arctan} for the case \eqref{eq:sde_atan_drift} in Section \ref{sec:numerics}).

\subsubsection*{The $\varepsilon$-Euler-Maruyama scheme without substitute}

As already mentioned in the introduction, the first approach introduced for SDEs driven by Lévy noise consists in fully neglecting the small jumps. This scheme can be seen as sub-case in our strong and weak convergence analysis.  In that case,  the $\varepsilon$-Euler-Maruyama scheme is denoted $\ueX$, {defined by $\ueX_{t_0} = X_0$ and for $1 \leq i \leq n$,} 
\begin{equation}\label{eq:numerical_scheme_without_s_jumps}
    \begin{aligned} 
    \ueX_{t_i} = & \ueX_{t_{i - 1}} + \left( a (t_{i - 1}, \ueX_{t_{i -
   1}}) - \int_{t_{i-1}}^{t_i} \int_{\mathbb{R} \backslash B(\varepsilon)} c(s,\ueX_{t_{i - 1}},z) \rpc \right) \frac{T}{n}  + b (t_{i - 1}, \ueX_{t_{i - 1}}) \ (W_{t_i} - W_{t_{i-1}}) \\
   & \quad +   \sum_{j=N^\varepsilon(t_{i-1})+1}^{N^\varepsilon(t_{i})} c(T^\varepsilon(j),\ueX_{t_{i -
   1}},Z^\varepsilon(j)).
    \end{aligned}
\end{equation}
 with the same notations introduced for $\oeX$  and the same simulation techniques. 

\section{Strong convergence}\label{sec:strong_convergence}

The convergence analysis  is divided into two parts. In order to get an upper bound for the difference between the small jumps of the SDE \eqref{eq:intro_SDEs} and the approximating Gaussian random variables, one needs a time-discretisation setting to apply Proposition \ref{prop:wasserstein}. Therefore, we introduce as a pivot term the {process} with frozen coefficients (or Euler-Peano scheme) $\widetilde{X}$ defined by
\begin{align} \label{eq:frozen_sde}
  \widetilde{X}_t & = {X_0 +} \int_0^t a (\eta (s^{{-}}), \widetilde{X}_{\eta (s^{{-}})}) d s +
  \int_0^t b (\eta(s^{{-}}), \widetilde{X}_{\eta (s^{{-}})}) d W_s + \int_0^t \int_{-
  \infty}^{+ \infty} c ({s}, \widetilde{X}_{\eta (s^-)}, z)  \rpm.
\end{align}

Firstly, we will demonstrate the strong convergence of the Euler-Peano scheme $\widetilde{X}$ in the $L^p$ norm.  For this purpose, we need to reinforce the $L^1$-integrability condition of $t\mapsto\psi_p(t)$, defined in \eqref{eq:def_psi_p}. In addition, we need a time-Hölder regularity condition for the drift coefficient $a$ and the diffusion coefficient $b$ of \eqref{eq:intro_SDEs}. This is not required for the jump coefficient $c$, since in this work, the jumps integral simulation is proposed  without freezing the time, as outlined in \ref{point_c}.
This leads us to introduce the following hypothesis.
\begin{hyp}{\hspace{-0.15cm}\Blue{\bf  \textrm Condition for Peano approximation. (\ref{hyp:peano}).}}
\makeatletter\def\@currentlabel{ {\bf\textrm H}$_{\mbox{\scriptsize\bf\textrm{Peano$\bm{(p)}$}}}$}\makeatother 
\label{hyp:peano} 
One has $L_{a,b}(\cdot) \in L^\infty([0,T])$, and there exists a real number $\gamma \in (0,1] $ such that
\begin{equation*}
\begin{aligned}
| a (t, x) - a (s, y)| + | b (t, x) - b (s, y) | & \leq \sup_{t\in [0,T]} L_{a,b}(t)  \ \left( |t - s|^{\gamma} + | x - y | \right), &\quad x,y \in \mathbb{R},& \; s,t \in [0,T].
\end{aligned}
\end{equation*}
We define
\[\overline{L}_{a,b}  = \sup_{t\in [0,T]} L_{a,b}(t) \vee \sup_{t\in[0,T]} \Big(|a(t,0)| + |b(t,0)|\Big).\]

In addition, there exists $\zeta \in (0,1]$ such that $\psi_p$, defined in \eqref{eq:def_psi_p} is in $L^{1+\zeta}([0,T])$. 
\end{hyp}

We provide a convergence rate for $|X-\widetilde{X}|$ in the following lemma. 
\begin{lemma}[\bf Strong convergence of the Euler-Peano scheme]\label{lem_peano}
  Let $p\geq 2$.  Assume that \ref{hyp:lipschitz} and \ref{hyp:peano} hold for that $p$. Then for all $n \in \mathbb{N}^{\ast}$,
  \begin{align} \label{euler_rate}
    \Big\| \underset{t \in [0, T]}{\sup} | X_t - \widetilde{X}_t
    | \Big\|_{L^{p}(\Omega)} & \leqpt  n^{-
    \left\{\gamma \wedge \left( \frac{2 \zeta }{p (1+\zeta)} \right) \right\}} 
  \end{align}
  where $X$ and $\widetilde{X}$ are the respective solution of \eqref{eq:intro_SDEs}  and \eqref{eq:frozen_sde}, and $\gamma, \zeta  \in (0,1]$ are  given in \ref{hyp:peano}.
\end{lemma}

\begin{remark}   
    Note that if the maps $s \mapsto \int_{-\infty}^{+\infty} |L_c(s,z) \vee c(s,0,z)|^2 \nu_s(dz)$ and $s \mapsto \int_{-\infty}^{+\infty} |L_c(s,z) \vee c(s,0,z)|^p \nu_s(dz)$ are in $L^{\infty}([0,T])$, then the choice $\zeta = 1$ in \ref{hyp:peano}, (\ie $\psi_p\in L^{2}([0,T])$), leads to the convergence rate ${n^{- (\gamma \wedge \frac{1}{p}) }}$. 
    Imposing such stronger condition would nonetheless exclude some relevant examples such as time-dependent truncated $\alpha$-stable-type compensator producing anomalous subdiffusion. For instance if $c(s,x,z) = z$ and $\nu_s(dz) = \frac{1}{\sqrt{s}} |z|^{-1-\alpha} \indic{|z|<1} {dz}$; taking $p=2$ we have in this case $\psi_2(t) = \tfrac{4}{(2-\alpha)}\tfrac{1}{ \sqrt{t}}$, which verifies \ref{hyp:peano} for any $\zeta \in (0,1)$ but not for $\zeta = 1$ since $\psi_2$ does not belong to $L^{2}([0,T])$. 
\end{remark}

The  proof of Lemma   \ref{lem_peano} (in Section \ref{sec:peano}) follows a standard approach, utilising Kunita's inequality instead of Burk{ho}lder--Davis--Gundy inequality and taking advantage of the Lipschitz conditions imposed on the coefficients. Nevertheless, we have taken special care to address the additional techniques required by the minimal regularity assumptions introduced by time-inhomogeneity. Specifically, the integrability condition in \ref{hyp:peano} with respect to time for the jump coefficient $c(s,\cdot,\cdot)$ and the compensator $\nu_s(dz)$ presents new difficulties when dealing with the local error of the Euler-Peano scheme, which we handle  {using an extended Gronwall Lemma}.

{
\begin{remark}[On the optimality of the $L^p$-convergence rate in Lemma \ref{lem_peano}] 
In the most standard case where $\zeta = 1$ and $\gamma \geq \tfrac{1}{p}$, the obtained  rate of convergence for the Euler-Peano scheme is $\tfrac{1}{p}$. Note that in the specific case of SDE driven by additive Lévy noise (with additional assumptions on the characteristic function that include truncated $\alpha$-stable processes), a  study by Butkovsky, Dareiotis and Gerencsér \cite{dareiotis} establishes  a $L^p$-convergence rate of order $1/2^{-}$ for any value of $p$. However, an application of the Itô formula shows that the local error $\| X_t - X_{\eta(t)} \|_{L^p(\Omega)}$ cannot behave better than $(t-\eta(t))^{1/p}$ when the driven process has symmetric or positive Lévy measure. 

In the subordinator case, the positivity of the driven process allows to decouple the local error from the frozen error $\| X_{\eta(t)} - \wtX_{\eta(t)} \|_{L^p(\Omega)}$ in the analysis. The following Lemma (proved in  Appendix \ref{sec:proof_lemma_rate}) shows that one cannot hope for a better rate of convergence than $\tfrac{1}{p}$ for multiplicative SDEs driven by a subordinator. Extending this theoretical statement to the case where the driven  process has a symmetric Lévy measure appears challenging. However, our numerical simulations on the SDE $dX_t = \cos(X_t) dt + \int_{\mathbb{R}} \sin(X_{t^{-}}) d L_t$ with $L$ being a truncated $\alpha$-stable Lévy process,  suggest a positive answer to this question (see Section \ref{sec:num_for_rem_optimal}). 
\end{remark}

\begin{lemma}[Example of optimal $\frac{1}{p}$-convergence rate for the Euler-Peano scheme]  \label{lm:optimality_of_rate}
Consider the strong solution $X$ of the multiplicative subordinator-driven  SDE
\begin{equation*}
X_t = 1 + \int_0^t X_{s^-} d L_s, \quad t \in [0, T],
\end{equation*}
where $L$ is a truncated subordinator with characteristic function $\EE[e^{i x L_t}] = \exp(t \int_0^1 (e^{i z x}-1) \nu(dz))$,  its Lévy measure $\nu$ is supported in $(0,+\infty)$ and satisfies $\int_0^1 z \nu(dz) < +\infty $. Then for any integer $p \geq 2$,  there exists a constant $C_{p, T} > 0$ such that
\begin{equation*}
   \Big\| \underset{t \in [0, T]}{\sup} |X_t - \wtX_t| \Big\|_{L^{p} (\Omega)} \geq C_{p, T} n^{- \frac{1}{p}}.
\end{equation*}
\end{lemma}
}

\subsection{Main result}\label{sec:main_strong}

Let us discuss first the case where the family of compensator measures has infinite activity, that is that for any $t\in [0,T]$: 
\begin{align}\label{eq_prerequit_AR}
    \int_{\mathbb{R}} \nu_t(dz) = +\infty.
\end{align}
Unless we neglect the influence of small jumps, standard substitution methods are using  thoughtfully chosen sequence of Gaussian distributions to define the scheme $\oeX$ in \eqref{eq:numerical_scheme}.   As a consequence, the scheme is  employing a substitute for the small jumps that is inherently decoupled from the noise of \eqref{eq:intro_SDEs}.
For this reason, in order to provide a strong convergence rate for the numerical scheme $\oeX$, one needs to construct a version $\oeeX$ of $\oeX$ where the approximation of small jumps by Gaussians lives on the same probability space as the solution $X$ to the SDE \eqref{eq:intro_SDEs}. 
This approximation can be quantified in term of the $L^p(\Omega)$-norm using Proposition~\ref{prop:wasserstein}, assuming the appropriate moments condition. However, this step requires to identify an optimal coupling with the actual small jumps  with respect to the $\mathcal{W}_p$ distance. In the class of SDE we are examining, where the jump coefficient depends on both the state and the jumps, this operation is more intricate compared to situations where state dependence is purely multiplicative. In particular, we need the existence of an optimal transport map between the distribution of the small jumps $\mathcal{L}_1$ at fixed discretisation step and its Gaussian approximation $\mathcal{L}_2$. According to \cite{fournier_giet}, this condition is satisfied as far as $\mathcal{L}_1$ is absolutely continuous with respect to the Lebesgue measure (see \cite{fontbona} and \cite{villani}). This leads us to introduce the following hypothesis.

\begin{hyp}{\hspace{-0.15cm}\Blue{\bf  \textrm Assmussen-Rosiński approximation for small jumps (\ref{hyp:ar}).}}
\makeatletter\def\@currentlabel{ {\bf\textrm H}$_{\mbox{\scriptsize\bf\textrm{AR\_method}}}$}\makeatother 
\label{hyp:ar}
Assume that the time-compensator measures satisfy \eqref{eq_prerequit_AR}. In addition, assume that there exists $\varepsilon^*>0$ such that:
\begin{itemize}
    \item {\bf (Moments) } for every $x\in \mathbb{R}$ and $t\in[0,T]$,  $c(t,x,\cdot)$ is non identically zero on $B(\varepsilon^*)$ and 
\begin{align*}
  \int_0^T \int_{B(\varepsilon^*)} | c(t,x,z) |^{p+2} \nu_t (dz) d t  \ + \ \int_0^T \left( \int_{B(\varepsilon^*)} | c(t,x,z) |^2 \nu_t (dz) \right) ^{\frac{p}{2}+1} d t  \ < +\infty.
  \end{align*}
\item {\bf (Coupling) } for every $x\in\mathbb{R}$ and $t\in[0,T]$, the image measure of $\indic{z\in B(\varepsilon^*)} \nu_t(dz)$ by the map $z \mapsto c(t,x,z)$ has a density with respect to the Lebesgue measure on $\mathbb{R}$. 
\end{itemize}
\end{hyp}

In particular, the coupling condition in the hypothesis \ref{hyp:ar} is satisfied if the map $z\mapsto c(t,x,z)$ is continuously differentiable on $\RR$ and $\int_{\RR} \indic{\frac{\partial c}{\partial z}(t,x,z)\neq 0} \nu_t(dz) = +\infty$ for all  $x\in \RR$ and $t\in[0,T]$, as outlined in \cite{fournier_giet}. For example, this is trivially satisfied in the case of an SDE driven by a truncated stable process (i.e $c(t,x,z) = \bar{c}(t,x) z$ and $\nu_t(dz) = \indic{|z| \leq 1} \frac{dz}{|z|^{1+\alpha}}$).

We may now state the main result of this section, which gives a bound for the $L^p$-distance between the solution of \eqref{eq:intro_SDEs} and the version $\oeeX$ of the scheme \eqref{eq:numerical_scheme}:
\begin{equation}
\begin{aligned}
    \oeeX_{t_i} = \oeeX_{t_{i-1}} + a(t_{i-1},\oeeX_{t_{i-1}}) (t_{i}-t_{i-1}) + b(t_{i-1},\oeeX_{t_{i-1}}) (W_{t_i} - W_{t_{i-1}}) \\ + T_i \left(\oeeX_{t_{i-1}}\right) + \int^{t_i}_{t_{i-1}} \int_{\mathbb{R} \backslash B(\varepsilon)} c({s},\oeeX_{t_{i-1}},z) \rpm,
\end{aligned}
\end{equation}
where $(T_i)_{0 \leq i \leq n}$ is a collection of measurable applications wisely chose to transport the small jump part of the process to a Gaussian distribution with the same variance in an optimal way with respect to the $\mathcal{W}_p$ distance.
The detailed construction, as well as the proof of the theorem are postponed to Section \ref{sec:maiTh_proof}. 

\begin{theorem}[\bf Strong convergence of the $\varepsilon$-EM scheme \eqref{eq:numerical_scheme}] \label{theoreme_maurer} Let $p\geq 2$. Assume that \ref{hyp:lipschitz}, \ref{hyp:peano} and \ref{hyp:ar} hold for that $p$ and a given $\varepsilon^*$. 
\begin{itemize}
\item[$(i.)$ ] For any $\varepsilon \in (0,\varepsilon^*]$, there exists a  sequence $(\oeeX_{t_i})_{0 \leq i \leq n}$ of random variables on $(\Omega,\mathcal{F})$, such that for any $0 \leq i \leq n$, $\oeeX_{t_i}$ is $\mathcal{F}_{t_i}$- measurable, and verifies $\tlaw(\oeeX_{t_i}) = \tlaw(\oeX_{t_i})$. 
Moreover, there exists $\overline{m}(p,T) > 0$ such that 
\begin{align} \label{eq:moment_bound_scheme}
\sup_{\varepsilon  \in [0,\varepsilon^*)}\mathbb{E} \big [ \sup_{0 \leq i \leq n}| \oeeX_{t_i} |^{ p} \big ] \leq  \overline{m}(p,T). 
\end{align}

\item[$(ii.)$ ] The following inequality stands true for any $\varepsilon \in (0,\varepsilon^*]$ and $\gamma, \zeta  \in (0,1]$ given in \ref{hyp:peano}: 
\begin{equation} \label{eq:rate_of_convergence_scheme}
        \Big\| \sup_{0 \leq i \leq n} \big| X_{t_i} - \oeeX_{t_i} \big| \Big\|_{L^p(\Omega)} \leqpt  n^{- \left\{ \gamma \wedge \left( \frac{2 \zeta }{p (1+\zeta)} \right) \right\} } + \delta_p^n(\varepsilon),
\end{equation}
\begin{equation}\label{eq:Delta_n_epsilon}
\text{where} \qquad \delta_p^n(\varepsilon)  \coloneqq \left(
\displaystyle \sum_{k = 1}^n \mathbb{E} \left[ \left( \frac{\int_{t_{k - 1}}^{t_k} \int_{B
(\varepsilon)} | c (s, \oeeX_{t_{k - 1}}, z) |^{p + 2} \nu_s
(d z) d s}{\int_{t_{k - 1}}^{t_k} \int_{B (\varepsilon)} 
| c (s,\oeeX_{t_{k - 1}}, z) |^2 \nu_s (d z) d s} \right) ^ {\frac{1}{p}} \right]^2 \right) ^{\frac{1}{2}}
\end{equation} 
satisfies  $\lim_{\varepsilon \rightarrow 0} \delta_p^n(\varepsilon) = 0$ when the following sufficient condition holds: 
\begin{align} \label{eq:c_tend_vers_zero}
   \lim_{|z| \rightarrow 0} \sup_{t \in [0,T]} |\overline{L}_c(t,z)| = 0.
\end{align}
\end{itemize}
\end{theorem}

The result obtained in Theorem \ref{theoreme_maurer} can therefore be viewed as a path convergence \textit{in law}, which nonetheless precisely informs us about the behaviour of the scheme $\oeX$ in terms of the parameters $\varepsilon$ and $n$, as confirmed by the numerical simulations presented in Section \ref{sec:numerics}.\smallskip

Condition \eqref{eq:c_tend_vers_zero} imposes on the jump  coefficient $c(t,x,z)$ in the SDE \eqref{eq:intro_SDEs} to not amplify both negative and positive small jumps of the {Poisson random} measure $\widetilde{N}$, a behaviour  which appears to be essential to provide the convergence in $\varepsilon$. For instance, one may consider the non-symmetric time-homogeneous Lévy measure $\nu_t(dz) = z^2 \indic{0<z\leq1} dz + \frac{1}{z^2} \indic{-1 \leq z < 0} dz$, with the jump coefficient $c(t,x,z) = \frac{1}{z^{1/2}} \indic{0<z\leq1} + |z|^{3/2} \indic{-1 \leq z \leq 0} $. In this case, the SDE \eqref{eq:intro_SDEs} is well-posed in $L^2(\Omega)$ and satisfies the assumptions of Theorem \ref{theoreme_maurer} while, one can easily observe that $\delta_2^n(\varepsilon)$ diverges:
\begin{equation*}
    \left(\delta_2^n(\varepsilon)\right)^2 \geq 
n \left( \frac{ \int_{0}^{\varepsilon} \frac{1}{z^2} z^2
d z }{ \int_{0}^{\varepsilon} \frac{1}{z} z^2 dz  + \int_{-\varepsilon}^{0} |z|^{3} \frac{1}{z^2} dz} \right) =  {\frac{n}{\varepsilon}} \xrightarrow[\varepsilon\to 0^+]{} +\infty.
\end{equation*}

\medskip
An interesting case arises when the jump coefficient $c$ is sublinear with respect with $z$, allowing to determine an optimal choice for $\varepsilon$ that aligns with the  convergence rate of the Euler-Peano scheme  in Lemma \ref{lem_peano}. More precisely, we have the following corollary, the proof of which is postponed in Section \ref{sec:proof_coro_main}.
\begin{corollary}[\bf Optimal choice of $\varepsilon$] \label{cor:mainresult}
In addition to the assumptions of Theorem \ref{theoreme_maurer}, we assume that there exists a constant $C_T$ satisfying 
\begin{equation} \label{eq:bound_on_the_jump_coefficient}
\forall (s,x,z) \in [0,T] \times \mathbb{R} \times B(\varepsilon^*), \quad |c(t,x,z)| \leq  C_T |z| (1+|x|).
\end{equation} 
Then, for any $\varepsilon \in (0,\varepsilon^*]$,  the  $L^{p}$-strong error of the ($(\Omega, \mathcal{F})$-representation of the)   $\varepsilon$-EM scheme  \eqref{eq:numerical_scheme} satisfies:
\begin{equation} \label{eq:explicit_varepsilon}
\Big\| \sup_{0 \leq i \leq n} \big| X_{t_i} - \oeeX_{t_i} \big| \Big\|_{L^p(\Omega)} \leqpt  n^{- \left\{ \gamma \wedge \left( \frac{2 \zeta }{p (1+\zeta)} \right) \right\} } +  \varepsilon \sqrt{n}. 
\end{equation}
Moreover, suppose we have $\psi_p \in L^2([0,T])$ (i.e. $\zeta=1$ in \ref{hyp:peano}). With $\varepsilon$  taken such that 
\begin{equation}\label{eq:optimal_varepsilon}
    \varepsilon \leq n^{ - \left( \frac{1}{2}+\left(\gamma \wedge \frac{1}{p}\right)\right)} \wedge \varepsilon^*,
\end{equation} 
we obtain the following convergence rate for the $L^{p}$-strong error of the scheme \eqref{eq:numerical_scheme}:
\begin{equation} \label{strong_error_rate}
    \Big\| \sup_{0 \leq i \leq n} \left| X_{t_i} - \oeeX_{t_i} \right| \Big\|_{L^{p}(\Omega)} 
    \leqpt  n^{- \left( \gamma \wedge  \frac{1}{p} \right) }.
  \end{equation}
\end{corollary}

This recovers the result announced  in the introduction  when the time-H{\"o}lder regularity exponent $\gamma$ is  at least $\tfrac{1}{p}$. 
For numerical efficiency,  $\varepsilon$ must be chosen larger as possible (saving the simulation effort of very small jumps of the stochastic Poisson integral), and the above corollary states that the choice $\varepsilon = \mathcal{O}\big( n^{ - \left( \frac{1}{2}+\left(\gamma \wedge \frac{1}{p}\right)\right)}\big)$ is optimal to guarant{ee} the rate $\mathcal{O}\big(n^{- \left( \gamma \wedge  \frac{1}{p} \right)} \big)$. 

\begin{remark} \label{rem:low_integrability}
One may relax the hypothesis  $\psi_p \in L^2([0,T])$ in Corollary \ref{cor:mainresult}. In particular, if $\psi_p(t) \sim \frac{1}{t^\varrho}$ at the neighbourhood of zero with $\varrho \in (0,\frac{1}{2}]$, one can derive from Theorem \ref{theoreme_maurer}  the upper-bound
 \begin{align} \label{strong_error_rate_remarck}
    \Big\| \sup_{0 \leq i \leq n} \Big| X_{t_i} - \oeeX_{t_i} \Big| \Big\|_{L^{p}(\Omega)} &  \leqpt  n^{ - \left(\gamma \wedge \left( \frac{2(1-\varrho)}{p} - \kappa \right) \right)}
  \end{align}
   for any $\kappa > 0$. 

   {Indeed, since $\psi_p(t) \sim \frac{1}{t^\varrho}$, one gets $\psi_p \in L^{1+\zeta}([0,T])$ with $\zeta = \frac{1}{\varrho} - 1 -r$ for any $r > 0$.
   Theorem \ref{theoreme_maurer} gives the  $L^p$-error upper-bound  $n^{-\left(\gamma \wedge \frac{2 \zeta}{p (1+\zeta)}\right)} + \delta_p^n(\varepsilon).$
    Observing that  $ \frac{2 \zeta}{p (1+\zeta)} \geq \frac{2(1-\varrho - r  \varrho)}{p}$, and by setting $\kappa = \frac{2r \varrho}{p}$,  we get the upper-bound  $n^{-\gamma \wedge \left(\frac{2(1-\varrho)}{p} - \kappa\right)} + \delta_p^n(\varepsilon)$ for any $\kappa > 0$.
    Thanks to \eqref{eq:bound_on_the_jump_coefficient}, choosing $\varepsilon \leq n^{-\frac{1}{2} - \gamma \wedge \left( \frac{2(1-\varrho)}{p} - \kappa \right)}$, following the proof of Corollary \ref{cor:mainresult}, the second term,  $\delta_p^n(\varepsilon)\leqpt \sqrt{n} \varepsilon $,  can also be made smaller than $n^{-\gamma \wedge \left(\frac{2(1-\varrho)}{p} - \kappa\right)}$.} 
    In Section \ref{ssec:numerical_low_integrability}, we  provide numerical tests that confirm the influence of the parameter $\varrho$ on  the convergence rate. 
\end{remark}

The proof of Theorem \ref{theoreme_maurer} can be easily adapted to derive a convergence rate for the $\varepsilon$-EM scheme without Gaussian compensation (i.e when the jumps smaller than $\varepsilon$ are neglected). In this case, there is no need for a Wasserstein coupling and one can work directly with $\ueX$. The proof of this theorem can be found in Section \ref{proof:thmnoapprox}. 

\begin{theorem}[\bf Strong convergence of the $\varepsilon$-EM scheme without  substitute \eqref{eq:numerical_scheme_without_s_jumps}]  \label{thm:noapprox}
Let $p\geq 2$. Assume that \ref{hyp:lipschitz} and \ref{hyp:peano} holds.
Then the $\varepsilon$-EM scheme without Gaussian substitute $\ueX$ satisfies: 
\begin{align}
        \sup_{\varepsilon > 0} \EE[ \sup_{0 \leq i \leq n} |\ueX_{t_i}|^p ] \leq \underline{m}(p,T),
\end{align}
and for any $\varepsilon > 0$, we have the convergence rate
\begin{align}
   \left\| \sup_{0 \leq i \leq n} |X_{t_i} - \ueX_{t_i}| \right\|_{L^p(\Omega)} \leqpt n^{- \left\{ \gamma \wedge \left( \frac{2 \zeta }{p (1+\zeta)} \right) \right\} } + \left( \int_0^{T} \psi_{p,\varepsilon}(s) ds\right)^\frac{1}{p},
\end{align}
where $\psi_{p,\varepsilon}$ stands for the $\psi_p$ function integrated only on $B(\varepsilon)$:
\begin{align*}
    \psi_{p,\varepsilon}(s) =  \left( \int_{B(\varepsilon)} \overline{L}_c^2(s,z) \nu_s(dz) \right)^{p / 2} +  \int_{B(\varepsilon)} \overline{L}_c^p(s,z) \nu_s(dz).
\end{align*}
\end{theorem}

For this scheme as well, we can derive a simpler convergence rate when the coefficient $c$ is sub-linear with respect to the $z$ variable.

\begin{corollary}[\bf Optimal choice of $\varepsilon$]\label{cor:scheme_cut_BG_index}
Assume that \eqref{eq:bound_on_the_jump_coefficient} holds in addition of the hypothesis of Theorem \ref{thm:noapprox}. In addition, assume that the family of measures $(\nu_t)_{t \in [0,T]}$ is dominated by a Lévy measure $\mu$ (i.e $\nu_t(A) \leq \mu(A)$ for $t \in [0,T]$ and $A \in \mathcal{B}(\RR)$). Let $\beta = \inf\{\alpha > 0, \int_{|z|<1} |z|^{\alpha} \mu(dz) < +\infty \}$  the Blumenthal--Getoor index of the indivisible distribution characterised by the Lévy measure $\mu$. Then the following inequality stands true:
\begin{align*}
    \left\| \sup_{0 \leq i \leq n} |X_{t_i} - \ueX_{t_i}| \right\|_{L^p(\Omega)} \leqpt n^{- \left\{ \gamma \wedge \left( \frac{2 \zeta }{p (1+\zeta)} \right) \right\} } + \varepsilon^{1-\frac{\beta^+}{2}}.
\end{align*}
In the convergence above, the term $\varepsilon^{1-\frac{\beta^+}{2}}$ stands for $\sup_{\theta>0} \{ \varepsilon^{1-\frac{\theta}{2}}; \int_{|z|<1} |z|^{\theta} \mu(dz) < +\infty \}$. In particular, when $\mu$ is a truncated  alpha stable distribution, this term is $\varepsilon^{1-\frac{\alpha}{2}}$.

Moreover, suppose we have $\psi_p \in L^2([0,T])$ (i.e. $\zeta=1$ in \ref{hyp:peano}). With $\varepsilon$  taken such that 
\begin{equation}\label{eq:optimal_varepsilon_nosubstitute}
    \varepsilon \leq n^{ - \left(\gamma \wedge \frac{1}{p}\right) \frac{2}{2-\beta^+}},
\end{equation} 
we obtain the following convergence rate for the $L^{p}$-strong error of the scheme \eqref{eq:numerical_scheme_without_s_jumps}:
\begin{equation} \label{strong_error_rate_nosubstitute}
    \Big\| \sup_{0 \leq i \leq n} \left| X_{t_i} - \ueX_{t_i} \right| \Big\|_{L^{p}(\Omega)} 
    \leqpt  n^{- \left( \gamma \wedge  \frac{1}{p} \right) }.
  \end{equation}

\begin{remark}
If the family of compensator measures $(\nu_t)_{t\in[0,T]}$ has finite activity, i.e if $\int_{\mathbb{R}} \nu_t(dz) < +\infty$ for all $t\in[0,T]$, then the scheme $\ueX$ with $\varepsilon = 0$ is well defined and the results of Theorem \ref{thm:noapprox} and Corollary \ref{cor:scheme_cut_BG_index} are valid with $\varepsilon = 0$. 
\end{remark} 
\end{corollary}

\subsection{Proofs preamble}
To some exten{t}, the proofs make use of martingale $L^p$-estimates. For continuous-time martingales, one usually rely on the Burk{ho}lder--Davis--Gundy inequality
\begin{align} \label{eq:bdg}
    \mathbb{E}\left[\sup_{0 \leq s \leq t} |M_s|^p\right] \leq \cpbdg  \  \mathbb{E}[[M,M]_t^{p/2}],
\end{align}
where $[M,M]$ is the quadratic variation of the càdlàg martingale $M$. When $M$ is continuous, $[M,M]$ is equal to the predictable quadratic variation $\langle M,M \rangle$. However, the inequality
\begin{align*}
    \mathbb{E}\left[\sup_{0 \leq s \leq t} |M_s|^p\right] \leq \cpbdg \  \mathbb{E}[\langle{M,M}\rangle_t^{p/2}]
\end{align*}
is not valid for martingales with jumps (see the introduction in Breton and Privault \cite{breton}, which refers to Remark 357 in Situ \cite{situ}). It follows that the  Burkh{o}lder--Davis--Gundy inequality \eqref{eq:bdg} is not useful in practice to estimate {stochastic Poisson} integrals, since the quadratic variation of the latter cannot be expressed in term of integrals with respect to its compensator measure in general. 
    Instead, we will use the following inequality stated by  Kunita \cite[Theorem 2.11]{kunita}, and extended in the time-inhomogeneous setting by Breton and Privault  \cite{breton}.
\begin{proposition}[\bf Kunita's BDG inequality for jump processes, Corollary 2.1 in \cite{breton}.] \label{prop:kunita}
  Consider the stochastic Poisson integral
  \begin{align*}
    K_t & = \int_0^t \int_{-\infty}^{+ \infty} w ({s}, z)  \rpm,
  \end{align*}
  with $(w (s, z))_{(s, z) \in \mathbb{R}_+ \times \mathbb{R}}$ being a
  predictable process. Then for any $p \geq 2$, there exists $\ckunp > 0$ such that
   \begin{align} \label{kunita}
 \mathbb{E} \left[ \Big( \sup_{t \in [0, T]} |K_t| \Big)^p
     \right] \leq & \ckunp
  \Biggr\{ \mathbb{E} \left[  \left( \int_0^T  \int_{- \infty}^{+ \infty} | w (t,z) |^2 \nu_t (d z) d t  \right)^{p / 2} \right] + \mathbb{E} \left[ \int_0^T
    \int_{- \infty}^{+ \infty} | w (t, z) |^p \nu_t (d z) d t \right] \Biggr\}.
   \end{align}
Alternatively: 
\begin{align} \label{kunita_2}
 \mathbb{E} \left[ \Big(\sup_{t \in [0, T]} |K_t| \Big)^p
    \right] \leq & \ckunpT  \int_0^T \mathbb{E} \left[  \left( \int_{- \infty}^{+ \infty} | w (t,z) |^2 \nu_t (d z) \right)^{p / 2} \right] dt + \ckunp \mathbb{E} \left[ \int_0^T
    \int_{- \infty}^{+ \infty} | w (t, z) |^p \nu_t (d z) d t \right]
  \end{align}
with $\ckunpT = T^{\frac{p}{2}-1} \ckunp$.
\end{proposition}

While quantifying the error made by the approximation of the small jumps, the proofs deal with  discrete-time martingales, that are handled with the  discrete version of  Burkh{o}lder--Davis--Gundy inequality:

\begin{theorem}[\bf Theorem 1.1 in \cite{burkholder1972integral}] \label{thm:discrete_BDG}
    Let $p \geq 2$, $m$ be an integer, and $(\mathcal{F}_k)_{0\leq k \leq m}$ be a discrete-time filtration and $(d_n)_{0\leq k \leq m}$ a sequence of $\mathcal{F}_m$-martingale differences. The following inequality holds:
    \begin{align*}
        \Big\| \sup_{0\leq k \leq m} \Big| \sum_{l=0}^k d_l \Big| \Big\|_{L^p(\Omega)} \leq \kappa_p \Big\| \sum_{k=0}^m d_k^2 \Big\|_{L^{\frac{p}{2}}(\Omega)}^\frac{1}{2},  
    \end{align*}
    where the constant $\kappa_p$ only depends on $p$.
\end{theorem}

\subsection{Proof of Proposition \ref{prop:wasserstein}} \label{sub:proof_wasserstein}

The proof relies on Theorem 1.1 in Bobkov  \cite{bobkov}, which extends the results of Rio in \cite{rio} on the Berry--Essen bounds for sums of independent random variables in Wasserstein
distance $\mathcal{W}_q$ for any $q \geq 1$. We reproduce it here for the sake of completeness. 
\begin{theorem}[\bf Theorem 1.1 in \cite{bobkov} ] \label{thm:Berry--Essen}
For $q\geq 1$, there exists $c_q > 0$ depending only on $q$ such that if $X_1, \ldots,
  X_m$ are independent random variables with $\sum_{j = 1}^m {\Var} (X_j) =
  1$, then
  \begin{align*}
    \mathcal{W}_q \left( \law{\sum_{j = 1}^m X_j}, \ 
    \mathcal{N} (0, 1) \right) & \leq c_q  \left( \sum_{j = 1}^m \mathbb{E}
    [| X_j |^{q + 2}] \right)^{1/q}.
  \end{align*}
\end{theorem}

After a re-normalisation, we can get rid of the variance condition and derive from Theorem \ref{thm:Berry--Essen} that for every $q \geq 1$, 
\begin{align} \label{eq:bobkov_renorm}
    \mathcal{W}_{q} \left( \tlaw\left(  \sum_{j = 1}^m X_j \right),
    \mathcal{N} \left( 0, \sum_{j = 1}^m \Var  (X_j) \right) \right) & \leq c_{q}\  \left(\frac{\sum_{j = 1}^m \mathbb{E} [| X_j
    |^{q + 2}]}{\sum_{j = 1}^m \Var  (X_j)} \right)^{1/q}.
\end{align}
Let $\varepsilon  \in (0,\varepsilon_F]$, $t \in [0, T]$, $m  \in \mathbb{N}^{\ast}$, and set $\tau_j = t_0 + \frac{j}{m}  {(t-t_0)}$ for $j \in \{0,\ldots,m\}$. We apply the latter inequality to the random variables 
\[X_j = \int_{\tau_{j-1}}^{\tau_j} \int_{B(\varepsilon)} F({s},z) \rpm, \quad \mbox{for $j \in \{ 1, \ldots, m \}$},\]
which are independents as the increments of the additive process $\left(\int_{t_0}^{t} \int_{B(\varepsilon)} F({s},z) \rpm \right)_{t\in[0,T]}$ and verify
\begin{align*}
\sum_{j = 1}^m X_j = \int_{t_0}^{t} \int_{B(\varepsilon)} F({s},z) \rpm \quad \text{ and } \quad
\sum_{j = 1}^m \Var  (X_j) = \int_{t_0}^t \int_{B(\varepsilon)} | F(s,z) |^2  \ \nu_s (d z) d s. 
\end{align*}
Now we evaluate $\sum_{j = 1}^m \mathbb{E} [| X_j |^{ q + 2}]$. To that extend, we make use of Kunita's inequality \eqref{kunita} to obtain
\begin{align} \label{eq:moment_2p_plus_2}
      \mathbb{E} \left[ |X_j|^{q+2} \right] & \leq {C^{\tiny{\textit{Kun}}}_{q+2}} \int_{\tau_{j-1}}^{\tau_j} \int_{B(\varepsilon)} |F(s,z)|^{q+2} \rpc +  {C^{\tiny{\textit{Kun}}}_{q+2}}  \left( \int_{\tau_{j-1}}^{\tau_j} \int_{B(\varepsilon)} |F(s,z)|^{2} \rpc \right) ^ \frac{q+2}{2}.
\end{align}
By hypothesis \eqref{eq:hypo_propo:wasserstein}, for any $j \in \{0,\ldots,m\},$
\begin{equation*}
      \int_{\tau_{j-1}}^{\tau_j} \left( \int_{B(\varepsilon)} |F(s,z)|^{2} \nu_s(dz) \right) ^ \frac{q+2}{2} ds \leq  \int_{0}^{T} \left( \int_{B(\varepsilon^*)} |F(s,z)|^{2} \nu_s(dz) \right) ^ \frac{q+2}{2} ds  < +\infty,
  \end{equation*}
so one may use Hölder inequality on the second integral in \eqref{eq:moment_2p_plus_2} to get
  \begin{align*}
   \int_{\tau_{j-1}}^{\tau_j} \int_{B(\varepsilon)} |F(s,z)|^{2} \rpc & \leq \left( \int_{\tau_{j-1}}^{\tau_j} ds \right) ^ {\frac{q}{q+2}} \left( \int_{\tau_{j-1}}^{\tau_j} \left( \int_{B(\varepsilon)} |F(s,z)|^{2} \nu_s(dz) \right) ^ \frac{q+2}{2} ds \right) ^ {\frac{2}{q+2}}.   
  \end{align*}
Then, plugging this bound into \eqref{eq:moment_2p_plus_2}, we obtain 
  \begin{align*} 
      \mathbb{E} \left[ |X_j|^{q+2} \right] & \leq \int_{\tau_{j-1}}^{\tau_j} \int_{B(\varepsilon)} |F(s,z)|^{q+2} \rpc + \left(\frac{t-t_0}{m} \right) ^ \frac{q}{2} \int_{\tau_{j-1}}^{\tau_j} \left( \int_{B(\varepsilon)} |F(s,z)|^{2} \nu_s(dz) \right) ^ \frac{q+2}{2} ds.
  \end{align*}
 Taking the sum from $1$ to $m$ and dividing by the variance, we derive from \eqref{eq:bobkov_renorm} the following bound for the Wasserstein distance:
  \begin{multline*}
      \mathcal{W}_{q} \left( \law{ \int_{t_0}^{t} \int_{B(\varepsilon)} F({s},z) \rpm },
    \mathcal{N} \left( 0, \int_{t_0}^t \int_{B(\varepsilon)} | F(s,z) |^2 \nu_s (d z) d s \right) \right)^{
    q} \\ \leq (c_{q})^{q}  \ \left( \frac{\int_{t_0}^t \int_{B(\varepsilon)}  |F(s,z)|^{q+2} \nu_s (d z) d s}{\int_{t_0}^t
    \int_{B(\varepsilon)}  | F(s,z) |^2 \nu_s (d z) d s} +
    \dfrac{T^\frac{q}{2}}{m^{\frac{q}{2}}} \frac{\int_{t_0}^{t} \left( \int_{B(\varepsilon)} |F(s,z)|^{2} \nu_s(dz) \right) ^ \frac{q+2}{2} ds}{\int_{t_0}^t
    \int_{B(\varepsilon)}  | F(s,z) |^2 \nu_s (d z) d s} \right).
\end{multline*}
  We end the proof by setting $\mathcal{A}(q) = (c_{q})^{q}$ and letting $m \rightarrow +  \infty$. The second term in the right hand side converges to zero since $q \geq 1$ and we obtain the desired result.
  
\subsection{Convergence of the Euler-Peano scheme (proof of Lemma \ref{lem_peano})} \label{sec:peano}

Set $\mathcal{E}_t \coloneqq |X_t - \widetilde{X}_t|$ and $\mathcal{E}_t^{\ast} = \sup_{t\in[0,t]} \mathcal{E}_s $. We now estimate $ \| \mathcal{E}_t^{\ast} \|_{L^p(\Omega)} $. Using Minkowski's integral inequality along with Burkh{o}lder--Davis--Gundy inequality and Kunita  inequality in Proposition \ref{prop:kunita}, one has  
\begin{align*}
    \| \mathcal{E}_t^{\ast} \|_{L^p(\Omega)} \leq & \int_0^t \left\| a (s, X_s) - a (\eta (s^{{-}}), \widetilde{X}_{\eta
    (s^{{-}})})  \right\|_{L^p(\Omega)} d s  + \left(\cpbdg\right)^\frac{1}{p}  \left( \int_0^t \left\| b (s, X_{s}) - b (\eta (s^{{-}}), \widetilde{X}_{\eta
    (s^{{-}})}) \right\|_{L^p(\Omega)}^2  d s  \right)^\frac{1}{2} \\
    & + \left(\ckunpT\right)^\frac{1}{p} \left( \int_0^t \mathbb{E} \left[ \left( \int_{-\infty}^{+\infty}  (c (s, X_{s^{{-}}},z) - c (s, \widetilde{X}_{\eta
    (s^{{-}})},z))^2 \nu_s(dz) \right)^\frac{p}{2} \right] ds \right)^\frac{1}{p} \\ & + \left(\ckunp\right)^\frac{1}{p} \left( \int_0^t \int_{-\infty}^{+\infty}  \mathbb{E} \left[|c (s, X_{s^{{-}}},z) - c (s, \widetilde{X}_{\eta
    (s^{{-}})},z)|^p \right] \nu_s(dz)  ds \right)^\frac{1}{p}.
\end{align*}
Then we apply \ref{hyp:lipschitz} and \ref{hyp:peano} to derive
\begin{align*}
    \| \mathcal{E}_t^{\ast} \|_{L^p(\Omega)} \leq & {\overline{L}_{a,b}} \int_0^t \left( |s - \eta(s)|^\gamma + \left\| X_s - \widetilde{X}_{\eta
    (s)}  \right\|_{L^p(\Omega)} \right) ds \\ & + \left(\cpbdg\right)^\frac{1}{p} {\overline{L}_{a,b}}   \left( \int_0^t \left( |s - \eta(s)|^\gamma + \left\| X_s - \widetilde{X}_{\eta
    (s)}  \right\|_{L^p(\Omega)} \right)^2 ds \right)^\frac{1}{2} \\
    & + \left(\ckunpT\right)^\frac{1}{p} \left( \int_0^t  \left( \int_{-\infty}^{+\infty}  L_c(s,z)^2 \nu_s(dz) \right)^\frac{p}{2} \left\| X_s - \widetilde{X}_{\eta(s)} \right\|_{L^p(\Omega)}^p ds \right)^\frac{1}{p} \\ & + \left(\ckunp\right)^\frac{1}{p} \left( \int_0^t \left( \int_{-\infty}^{+\infty}  L_c(s,z)^p \nu_s(dz) \right) \left\| X_s - \widetilde{X}_{\eta(s)} \right\|_{L^p(\Omega)}^p ds \right)^\frac{1}{p}.
\end{align*}
Recalling the definition of $\psi_p(t)$ in \eqref{eq:def_psi_p},  we rewrite this estimate as
\begin{align*}
    \| \mathcal{E}_t^{\ast} \|_{L^p(\Omega)} \leq & {\overline{L}_{a,b}}  \int_0^t \left( |s - \eta(s)|^\gamma + \left\| X_s - \widetilde{X}_{\eta
    (s)}  \right\|_{L^p(\Omega)} \right) ds \\ & + \left(\cpbdg\right)^\frac{1}{p} {\overline{L}_{a,b}}   \left( \int_0^t \left( |s - \eta(s)|^\gamma + \left\| X_s - \widetilde{X}_{\eta
    (s)}  \right\|_{L^p(\Omega)} \right)^2 ds \right)^\frac{1}{2} \\
    & + 2 \ \left(\ckunpT\right)^\frac{1}{p} \left( \int_0^t  \psi_p(s) \left\| X_s - \widetilde{X}_{\eta(s)} \right\|_{L^p(\Omega)}^p ds \right)^\frac{1}{p}.
\end{align*}
Then the triangle inequality $\| X_t - \widetilde{X}_{\eta
    (t)}\|_{L^p(\Omega)} \leq \| X_t - \widetilde{X}_t \|_{L^p(\Omega)} + \| \widetilde{X}_t - \widetilde{X}_{\eta (t)} \|_{L^p(\Omega)} $ allows to recover the global error term $\mathcal{E}_t$ and the local error term $\widetilde{X}_t - \widetilde{X}_{\eta(t)}$.
Using $|s-\eta(s)|^\gamma \leq (\frac{T}{n})^\gamma $ and the $L^2([0,t])$ and $L^p([0,t])$-Minkowski inequalities, we obtain
\begin{equation} \label{eq:euler_error_expansion} 
\begin{aligned}
    \| \mathcal{E}_t^{\ast} \|_{L^p(\Omega)} \leq & {\overline{L}_{a,b}}  \left( \frac{T^{\gamma+1}}{n^\gamma} + \int_0^t \| \widetilde{X}_s - \widetilde{X}_{\eta(s)} \|_{L^p(\Omega)} ds + \int_0^t \| \mathcal{E}_s^{\ast} \|_{L^p(\Omega)} ds \right) \\
    & + \left(\cpbdg\right)^\frac{1}{p} {\overline{L}_{a,b}}  \left( \frac{T^{\gamma+{\frac{1}{2}}}}{n^\gamma} + \left( \int_0^t \| \widetilde{X}_s - \widetilde{X}_{\eta(s)} \|_{L^p(\Omega)}^2 ds \right)^\frac{1}{2} + \left( \int_0^t \| \mathcal{E}_s^{\ast} \|_{L^p(\Omega)}^2 ds \right)^\frac{1}{2} \right) \\
    & + 2 \ \left(\ckunpT\right)^\frac{1}{p} \left( \left( \int_0^t \psi_p(s) \| \widetilde{X}_s - \widetilde{X}_{\eta(s)} \|_{L^p(\Omega)}^p ds \right)^\frac{1}{p}+\left( \int_0^t \psi_p(s) \| \mathcal{E}_s^\ast \|_{L^p(\Omega)}^p ds \right)^\frac{1}{p}\right).
\end{aligned}
\end{equation}
To end the proof with a Gronwall lemma, one needs to bound the local error terms, namely $\int_0^t \| \widetilde{X}_s - \widetilde{X}_{\eta(s)} \|_{L^p(\Omega)} ds$, $\left( \int_0^t \| \widetilde{X}_s - \widetilde{X}_{\eta(s)} \|_{L^p(\Omega)}^2 ds \right)^\frac{1}{2}$ and $\left( \int_0^t \psi_p(s) \| \widetilde{X}_s - \widetilde{X}_{\eta(s)} \|_{L^p(\Omega)}^p ds \right)^\frac{1}{p}$ by some power of the discretisation step $\frac{T}{n}$, where $\psi_p$ is integrable on $[0,T]$ thanks to \ref{hyp:peano}.

\paragraph{$\bullet$ Estimation of $\|  \widetilde{X}_t - \widetilde{X}_{\eta(t)}  \|_{L^p (\Omega)}$.}
With a triangle inequality,
\begin{align*}
  \|  \widetilde{X}_t - \widetilde{X}_{\eta(t)}  \|_{L^p (\Omega)} \leq & \left\|
  \int_{\eta(t)}^t a (\eta (s), \widetilde{X}_{\eta (s)}) d s \right\|_{L^p (\Omega)} 
  + \left\| \int_{\eta(t)}^t b (\eta (s), \widetilde{X}_{\eta (s)}) d W_s
  \right\|_{L^p (\Omega)}\\
  & + \left\| \int_{\eta(t)}^t \int_{- \infty}^{+\infty} c ({s}, \widetilde{X}_{\eta
  (s^{-})}, z)  \widetilde{N} (d s, d z) \right\|_{L^p (\Omega)}.
\end{align*}
Then using \ref{hyp:lipschitz} and standard arguments, we upper-bound the three terms above as follows:
\begin{gather}\label{eq:local_estimate_a_et_b}
        \left\| \int_{\eta(t)}^t a (\eta (s), \widetilde{X}_{\eta (s)}) d s \right\|_{L^p (\Omega)} \leqpt n^{-1}, \qquad 
        \left\| \int_{\eta(t)}^t b (\eta (s), \widetilde{X}_{\eta (s)}) d W_s \right\|_{L^p (\Omega)} \leqpt n^{-1/2} \notag \\
         \left\| \int_{\eta(t)}^t c ({s}, \widetilde{X}_{\eta(s^-)}, z) \rpm\right\|_{L^p (\Omega)} ^p  \leqpt \int_{\eta(t)}^t \psi_p(s) ds, 
\end{gather}
leading to
\begin{align}\label{eq:estimate_of_local_error}
    \|  \widetilde{X}_t - \widetilde{X}_{\eta(t)}  \|_{L^p (\Omega)} \leqpt n^{-1/2} + \left( \int_{\eta(t)}^t \psi_p(s) ds \right)^{\frac{1}{p}}.
\end{align}

\paragraph{$\bullet$ Estimation of $\int_0^t \| \widetilde{X}_s - \widetilde{X}_{\eta(s)} \|_{L^p(\Omega)} ds$ and $\left( \int_0^t \| \widetilde{X}_s - \widetilde{X}_{\eta(s)} \|_{L^p(\Omega)}^2 ds \right)^\frac{1}{2}$.}

We derive from \eqref{eq:estimate_of_local_error} that 
\begin{align*}
    \int_0^t \| \widetilde{X}_s - \widetilde{X}_{\eta(s)} \|_{L^p(\Omega)} ds & \leqpt T n^{-1/2} + \int_0^t \left( \int_{\eta(s)}^s \psi_p(u) du \right)^{\frac{1}{p}} ds, \\
    \int_0^t \| \widetilde{X}_s - \widetilde{X}_{\eta(s)} \|_{L^p(\Omega)}^2 ds & \leqpt T^2 n^{-1} + \int_0^t \left( \int_{\eta(s)}^s \psi_p(u) du \right)^{\frac{2}{p}} ds.
\end{align*}
It follows from Jensen inequality (since $p \geq 2$) that 
\begin{align} \label{eq:jensen_le}
    \int_0^t \| \widetilde{X}_s - \widetilde{X}_{\eta(s)} \|_{L^p(\Omega)} ds & \leqpt T n^{-1/2} + T^{1-\frac{1}{p}} \left( \int_0^t  \int_{\eta(s)}^s \psi_p(u) du \ ds\right)^{\frac{1}{p}} \\
    \int_0^t \| \widetilde{X}_s - \widetilde{X}_{\eta(s)} \|_{L^p(\Omega)}^2 ds & \leqpt T^2 n^{-1} + T^{1-\frac{2}{p}} \left( \int_0^t  \int_{\eta(s)}^s \psi_p(u) du \ ds\right)^{\frac{2}{p}}.
\end{align}
Integrating by parts, we get 
\begin{align*}
    \int_0^t  \int_{\eta(s)}^s \psi_p(u) du \ ds = &  \int_{0}^{t} (s - \eta(s)) \psi_p(s) ds \leq n^{-1} \int_{0}^{t} \psi_p(s) ds,
\end{align*}
leading to
\begin{align} \label{eq:ipp}
\begin{aligned}
   \left( \int_0^t  \int_{\eta(s)}^s \psi_p(u) du \ ds \right)^\frac{1}{p} \leq & \ n^{-1/p} \ \left( \int_0^t \psi_p(s) ds \right)^\frac{1}{p}, \\ 
\text{and} \quad   \left( \int_0^t  \int_{\eta(s)}^s \psi_p(u) du \ ds \right)^\frac{2}{p} \leq & \ n^{-2/p} \ \left( \int_0^t \psi_p(s) ds \right)^\frac{2}{p}.
\end{aligned}
\end{align}
By \ref{hyp:peano}, the integral $\int_0^t  \psi_p(s) ds$ is finite, and plugging the estimates \eqref{eq:ipp} into \eqref{eq:jensen_le}, it follows that 
\begin{align*}
    \int_0^t \| \widetilde{X}_s - \widetilde{X}_{\eta(s)} \|_{L^p(\Omega)} ds +
    \left( \int_0^t \| \widetilde{X}_s - \widetilde{X}_{\eta(s)} \|_{L^p(\Omega)}^2 ds \right)^\frac{1}{2} \leqpt n^{-1/p}.
\end{align*}

\paragraph{$\bullet$ Estimation of $\left(\int_0^t \psi_p(s) \ \| \widetilde{X}_s - \widetilde{X}_{\eta(s)} \|_{L^p(\Omega)}^p ds\right)^{\frac{1}{p}}$. } Using \eqref{eq:estimate_of_local_error} again, we get directly
\begin{align*}
    \int_0^t \psi_p(s) \| \widetilde{X}_s - \widetilde{X}_{\eta(s)} \|_{L^p(\Omega)}^p ds \leqpt n^{-p/2} + \int_0^t \psi_p(s) \left(  \int_{\eta(s)}^s \psi_p(u) du \right) ds.
\end{align*}
Then the Chasles relation yields
\begin{align} \label{eq:chasles}
    \int_0^t \psi_p(s)  \left(  \int_{\eta(s)}^s \psi_p(u) du \right) ds = & \sum_{i=0}^{n-1} \int_{t_i}^{t_{i+1}} \psi_p(s) \left(  \int_{t_i}^s \psi_p(u) du \right) ds  = \frac{1}{2} \sum_{i=0}^{n-1} \left( \int_{t_i}^{t_{i+1}} \psi(s) ds \right)^2.
\end{align}
The next step justifies our requirement for $\psi_p$ to have a slightly higher integrability than $L^1([0,T])$.  At this point, we apply Hölder inequality to derive
\begin{equation} \label{eq:holder}
\begin{aligned}
    \Bigg(  \int_{t_i}^{t_{i+1}} \psi_p(s) ds \Bigg) ^ 2 \leq & (t_{i+1}-t_i) ^ {2\left(1-\frac{1}{1+\zeta}\right)} \Bigg( \int_{t_i}^{t_{i+1}} \psi_p(s)^{1+\zeta} ds \Bigg) ^ {\frac{2}{1+\zeta}}
    \\ & \leq n^{-\frac{2\zeta}{1+\zeta}} \Bigg( \int_{t_i}^{t_{i+1}} \psi_p(s)^{1+\zeta} ds\Bigg) \Bigg( \int_{0}^{T} \psi_p(s)^{1+\zeta} ds \Bigg) ^ \frac{1-\zeta}{1+\zeta},
\end{aligned}
\end{equation}
where the integrals are finite by \ref{hyp:peano}, and $\frac{1-\zeta}{1+\zeta} \geq 0$. Note that this imposes the constraint  $\frac{2}{1+\zeta} \geq 1$ in \eqref{eq:holder} in order to obtain an upper-bound for the sum when coming back to \eqref{eq:chasles}. In particular,  a greater integrability than $L^2([0,T])$ for $\psi_p$ does not improve the convergence rate in our proof. We obtain 
\begin{align*}
    \int_0^t \psi_p(s) \| \widetilde{X}_s - \widetilde{X}_{\eta(s)} \|_{L^p(\Omega)}^p ds \leqpt n^{-p/2} + n^{-\frac{2\zeta}{1+\zeta}} \  \| \psi_p\|_{L^{1+\zeta}([0,T])}^{2},
\end{align*}
which finally resumes to
\begin{align*}
    \left( \int_0^t \psi_p(s) \| \widetilde{X}_s - \widetilde{X}_{\eta(s)} \|_{L^p(\Omega)}^p ds \right)^\frac{1}{p} \leqpt  n^{-\frac{2\zeta}{p(1+\zeta)}}.
\end{align*}
To come to an end, we may plug the local errors estimates back into \eqref{eq:euler_error_expansion} to obtain
\begin{align*}
    \| \mathcal{E}_t^{\ast} \|_{L^p(\Omega)} \leqpt& \ n^{- \left(\gamma \wedge \left( \frac{2 \zeta}{p(1+\zeta)}\right)\right)} +  \int_0^t \| \mathcal{E}_s^{\ast} \|_{L^p(\Omega)} ds +  \left( \int_0^t \| \mathcal{E}_s^{\ast} \|_{L^p(\Omega)}^2 ds \right)^\frac{1}{2} + \left(\int_0^t \psi_p(s) \  \| \mathcal{E}_s^\ast \|_{L^p(\Omega)}^p ds \right)^\frac{1}{p}.
\end{align*}
Using the Gronwall lemma \ref{lem:generalized_gronwall}, we conclude that
\begin{align*}
    \| \mathcal{E}_t^{\ast} \|_{L^p(\Omega)} \leqpt {n^{-{\gamma \wedge \left( \frac{2 \zeta}{p(1+\zeta)}\right)}}},
\end{align*}
which ends the proof. 

\subsection{Proof of Theorem \ref{theoreme_maurer}}\label{sec:maiTh_proof}

\subsubsection{Construction of $\oeeX$} \label{ssec:construction_schema}

Let $\varepsilon  \in [0,\varepsilon^*)$, and $1 \leq i \leq n$. For $x \in \mathbb{R}$ we set $Y_i(x) = \int_{t_{i-1}}^{t_i} \Eint c({s},x,z) \rpm$. Our aim consists in constructing a map $T_i : \mathbb{R}\times \mathbb{R} \rightarrow \mathbb{R}$ that optimally transports (for the $\mathcal{W}_p$-distance) $\law{Y_i(x)}$ to the centred normal distribution $\mathcal{N}_i (x)$ of variance $\int_{t_{i-1}}^{t_i} \Eint c(s,x,z)^2 \rpc$. In addition, this map must be measurable with respect to the parameter $x$.

To that extend, we rely on Theorem 1.1 in \cite{fontbona}.
As expounded within the same source, a prerequisite to the  application of this theorem is  the existence, for every $x\in\mathbb{R}$, of an optimal transport map from $\law{Y_i(x)}$ to $\mathcal{N}_i (x)$. In the case of the $\mathcal{W}_p$-distance, this requisite is satisfied as long as  $\law{Y_i(x)}$ admits a Lebesgue density, as detailed in \cite{fontbona}.  Considering the case of an SDE driven by time-homogeneous Poisson measure, the absolute continuity of the marginal laws of  the solution have been studied in     \cite{fournier_giet}. We apply \cite[Theorem 2.2]{fournier_giet} to our situation where $Y_i(x)$ can be seen as the solution of an SDE with constant coefficients, straightforwardly adapting the proof to the case of time-in homogeneity of $\nu_s(dz)$. This leads us to formulate the coupling hypothesis \ref{hyp:ar}.

We are now in position to construct $\oeeX$ recursively. We set $\oeeX_{t_0} = X_0$. For $1 \leq i \leq n$, let $\mathbb{Q}_i = \tlaw{(\oeX_{t_{i-1}})}$. Applying \cite[Theorem 1.1]{fontbona},\footnote{In the notations of \cite{fontbona}: $(E,\Sigma,m) = (\mathbb{R}, \mathcal{B}(\mathbb{R}),\mathbb{Q}_i)$, $\mu_\lambda = \tlaw{(Y_i(\lambda))}$ and $\nu_\lambda = \mathcal{N}_i (\lambda)$.} there exists an application $T_i : \mathbb{R}\times \mathbb{R} \rightarrow \mathbb{R}$ which is $(\mathcal{B}(\mathbb{R}) \otimes \mathcal{B}(\mathbb{R}) , \mathcal{B}(\mathbb{R}))$-measurable such that for $\mathbb{Q}_i$-almost every $x \in \mathbb{R}$, one has 
\begin{align} \label{eq:def_gj}
    \EE[|Y_i(x) - T_i (x,Y_i(x))|^p] = \mathcal{W}_p(\law{Y_i(x)},\mathcal{N}_i(x))^p.
\end{align}
Then, given $\oeeX_{t_{i-1}}$, we can set
\begin{equation}\label{eq:def_x_chapeau}
\begin{aligned}
    \oeeX_{t_i} =  & \oeeX_{t_{i-1}} + a(t_{i-1},\oeeX_{t_{i-1}}) (t_{i}-t_{i-1}) + b(t_{i-1},\oeeX_{t_{i-1}}) (W_{t_i} - W_{t_{i-1}}) \\ & \qquad + T_i \left(\oeeX_{t_{i-1}},Y_i(\oeeX_{t_{i-1}})\right) + \int^{t_i}_{t_{i-1}} \int_{\mathbb{R} \backslash B(\varepsilon)} c({s},\oeeX_{t_{i-1}},z) \rpm.
\end{aligned}
\end{equation} 
Note that for each $i$, the application $(\omega,x) \mapsto T_i(x,Y_i(x,\omega))$ is $(\mathcal{F}_{t_i} \otimes \mathcal{B}(\mathbb{R}),\mathcal{B}(\mathbb{R}))$-measurable as a composition of the measurable applications $T_i$ and $(\omega,x) \mapsto Y_i(x,\omega)$, since $Y_i(x)$ is a $\mathcal{F}_{t_i}$-measurable random variable and $x \mapsto Y_i(x,\omega)$ is a continuous function for all $\omega \in \Omega$. 
This ensures in particular that the sequence $(\oeeX_{t_i})_{0 \leq i \leq n}$ is $(\mathcal{F}_{t_i})_{0 \leq i \leq n}$-adapted.
Comparing the term-by-term construction in \eqref{eq:numerical_scheme} and \eqref{eq:def_x_chapeau}, it is straightforward to conclude that $\oeeX_{t_i}$ has the same law as $\oeX_{t_i}$ for each $0 \leq i \leq n$.

\subsubsection{Uniform moment bound for $\oeeX$} \label{ssec:uniform_moment_bound} Thanks to the Chasles relation, one can write for any $t \in [0,T]$:
\begin{equation} \label{eq:x_chapeau_eta}
\begin{aligned}
    \oeeX_{\eta(t)} = & \  X_0 + \int_0^{\eta(t)} a(\eta(s),\oeeX_{\eta(s)}) ds + \int_0^{\eta(t)} b(\eta(s),\oeeX_{\eta(s)}) dW_s \\ & + \int_0^{\eta(t)} \WEint c({s},\oeeX_{\eta(s^-)},z) \rpm + \sum_{i=0}^{\rho(t)} T_i(\oeeX_{t_{i-1}},Y_i(\oeeX_{t_{i-1}})).
\end{aligned}
\end{equation}
We want to estimate $\normPt{\sup_{0 \leq s \leq t}|\oeeX_{\eta(s)}|}$ with a Gronwall argument. While we can deal with the first terms of \eqref{eq:x_chapeau_eta} with the standard BDG and Kunita's BDG inequalities, the last one require a special attention. We now show that $(T_i(\oeeX_{t_{i-1}},Y_i(\oeeX_{t_{i-1}})))_{0 \leq i \leq n}$ is a sequence of martingale differences, allowing to apply the discrete BDG inequality mentioned in the preamble (Theorem \ref{thm:discrete_BDG}). 

From Section \ref{ssec:construction_schema}, we know that $(T_i(\oeeX_{t_{i-1}},Y_i(\oeeX_{t_{i-1}})))_{0 \leq i \leq n}$ is adapted to the filtration $(\mathcal{F}_{t_i})_{0 \leq i \leq n}$. For each $i \geq 1$, we now compute 
$\mathbb{E} [ T_i(\oeeX_{t_{i-1}},Y_i(\oeeX_{t_{i-1}}))
  | \mathcal{F}_{t_{i - 1}}]$ with the properties of the conditional expectation. 
Indeed, we know that:
\begin{itemize}
    \item the random variable $\oeeX_{t_{i-1}}$ is $\mathcal{F}_{t_{i-1}}$-measurable by construction,
    \item since the increments of $\widetilde{N}$ on $(t_{i-1},t_i]$ are independent of $\mathcal{F}_{t_{i-1}}$, the random variable $Y_i = (Y_i(x))_{x\in\mathbb{R}}$ is independent of $\mathcal{F}_{t_{i-1}}$,
    \item and the map $(x,f)\ni \mathbb{R}\times\mathbb{R}^\mathbb{R} \mapsto T_i(x,f(x))$ is $(\mathcal{B}(\mathbb{R}) \otimes \mathcal{B}(\mathbb{R}^\mathbb{R}))$-measurable, where $\mathcal{B}(\mathbb{R}^\mathbb{R})$ designate the $\sigma$-algebra generated by the cylindrical subsets of $\mathbb{R}^\mathbb{R}$, since it can be seen as the composition of $T_i$ with $(x,f) \mapsto (x,f(x))$. 
\end{itemize}
Using for example \cite[Proposition A.2{.}5]{lamberton2007} we deduce that 
\begin{align} \label{eq:contionnement_mesurable}
    \mathbb{E} \left[ T_i(\oeeX_{t_{i-1}},Y_i(\oeeX_{t_{i-1}}))
  \Big| \mathcal{F}_{t_{i - 1}} \right] &= v(\oeeX_{t_{i-1}}),
\end{align}
where for every $x\in \mathbb{R}$,  it holds that $v(x) := \EE[ T_i(x,Y_i(x)) ] = 0$,  since $T_i(x,Y_i(x))$ is a centred Gaussian random variable. This shows that for any $i\geq 1$, $ T_i(\oeeX_{t_{i-1}},Y_i(\oeeX_{t_{i-1}}))$ is centred conditionally to $\mathcal{F}_{t_{i-1}}$, hence the sequence $(T_i(\oeeX_{t_{i-1}},Y_i(\oeeX_{t_{i-1}})))_{0\leq i \leq n}$ is a sequence of martingale differences.

We may now use the discrete BDG inequality (Theorem \ref{thm:discrete_BDG})  to obtain:
\begin{equation}\label{eq:first_discrete_bdg}
\begin{aligned}
    \normP{\sup_{0 \leq s \leq t} \left|\sum_{i=0}^{\rho(s)} T_i(\oeeX_{t_{i-1}},Y_i(\oeeX_{t_{i-1}}))\right|} 
= & \normP{\sup_{k \in \{0,\ldots,\rho(t)\}} \left|\sum_{i=0}^{k} T_i(\oeeX_{t_{i-1}},Y_i(\oeeX_{t_{i-1}}))\right|} \\
 & \leq \kappa_p \left\|   \sum_{i = 0}^{\rho (t)} | T_i(\oeeX_{t_{i-1}},Y_i(\oeeX_{t_{i-1}})) |^2 \right\|_{L^{\frac{p}{2}}   (\Omega)}^{\frac{1}{2}}\\
    & \leq \kappa_p \left( \sum_{i = 0}^{\rho (t)} \|  T_i(\oeeX_{t_{i-1}},Y_i(\oeeX_{t_{i-1}})) \|_{L^p (\Omega)}^2 \right)
    ^{\frac{1}{2}}
\end{aligned}
\end{equation}
where $\kappa_p$ stands for the constant in the discrete BDG inequality.  Using integral Minkowski inequality along with standard BDG inequality, Kunita's BDG inequality and \eqref{eq:first_discrete_bdg}, we estimate $\normPt{\sup_{0 \leq s \leq t}|\oeeX_{\eta(s)}|}$ as follows:
\begin{equation}\label{eq:estimation_moment_epsilon}
\begin{aligned}
\normPt{\sup_{0 \leq s \leq t}|\oeeX_{\eta(s)}|}
\leqpt & \normP{X_0} + \int_0^{\eta(t)} \normP{a(\eta(s),\oeeX_{\eta(s)})} ds 
    + \left( \int_0^{\eta(t)} \normP{b(\eta(s),\oeeX_{\eta(s)})}^2 ds \right)^\frac{1}{2} \\
    & + \left( \int_0^{\eta(t)} \EE\left[ \left( \WEint |c(s,\oeeX_{\eta(s^{{-}})},z)|^{2} \lm \right) ^\frac{p}{2} \right] ds \right)^\frac{1}{p} \\
    & + \left( \int_0^{\eta(t)} \WEint \EE\left[ |c(s,\oeeX_{\eta(s^{{-}})},z)|^{p} \lm \right] ds \right)^\frac{1}{p} + \left( \sum_{i = 0}^{\rho (t)} \|  T_i(\oeeX_{t_{i-1}},Y_i(\oeeX_{t_{i-1}})) \|^2_{L^p (\Omega)} \right)
    ^{\frac{1}{2}}.
\end{aligned}
\end{equation}

Next, we provide an upper-bound for the $L^p$-norm of $T_i(\oeeX_{t_{i-1}},Y_i(\oeeX_{t_{i-1}}))$. For $i \geq 1$,  using 
again \cite[Proposition A.2{.}5]{lamberton2007}, it holds that
\begin{align} \label{eq:conditionnement}
     \mathbb{E} \left[ \left| T_i(\oeeX_{t_{i-1}},Y_i(\oeeX_{t_{i-1}})) \right|^{p}
  \Bigg| \mathcal{F}_{t_{i - 1}} \right] = w(\oeeX_{t_{i-1}}),
\end{align}
where $w(x) := \EE[|T_i(x,Y_i(x))|^p]$ for every $x\in \mathbb{R}$. 
Recalling that the $p^{\text{th}}$-moment of a centred Gaussian random variable with variance $\sigma$ is given by $K_p \sigma^{\frac{p}{2}}$, where $K_p = \frac{2^{p/2} \Gamma(\frac{p+1}{2})}{\sqrt{\pi}}$,\footnote{$\Gamma$ being the standard Gamma function.} we compute the  $p^{\text{th}}$-moment of the centred Gaussian variable $T_i(x,Y_i(x))$, whose variance is $\int_{t_{i - 1}}^{t_i} \int_{B(\varepsilon)} | c (s, y, z) |^2 \nu_s (d z) d s$:
\begin{align*}
    \mathbb{E} [| T_i(x,Y_i(x)) |^{p}] = K_p \ \left(\int_{t_{i - 1}}^{t_i} \int_{B(\varepsilon)} | c (s, x, z) |^2 \nu_s (d z) d s\right)^{p/2}.
\end{align*}
In addition,  by \ref{hyp:lipschitz}, $| c (s, x, z) |^2 \leq \overline{L}_c (s,z) ^ 2 (1 + | x |)^2$,
hence 
\begin{align*}
  w(x) \leq & K_p \ \left(\int_{t_{i-1}}^{t_i} \int_{B(\varepsilon)} | \overline{L}_c(s,z) |^2 \nu_s (d z) d s\right)^{p/2} \ (1+|x|)^p.
\end{align*}
Coming back to \eqref{eq:conditionnement}, we obtain
\begin{align} \label{eq:moment_of_G}
    \EE[|T_i(\oeeX_{t_{i-1}},Y_i(\oeeX_{t_{i-1}}))|^p] & \leqpt \left(\int_{t_{i-1}}^{t_i} \int_{B(\varepsilon)} | \overline{L}_c(s,z) |^2 \nu_s (d z) d s\right)^{p/2} \ (1+\EE[|\oeeX_{t_{i-1}}|]^p),
\end{align}
leading to
\begin{align*}
    \|  T_i(\oeeX_{t_{i-1}},Y_i(\oeeX_{t_{i-1}})) \|^2_{L^p (\Omega)} \leqpt
    &  \int_{t_{i-1}}^{t_i} \left( \Sint | \overline{L}_c(s,z) |^2 \nu_s (d z) \right) \left(1+\normPt{\oeeX_{\eta(s^{{-}})}}^2\right) d s. 
\end{align*}
Taking the sum from $i=0$ to $\rho(t)$ and using $t_{\rho(t)} = \eta(t)$, a Chasles relation gives
\begin{align*}
    \sum_{i = 0}^{\rho (t)} \|  T_i(\oeeX_{t_{i-1}},Y_i(\oeeX_{t_{i-1}})) \|^2_{L^p (\Omega)} \leqpt \int_{0}^{\eta(t)} \left( \Sint | \overline{L}_c(s,z) |^2 \nu_s (d z) \right) \left(1+\normPt{\oeeX_{\eta(s^{{-}})}}^2\right) d s.
\end{align*}
Finally, coming back to \eqref{eq:estimation_moment_epsilon} {and since $|\oeeX_{\eta(s^{-})}| \leq \sup_{0 \leq u \leq s}|\oeeX_{\eta(u^{-})}| \leq \sup_{0 \leq u \leq s}|\oeeX_{\eta(u)}| $}, we end up with 
\begin{equation} \label{eq:extended_extended_gronwall}
\begin{aligned}
\normPt{\sup_{0 \leq s \leq t}|\oeeX_{\eta(s)}|} \leqpt & \ 1 + \int_0^{\eta(t)} \normPt{\sup_{0 \leq u \leq s}|\oeeX_{\eta(u)}|}  ds 
    + \left( \int_0^{\eta(t)} \normPt{\sup_{0 \leq u \leq s}|\oeeX_{\eta(u)}|}^2 ds \right)^\frac{1}{2} \\ 
    & + \left( \int_0^{\eta(t)} \psi_p(s) \normPt{\sup_{0 \leq u \leq s}|\oeeX_{\eta(u)}|}^p ds \right)^\frac{1}{p} \\
    & + \left( \int_{0}^{\eta(t)} \left( \Sint | \overline{L}_c(s,z) |^2 \nu_s (d z) \right) \normPt{\sup_{0 \leq u \leq s}|\oeeX_{\eta(u)}|}^2 d s \right)^\frac{1}{2}.
\end{aligned}
\end{equation}
Using the extended Gronwall lemma \ref{lem:generalized_gronwall},  (with   $ g(t) = \ell(t) = 1$, $\tfrac{1}{2} h(t) = 1 \vee \Sint | \overline{L}_c(s,z) |^2 \nu_s (d z)$, $k(t) =\psi_p(t)$), we obtain a bound  which does not depend on $\varepsilon$ anymore. In particular, we have 
\[\sup_{\varepsilon  \in [0,\varepsilon^*)} \normPt{\sup_{0 \leq s \leq t}|\oeeX_{\eta(s)}|} \leq \overline{m}(p,T)\]
for some constant $\overline{m}(p,T) < +\infty$, which prove the first statement \eqref{eq:moment_bound_scheme} of Theorem \ref{theoreme_maurer}.

\subsubsection{Upper-bound for $ \|\sup_{0\leq s \leq t} |\wtX_{\eta(s)} - \oeeX_{\eta(s)}|\ \|_{L^p(\Omega)}$} \label{proof:convergence_eemscheme}

We introduce the sequence $(D_i^\varepsilon)_{0 \leq i \leq n}$ defined, for every $0 \leq i \leq n$, by
\[D_i^\varepsilon = \int_{t_{i-1}}^{t_i} \Eint c({s},\oeeX_{t_{i-1}},z) \rpm - T_i(\oeeX_{t_{i-1}},Y_i(\oeeX_{t_{i-1}})), \quad D_0^\varepsilon = 0,\]
that allows us to write, for $t \in [0,T]$:
\begin{equation}
\begin{aligned}
\oeeX_{\eta(t)} - \wtX_{\eta(t)} = & \int_0^{\eta(t)} \left(a(\eta(s^{{-}}),\oeeX_{\eta(s^{{-}})}) - a(\eta(s^{{-}}),\wtX_{\eta(s^{{-}})})\right) ds \  + \int_0^{\eta(t)} \left(b(\eta(s^{{-}}),\oeeX_{\eta(s^{{-}})}) - b(\eta(s^{{-}}),\wtX_{\eta(s^{{-}})})\right) dW_s \\ & + \int_0^{\eta(t)} \WEint \left(c(s,\oeeX_{\eta(s^{{-}})},z) - c({s},\wtX_{\eta(s^-)},z)\right) \rpm + \  \sum_{i=0}^{\rho(t)} D_i^\varepsilon.
\end{aligned}
\end{equation}
We  show that $(D_i^\varepsilon)_{0 \leq i \leq n}$ is a sequence of $(\mathcal{F}_{t_i})_{0 \leq i \leq n}$-martingale differences with the same arguments used in Section \ref{ssec:uniform_moment_bound} for $(T_i(\oeeX_{t_{i-1}},Y_i(\oeeX_{t_{i-1}})))_{0 \leq i \leq n}$.  Then, we set $\mathcal{E}(t,\varepsilon) = \normP{\sup_{0\leq s \leq t} |\wtX_{\eta(s)} - \oeeX_{\eta(s)}|}$ and use Minkowski, BDG, Kunita's BDG and the discrete BDG inequalities to obtain, similarly to Section \ref{ssec:uniform_moment_bound}:
\begin{align} \label{eq:majoration_E_t_epsilon}
    \mathcal{E}(t,\varepsilon) \leqpt \int_0^{\eta(t)} \mathcal{E}(s,\varepsilon) ds + \left( \int_0^{\eta(t)} \mathcal{E}(s,\varepsilon)^2 \right)^\frac{1}{2} + \left( \int_0^{\eta(t)} \psi_p(s) \mathcal{E}(s,\varepsilon)^p \right)^\frac{1}{p} + \left( \sum_{i=0}^{\rho(t)} \normP{D_i^\varepsilon}^2 \right)^\frac{1}{2}.
\end{align}
Let us detail the estimation of $\normP{D_i^\varepsilon}$ for $i\geq 1$. Conditioning by $\mathcal{F}_{t_{i-1}}$ and using \eqref{eq:def_gj}, it holds that 
\begin{align} \label{eq:dj_norme}
    \normP{D_i^\varepsilon} = \EE[\phi(\oeeX_{t_{i-1}})],
\end{align}
where for $\mathbb{Q}_{i}$-almost every $x \in \mathbb{R}$,
\begin{align} \label{eq:calcul_phi_x}
    \phi(x) & =  \mathcal{W}_{p} \Bigg( \law{\int_{t_{i - 1}}^{t_i} \int_{B(\varepsilon)} c (s, x, z) \rpm}, \mathcal{N} \left( 0, \int_{t_{i - i}}^{t_i} \int_{B(\varepsilon)} | c ({s}, x, z) |^2 \nu_s (d z) d s \right) \Bigg).
\end{align}
Then we apply Proposition \ref{prop:wasserstein} with $F(s,z) = c(s,x,z)$, $t_0 = t_{i-1}$ and $t=t_i$ to obtain:
\begin{align} \label{eq:Berry--Essenation_phi}
    \phi(x) \leq \left(\mathcal{A}(p)  \frac{\int_{t_{i-1}}^{t_i} \int_{B(\varepsilon)} |
    c(s,x,z) |^{p + 2} \nu_s (d z) d s}{\int_{t_{i-1}}^{t_i} \int_{B(\varepsilon)}  | c(s,x,z) |^2 \nu_s (d z) d s} \right)^{1/p}.
\end{align}
We recall that by construction, $\mathbb{Q}_i = \tlaw{(\oeeX_{t_{i-1}})}$. Since the equality \eqref{eq:calcul_phi_x} stands true for $\mathbb{Q}_i$-almost every $x$, we deduce from \eqref{eq:majoration_phi} that 
\begin{align*}
     \EE[\phi(\oeeX_{t_{i-1}})] \leq \mathcal{A}(p) ^{1/p}  \ \mathbb{E} \left[ \left( \frac{\int_{t_{i-1}}^{t_i} \int_{B(\varepsilon)}  |
    c(s,\oeeX_{t_{i-1}},z) |^{p + 2} \nu_s (d z) d s}{\int_{t_{i-1}}^{t_i} \int_{B(\varepsilon)}  |  c(s,\oeeX_{t_{i-1}},z) |^2 \nu_s (d z) d s} \right)^{1/p} \right].
\end{align*}
Hence
\begin{align} \label{eq:lp_control_of_D}
     \left( \sum_{i=0}^{\rho(t)} \normP{D_i^\varepsilon}^2 \right)^\frac{1}{2} \leqpt & \left( \sum_{i=1}^{n} \mathbb{E} \left[ \left( \frac{\int_{t_{i-1}}^{t_i} \int_{B(\varepsilon)}  |
    c(s,\oeeX_{t_{i-1}},z) |^{p + 2} \nu_s (d z) d s}{\int_{t_{i-1}}^{t_i} \int_{B(\varepsilon)}  |  c(s,\oeeX_{t_{i-1}},z) |^2 \nu_s (d z) d s} \right)^{1/p} \right]^2 \right)^\frac{1}{2}.
\end{align}
To conclude, we combine \eqref{eq:majoration_E_t_epsilon} and \eqref{eq:lp_control_of_D} (which right hand side is equal to $\delta_p^n(\varepsilon)$ according to \eqref{eq:Delta_n_epsilon}), and use Gronwall Lemma \ref{lem:generalized_gronwall} to get: 
\begin{align} \label{eq:erreur_schema}
    \normP{\sup_{0\leq s \leq t} |\wtX_{\eta(s)} - \oeeX_{\eta(s)}|} \leqpt  \delta_p^n(\varepsilon).
\end{align}
Finally, given \eqref{eq:erreur_schema} and the result of Lemma \ref{euler_rate}, we obtain the second statement of Theorem \ref{theoreme_maurer}.

\subsubsection{Convergence of $\delta_p^n(\varepsilon)$ with $\varepsilon \rightarrow 0$} \label{ssec:convergence_epsilon}

Fix $n$ and  $i \in \{1,\ldots, n\}$. By \ref{hyp:lipschitz}, we get the upper-bound 
\begin{equation*}
\sup_{s\in[t_{i-1},t_i)}|c(s,\oeeX_{t_{i-1}},z)|^{p} \leq \sup_{s\in[0,T]} \overline{L}_c(s,z)^p (1+|\oeeX_{t_{i-1}}|)^p.
\end{equation*}

Let $\eta > 0$. According to  \eqref{eq:c_tend_vers_zero}, there exists an $\varepsilon_0 > 0$ such that for any $|z| \leq \varepsilon_0$, we have $\sup_{s\in[0,T]} \overline{L}_c(s,z) \leq \eta$.
Then taking $\varepsilon \leq \varepsilon_0$, one can split the term  $|c(s,\oeeX_{t_{i-1}},z)|^{p+2}$ into two parts and obtain:
\begin{align*}
    \int_{t_{i-1}}^{t_i} \int_{B(\varepsilon)} |c(s,\oeeX_{t_{i-1}},z)|^{p+2} \rpc = & \int_{t_{i-1}}^{t_i} \int_{B(\varepsilon)} |c(s,\oeeX_{t_{i-1}},z)|^{p}  |c(s,\oeeX_{t_{i-1}},z)|^2 \rpc \\ & \leq \eta^p (1+|\oeeX_{t_{i-1}}|)^p \int_{t_{i-1}}^{t_i} \int_{B(\varepsilon)} |c(s,\oeeX_{t_{i-1}},z)|^2 \rpc.
\end{align*}
Hence we have the simplification
\begin{align*}
   \EE\left[ \left( \frac{\int_{t_{i-1}}^{t_i} \int_{B(\varepsilon)} |c(s,\oeeX_{t_{i-1}},z)|^{p+2} \rpc}{\int_{t_{i-1}}^{t_i} \int_{B(\varepsilon)} |c(s,\oeeX_{t_{i-1}},z)|^{2} \rpc} \right)^{\frac{1}{p}} \right] \leq \eta \  \left(1+\EE\left[\sup_{0\leq i\leq n} |\oeeX_{t_i}|\right]\right).
\end{align*}
Now thanks to the moment upper-bound  \eqref{eq:moment_bound_scheme} along with Hölder inequality it holds that
\begin{align*}
    \sup_{\varepsilon \in [0,\varepsilon_0]} \EE\left[\sup_{0\leq i\leq n} |\oeeX_{t_i}|\right] \leq {\overline{m}(p,T)}^{\frac{1}{p}},
\end{align*}
therefore we obtain that for any $\varepsilon \leq \varepsilon_0$, 
\begin{align*}
    \delta_p^n(\varepsilon) \leq \eta \  \sqrt{n}(1+{\overline{m}(p,T)}^{\frac{1}{p}}),
\end{align*}
which shows that $\delta_p^n(\varepsilon) \rightarrow 0$, when $\varepsilon$ goes to zero by definition of the limit. 
This concludes the proof of Theorem~\ref{theoreme_maurer}. \qed

\subsection{Proof of Corollary \ref{cor:mainresult}.}
\label{sec:proof_coro_main}

Using the sublinearity condition \eqref{eq:bound_on_the_jump_coefficient}, we directly  simplify 
\begin{align*}
   \EE\left[ \left( \frac{\int_{t_{i-1}}^{t_i} \int_{B(\varepsilon)} |c(s,\oeeX_{t_{i-1}},z)|^{p+2} \rpc}{\int_{t_{i-1}}^{t_i} \int_{B(\varepsilon)} |c(s,\oeeX_{t_{i-1}},z)|^{2} \rpc} \right)^{\frac{1}{p}} \right] \leq \varepsilon \  \left(1+\EE\left[\sup_{0\leq i\leq n} |\oeeX_{t_i}|\right]\right),
\end{align*}
leading to $\delta_p^n(\varepsilon) \leqpt \varepsilon \sqrt{n}$.
With the $L^2$-integrability of $\psi_p$, Theorem \ref{theoreme_maurer} applies  with $\zeta = 1$, and we get
\begin{align*}
 \left\| \sup_{t\in[0,T]} | X_{\eta(t)} - \oeeX_{\eta(t)}| \right\|_{L^p(\Omega)}   \leqpt   {n^{ - \left(\gamma \wedge \left( \frac{1}{p} \right) \right) }}   +  \varepsilon \sqrt{n}
\end{align*} 
for $\varepsilon \in [0,\varepsilon_0)$.
Finally, assuming in addition that $\varepsilon \leq  {n^{ - \frac{1}{2} - \left(\gamma \wedge \frac{1}{p}\right)}}$, we obtain the announced convergence rate
\begin{align*}
    \left\| \sup_{t\in[0,T]} | X_{\eta(t)} - \oeeX_{\eta(t)}| \right\|_{L^p(\Omega)}   \leqpt  {n^{ - \left(\gamma \wedge \left( \frac{1}{p} \right) \right) }} .
\end{align*}

\subsection{Proof of Theorem \ref{thm:noapprox}} \label{proof:thmnoapprox}

Since $\ueX$ can be seen as a simpler jump-SDE than $X$ itself, the proof that $\sup_{\varepsilon > 0}  \EE[\sup_{0 \leq i \leq n-1} |\ueX_{t_i}|^p]$ is finite poses no difficulties. Using again the Euler-Peano scheme as a pivot, we know from Lemma \ref{lem_peano} that 
\begin{align*}
    \Big\| \underset{t \in [0, T]}{\sup} | X_t - \widetilde{X}_t
    | \Big\|_{L^{p}(\Omega)} & \leqpt  n^{-
    \left\{\gamma \wedge \left( \frac{2 \zeta }{p (1+\zeta)} \right) \right\}}.
  \end{align*}
It {remains} to find an upper-bound for $ \|\sup_{0\leq t \leq T} |\wtX_{\eta(t)} - \ueX_{\eta(t)}|\ \|_{L^p(\Omega)}$. The arguments used are similar to the ones from Section \ref{proof:convergence_eemscheme}, instead that we now have to estimate the $L^p$-norm of the small jump integral itself, instead of the difference $D_i^\varepsilon$ between the latter and its Gaussian approximation.
For $t \in [0,T]$, one may write:
\begin{equation}
\begin{aligned}
\ueX_{\eta(t)} - \wtX_{\eta(t)} = & \int_0^{\eta(t)} \left(a(\eta(s^{{-}}),\ueX_{\eta(s^{{-}})}) - a(\eta(s^{{-}}),\wtX_{\eta(s^{{-}})})\right) ds \  + \int_0^{\eta(t)} \left(b(\eta(s^{{-}}),\ueX_{\eta(s^{{-}})}) - b(\eta(s^{{-}}),\wtX_{\eta(s^{{-}})})\right) dW_s \\ & + \int_0^{\eta(t)} \WEint \left(c({s},\ueX_{\eta(s^-)},z) - c(s,\wtX_{\eta(s^-)},z)\right) \rpm + \int_0^{\eta(t)} \Eint c({s},{\wtX}_{\eta(s^-)},z) \rpm.
\end{aligned}
\end{equation}
For  $\mathcal{E}(t,\varepsilon) = \normP{\sup_{0\leq s \leq t} |\wtX_{\eta(s)} - \ueX_{\eta(s)}|}$, the standard inequalities induce

\begin{equation}\label{eq:majoration_E_t_epsilon_noapprox}
\begin{aligned} 
    \mathcal{E}(t,\varepsilon) \leqpt & \int_0^{\eta(t)} \mathcal{E}(s,\varepsilon) ds + \left( \int_0^{\eta(t)} \mathcal{E}(s,\varepsilon)^2 \right)^\frac{1}{2} + \left( \int_0^{\eta(t)} \psi_p(s) \mathcal{E}(s,\varepsilon)^p \right)^\frac{1}{p} \\ & + \left\|\int_0^{\eta(t)} \Eint c({s},{\wtX}_{\eta(s^-)},z) \rpm \right\|_{L^{p}(\Omega)}.
\end{aligned}
\end{equation}
As said previously, we estimate the $L^p$-norm of the small jump integral, using Kunita's BDG inequality:
\begin{align*} 
    \left\|\int_0^{\eta(t)} \Eint c({s},{\wtX}_{\eta(s^-)},z) \rpm \right\|_{L^{p}(\Omega)}^p  = & \EE\left[\left(\int_{0}^{\eta(t)} \Eint c({s},{\wtX}_{t_{i-1}},z) \rpm \right)^p\right]\\
    \leq & \ckunpT \int_{0}^{\eta(t)}  \mathbb{E} \left[  \left( \int_{B(\varepsilon)} | c(s,{\wtX}_{t_{i-1}},z)|^2 \nu_s (d z) \right)^{p / 2} \right] ds \\ & + \ckunp \mathbb{E} \left[ \int_{0}^{\eta(t)}     \int_{B(\varepsilon)}|c(s,{\wtX}_{t_{i-1}},z)|^p \nu_s (d z) d s \right].
\end{align*}
Using \ref{hyp:lipschitz}, we have the upper bounds
\begin{align*}
    \mathbb{E} \left[  \left( \int_{B(\varepsilon)} | c(s,{\wtX}_{t_{i-1}},z)|^2 \nu_s (d z) \right)^{p / 2} \right] & \leq \left( \int_{B(\varepsilon)} \overline{L}_c^2(s,z) \nu_s(dz) \right)^{p / 2} \EE[(1+|{\wtX}_{t_{i-1}}|)^p], \\
    \EE \left[ \int_{0}^{\eta(t)}     \int_{B(\varepsilon)}|c(s,{\wtX}_{t_{i-1}},z)|^p \nu_s (d z) d s \right] & \leq \left( \int_{0}^{\eta(t)} \int_{B(\varepsilon)} \overline{L}_c^p(s,z) \nu_s(dz) ds \right) \EE[(1+|{\wtX}_{t_{i-1}}|)^p].
\end{align*}
Considering the $\varepsilon$-truncated $\psi_p$ function 
\begin{align*}
    \psi_{p,\varepsilon}(s) =  \left( \int_{B(\varepsilon)} \overline{L}_c^2(s,z) \nu_s(dz) \right)^{p / 2} +  \int_{B(\varepsilon)} \overline{L}_c^p(s,z) \nu_s(dz), 
\end{align*}
we obtain 
\begin{align*}
    \left\|\int_0^{\eta(t)} \Eint c({s},{\wtX}_{\eta(s^-)},z) \rpm \right\|_{L^{p}(\Omega)}^p \leqpt (1+ \EE[\sup_{0 \leq i \leq n-1} |{\wtX}_{t_i}|^p]) \int_0^{\eta(t)} \psi_{p,\varepsilon}(s) ds, 
\end{align*}
hence \eqref{eq:majoration_E_t_epsilon_noapprox} reduces to
\begin{align*}
    \mathcal{E}(t,\varepsilon) \leqpt & \int_0^{\eta(t)} \mathcal{E}(s,\varepsilon) ds + \left( \int_0^{\eta(t)} \mathcal{E}(s,\varepsilon)^2 \right)^\frac{1}{2} + \left( \int_0^{\eta(t)} \psi_p(s) \mathcal{E}(s,\varepsilon)^p \right)^\frac{1}{p} + \left( \int_0^{\eta(t)} \psi_{p,\varepsilon}(s) ds\right)^\frac{1}{p}.
\end{align*}
We end the proof using the Gronwall Lemma \ref{lem:generalized_gronwall}, leading to
\begin{align*}
\normP{\sup_{0\leq s \leq t} |\wtX_{\eta(s)} - \ueX_{\eta(s)}|} \leqpt \left( \int_0^{\eta(t)} \psi_{p,\varepsilon}(s) ds\right)^\frac{1}{p}.
\end{align*}

\section{Weak convergence}\label{sec:weak_convergence}

In this section, we give a rate of convergence bound for the weak error of the $\varepsilon$-EM scheme. It is well established that for a classical Euler scheme, the weak error analysis can be derived from the regularity of the stochastic flow of the SDE. Here, we consider the Lévy flow $(X_\theta(t,x), \theta\geq t)$,  solution of \eqref{eq:intro_SDEs} and starting at $x$ at time $t$. 
Once again, we make use of the  work of  Breton and Privault \cite{breton} that gives the existence in $L^p$ of the $k^{\text{th}}$ derivative processes  ${\partial_x^k X_\theta}(t,x)$ of this  Lévy flow. 

Working in a framework of sufficiently regular coefficients for the SDE \eqref{eq:intro_SDEs}, it is then straightforward to deduce the wellposedness of the Cauchy problem associated with the Kolmogorov PDE involving the integro-differential generator $\mathcal{L}_{(\cdot)}$ of $X$.  Precisely, for some given function $\varPhi$,  we consider  the solution $u(t,x)$ to the  Kolmogorov PDE
\begin{equation}\label{eq:kolmogorov_pde}
\begin{aligned}
\left\{
    \begin{array}{ll}
         \frac{\partial u}{\partial t} +\mathcal{L}_{(t)} u  = 0, \quad (t,x)\in[0,T)\times\mathbb{R}, \\
         u (T, x)  = \varPhi (x), \quad x \in \mathbb{R},
    \end{array}
\right.
\end{aligned}
\end{equation}
where 
\begin{align*}
 \mathcal{L}_{(t)}f (x) & = a (t,x)  \frac{\partial f}{\partial x} (x) + \frac{1}{2} b^2 (t,x)   \frac{\partial^2 f}{\partial x^2} (x) +
\int_{\mathbb{R}} \left\{ f (x + c (t, x, z)) - f (x) - c (t, x, z)  \frac{\partial f}{\partial x} (x) \right\} \nu_t (d z).
\end{align*}

Under \ref{hyp:lipschitz}, the map $(t,x) \mapsto (X_\theta(t,x), \theta\geq t)$  
is a two-parameters Lévy flow (in the sense of Applebaum \cite{apple}, from immediate extension of Theorem 6.4.2 to the time-non-homogeneous case).
The existence and integrability of the flow derivatives are more demanding in terms of smoothness. We summarise below the requirements needed from \cite{breton}, combined with the additional time-regularity and integrability needed for the weak error analysis. 

\begin{hyp}{\hspace{-0.1cm}\Blue{\bf  \textrm Stochastic Flow regularity condition. (\ref{hyp:smooth}).}}
\makeatletter\def\@currentlabel{ {\bf\textrm H}$_{\mbox{\scriptsize\bf\textrm{S\_Flow}}}$}\makeatother
\label{hyp:smooth}
\begin{itemize}
    \item[\textit{\textbf{(i)}}] There exists a Lévy measure $\mu$ on $\mathbb{R}$,  such that for any $t\in[0,T]$ and $A\in\mathcal{B}(\mathbb{R})$, $\nu_t(A) \leq \mu(A)$. 

    \item[\textit{\textbf{(ii)}}]  For every $t\in\mathbb{R}_{+}$, the functions $a(t,\cdot)$, $b(t,\cdot)$ and $c(t,\cdot,\cdot)$ are $\mathcal{C}^4$-differentiable and there exists some positive constant $C$ such that
\begin{align*}
    \sup_{(t,x) \in [0,T] \times \RR}\left( \left| \frac{\partial^k a}{\partial x^k}(t,x) \right| + \left| \frac{\partial^k b}{\partial x^k}(t,x) \right|  \right) \ +  \  \sup_{(t,x,z) \in [0,T] \times \RR^2}\left| \frac{\partial^{k+l} c}{\partial x^k \partial z^l}(t,x,z) \right| \leq C,
\end{align*}
for all $k,l$ in $\{1,\ldots,4\}$ with $1 \leq k+l \leq 4$. Moreover,  there exists a function $g$ in  $\bigcap_{q\geq 2} L^q(\mathbb{R},\mu)$ such that for any $z\in \RR$,
\begin{align*}
    \sup_{(t,x) \in [0,T] \times \RR} \left(\left| \frac{\partial^{k} c}{\partial x^k}(t,x,z) \right| +  \left| \frac{\partial^{k} c^2}{\partial x^k}(t,x,z) \right|\right)  \leq C  \ g(z), \ \mbox{ for }k=1,\ldots,4.
\end{align*}
\item[\textit{\textbf{(iii)}}]  For every $(x,z) \in\RR\times \RR$, the functions $a(\cdot,x)$, $b(\cdot,x)$, $c(\cdot,x,z)$ 
 are a.e. time-differentiable on $[0,T]$,  and there exists a function $\tilde{g}$ in $\bigcap_{q\geq 2} L^q(\mathbb{R},\mu)$ such that for any $z\in \RR$, 
\begin{align*}
\sup_{(t,x) \in [0,T] \times \RR} \left( \left| \frac{\partial a}{\partial t}(t,x) \right| \ 
+ \  \left| \frac{\partial b}{\partial t}(t,x) \right| \right) \leq \tilde{C},    
\qquad \sup_{(t,x) \in [0,T] \times \RR} \left(\left| \frac{\partial c}{\partial t}(t,x,z) \right| 
+ \left| \frac{\partial c^2}{\partial t}(t,x,z) \right|\right) \leq \tilde{C}  \  \tilde{g} (z). 
\end{align*}

\item[\textit{\textbf{(iv)}}]  As in \ref{hyp:peano}, $L_{a,b}(\cdot) \in L^{\infty}([0,T])$.
    
\item[\textit{\textbf{(v)}}] The sub-linearity condition \eqref{eq:bound_on_the_jump_coefficient} holds for the jump coefficient $c$. Moreover, for all $z \in \RR$ and $k \in \{1,2\},$
\begin{align*}
    \sup_{(t,x) \in [0,T] \times \RR} \left( \left| \frac{\partial^k c}{\partial x^k}(t,x,z) \right| \ 
+ \  \left| \frac{\partial c}{\partial t}(t,x,z) \right| \right) \leqpt |z|,  \quad \sup_{(t,x) \in [0,T] \times \RR} \left( \left| \frac{\partial^k c^2}{\partial x^k}(t,x,z) \right| \ 
+ \  \left| \frac{\partial c^2}{\partial t}(t,x,z) \right| \right) \leqpt |z|^2.
\end{align*}
\end{itemize}

\end{hyp}
In the hypotheses above, \text{\it{(i)}} and \text{\it{(ii)}} are reproduced from \cite{breton}, which  themselves come from \cite{bichteler-etal-87},  providing sufficient conditions for the  flow derivatives of the exact process up to the order 4 to exist and belong to $L^p(\Omega)$ uniformly in $(x, t) \in\mathbb{R}\times [0, T ]$, for any $p\geq 2$. Note that \text{\it{(i)}} is not excessively restrictive, since it covers for example the case of additive Markov processes (see Theorem 9.8 in \cite{sato}). In addition we require \text{\it{(iii)}} to the control the time derivative of the coefficients. This allows the two successive applications of the Itô formula used in the standard analysis of the weak error, based on the Kolmogorov PDE argument. We also impose the second part of \text{\it{(v)}} for compatibility reasons with \eqref{eq:bound_on_the_jump_coefficient}.

Assuming that $(\Omega,\mathcal{F},\mathbb{P},(\mathcal{F}_t)_{t\geq0})$ is rich enough to contain another Brownian motion $B$, independent of $W$ and adapted to $\mathcal{F}$, we introduce the time continuous version of the $\varepsilon$-EM scheme as 
\begin{equation}\label{eq:numerical_scheme_continous_times}
    \begin{aligned} 
    \oeX_{\theta} = & \oeX_{\eta(\theta)} +    \int_{\eta(\theta)}^\theta  a_\varepsilon (\eta(s), s, \oeX_{\eta(s)})  ds  + \int_{\eta(\theta)}^\theta  b (\eta(s),  \oeX_{\eta(s)}) dW_s \\ 
&   + \int_{\eta(\theta)}^\theta    d_{\varepsilon} (s, \oeX_{\eta(s)})   dB_s +  \sum_{j=N^\varepsilon(\eta(\theta))+1}^{N^\varepsilon(\theta)} c(T^\varepsilon(j),\oeX_{\eta(\theta^{{-}})},Z^\varepsilon(j)).
    \end{aligned}
\end{equation}
where for all $\varepsilon >0$, we  define  the $\varepsilon$-corrected coefficients 
\begin{equation}\label{eq:epsilon_corrected_coeffs}
\begin{aligned}
    a_\varepsilon(s,t,x):= a(s,x) - \WEint c(t,x,z) \nu_s(dz), \quad & \quad 
    d_\varepsilon(t,x):= \left( \Eint |c(t,x,z)|^2 \nu_t(dz) \right)^{\frac{1}{2}}.
\end{aligned}
\end{equation}
The process $\oeX$ also writes as the unique solution of 
\begin{equation}\label{eq:numerical_scheme_as_process}
\begin{aligned} 
\oeX_{\theta} = & X_0 + \int_{0}^\theta  a_\varepsilon (\eta(s), s, \oeX_{\eta(s)})  ds  + \int_{0}^\theta  b (\eta(s),  \oeX_{\eta(\theta)}) dW_s  \\ & +
 \int_{0}^\theta    d_{\varepsilon} (s, \oeX_{\eta(s)})   dB_s  +  \int_0^\theta \WEint c({s},\oeX_{\eta(s^{-})},z) N(ds,dz).
\end{aligned}
\end{equation}
In the same manner, the process $\ueX$, as the $\varepsilon$-EM scheme without substitute,  also writes as the unique solution of 
\begin{equation}\label{eq:numerical_scheme_as_process_ws}
\begin{aligned} 
\ueX_{\theta} = & X_0 + \int_{0}^\theta  a_\varepsilon (\eta(s), s, \ueX_{\eta(s)})  ds  + \int_{0}^\theta  b (\eta(s),  \ueX_{\eta(s)}) dW_s   +  \int_0^\theta \WEint c({s},\ueX_{\eta(s^{-})},z) N(ds,dz).
\end{aligned}
\end{equation}
With the use of standard inequalities (and Gronwall argument)  in the right-hand side of \eqref{eq:numerical_scheme_as_process} for  $\oeX$ (resp. \eqref{eq:numerical_scheme_as_process_ws} for $\ueX$),  we easily obtain moment bounds for the continuous schemes on the same condition than for $X$: 
\begin{lemma} \label{lem:uniform_bound_continuous_scheme} 
   Assume that \ref{hyp:lipschitz} holds for some $p\geq 2$. Then, 
   \begin{equation*}
       \sup_{ \varepsilon > 0  } \EE\left[\sup_{\theta\in[0,T]} \left( |\oeX_\theta|^p + |\ueX_\theta|^p \right) \right]  <  +\infty.
   \end{equation*}
\end{lemma}
\medskip

Our convergence result is as follows.
\begin{theorem}\label{thm:weak_error_bound}
For some $p\geq {3}$, assume \ref{hyp:lipschitz} and \ref{hyp:smooth}.
Let $\beta$  the Blumenthal--Getoor index of the indivisible distribution characterised by the Lévy measure $\mu$. Let $\varPhi$ be a function in $\mathcal{C}^4(\mathbb{R})$ such that 
\begin{equation}\label{eq:hypo_PHI}
    \Big| \frac{\partial^k \varPhi}{\partial x^k}  \Big| (x) \leq C (1+|x|^q),\  \mbox{for every $k \in \{0,1,\ldots,4\}$, with { $0\leq  q\leq (p-4)^+$}}. 
\end{equation}
\begin{itemize}
\item[\textit{\textbf{(i)}}]
For a given integer $n$ and  jump threshold $\varepsilon > 0$,   we have the following weak error upper-bound for the $\varepsilon$-Euler-Maruyama scheme with Gaussian substitute: 
\begin{align}\label{eq:weak_error_bound_with_substit}
| \EE[\varPhi(X_T)] - \EE[\varPhi(\oeX_T)] | \leqpt \frac{T}{n}  + \  \varepsilon^{3-\beta^+}. 
\end{align}

\item[\textit{\textbf{(ii)}}] In the same setting, the weak error upper-bound for the  $\varepsilon$-Euler-Maruyama scheme without substitute is 
\begin{align}\label{eq:weak_error_bound_without_substit}
| \EE[\varPhi(X_T)] - \EE[\varPhi(\ueX_T)] | \leqpt \frac{T}{n}  + \  \varepsilon^{2-\beta^+}. 
\end{align}
\end{itemize}
\end{theorem}
In the above inequality, the term $\varepsilon^{m-\beta^+}$ stands for $\sup_{\theta>0} \{ \varepsilon^{m-\theta}; \int_{|z|<1} |z|^{\theta} \mu(dz) < +\infty \}$. In particular, when $\mu$ is a truncated  alpha stable distribution, this term is   $\varepsilon^{m-\beta}$.

\begin{proof} {\bf Wellposedness and regularity of the Kolmogorov PDE \eqref{eq:kolmogorov_pde}. }
Let $u(t,x) \coloneqq \EE[\varPhi(X_T(t,x))]$. Applying \cite[Proposition 4.1]{breton} and \cite[Theorem 5.1]{breton}, we derive that for $k \in \{1,\ldots,4\}$, the $k^{\text{th}}$ derivative of the Lévy flow process $\partial^k_x X_\cdot(t,x)$ exists in the $L^p(\Omega)$-sense, and satisfy
\begin{align*}
    \sup_{x \in \mathbb{R}} \EE[ \sup_{\theta \in(t,T]}|\partial^k_x X_\theta(t,x)|^m ] < + \infty, 
\end{align*}
for any $m \geq 2$. Then, with a direct extension of \cite[Theorems 5.4 and 5.5]{friedman}, we deduce that $u$ admits four continuous spatial derivatives.   In particular, since $u$ is $C^2$-differentiable and $X$ is a Markov process w.r.t. the filtration $(\mathcal{F}_t)$, one can straightforwardly extend the proof of \cite[Theorem 6.1]{friedman}, ensuring that $\partial_t u(t,x)$ exists everywhere,  and hence that $u(t,x)$ verifies \eqref{eq:kolmogorov_pde}.  Moreover, $u(t,x)$ satisfies  
\begin{align}\label{eq:polynomial_growth_du}
    |\partial^k_x u(t,x)| + |\partial_t\partial^{\lfloor k/2\rfloor}_x u(t,x)| \leqpt 1 +  \EE[|X_T(t,x)|^{2q}]^{1/2} \leqpt (1+|x|^{q}), \quad k\leq 4. 
\end{align}

\paragraph{The weak error related to $\oeX$.}
Observing that $\EE[u(0,\oeX_0)] = \EE[\varPhi(X_T(0,X_0))] = \EE[\varPhi(X_T)]$ and $ \EE[u(T,\oeX_T)] = \EE[\varPhi(\oeX_T)]$, we can write 
\begin{align*}
  \mathbb{E} [\varPhi (\oeX_T) - \varPhi (X_T )] & =  \mathbb{E} [u
  (T,\oeX_T) - u (0,\oeX_0)].
\end{align*}
Applying the Itô formula to the function $(t,x) \mapsto u(t,x)$ and the process $(\oeX_t)_{t\in[0,T]}$, we derive that
\begin{equation}\label{eq:ito_1}
\begin{aligned} 
  \mathbb{E} [\varPhi (\oeX_T) - \varPhi (X_T )] = & \int_0^T \EE \left[ \frac{\partial
  u}{\partial t} (t,\oeX_t) + \overline{\mathcal{L}}_{\varepsilon, t, \oeX_{\eta(t^{{-}})}} u(t,\oeX_t) \right] dt + \EE \left[ \int_0^T b(\eta(t),\oeX_\eta(t)) \frac{\partial u}{\partial x}(t,\oeX_t) d W_t \right] 
  \\ & + \EE \left[ \int_0^T d_\varepsilon(t,\oeX_{\eta(t^{{-}})}) \frac{\partial u}{\partial x}(t,\oeX_t) d B_t \right], 
\end{aligned}
\end{equation}
where for $(t,y) \in [0,T]\times\RR$, the differential operator $\overline{\mathcal{L}}_{(\varepsilon, t, y)}$ is the frozen infinitesimal generator of $\oeX$, defined by
\begin{align*}
\overline{\mathcal{L}}_{(\varepsilon, t, y)} f (x)  =  &   \ 
a(\eta (t), y)  \frac{\partial f}{\partial x} (x)  + \frac{1}{2} b^2 (\eta (t), y) 
\frac{\partial^2 f}{\partial x^2} (x) \\
&  + \WEint \left\{f (x + c (t, y, z)) - f (x)   -  c(t,y,z) \frac{\partial  f}{\partial x} (x) \right\}  \nu_t(dz)  + \frac{1}{2}  \Eint c^2(t,y,z) \nu_t(dz) \frac{\partial^2 f}{\partial x^2} (x), 
\end{align*}
for any $f \in \mathcal{C}^2(\RR)$. The quadratic variation of the local martingale $M_t = \int_0^t b(\eta(s),\oeX_\eta(s)) \frac{\partial u}{\partial x}(\oeX_s,s) d W_s$ is given for any $t\in[0,T]$ by 
\begin{align*}
    \langle M \rangle_t = \int_0^t (b(\eta(s),\oeX_{\eta(s)}))^2 \left(\frac{\partial u}{\partial x}(s,\oeX_s)\right)^2 d s.
\end{align*}
Hence, using the fact that $b$ and $\partial_x u$ have at most respectively linear growth by \ref{hyp:lipschitz} and polynomial growth by \eqref{eq:polynomial_growth_du}, one derive that
\begin{align*}
    \langle M \rangle_t \leq C_1 \big(1 + \sup_{t \in [0,T]} |\oeX_t|^{2 q +2}\big),
\end{align*}
for some positive constant $C_1$. 
Applying Lemma \ref{lem:uniform_bound_continuous_scheme},  
it follows that $(M_t)_{t \in [0,T]}$ is a $L^2(\Omega)$-martingale. We proceed similarly to show that $\left(\int_0^t d_\varepsilon(s,\oeX_{\eta(s)}) \frac{\partial u}{\partial x}(s,\oeX_s) d B_s \right)_{t \in [0,T]}$ is a $L^2(\Omega)$-martingale, observing that the coefficient $d_\varepsilon$ verifies
\begin{align*}
    d_\varepsilon^2(t,x) & \leqpt \Eint c^2(t,x,z) \nu_t(dz) \leq \left( \Eint  \overline{L}_c^2(t,z) \nu_t(dz) \right) (1+|x|)^2, 
\end{align*}
and the function $t\mapsto \Eint \overline{L}_c^2(t,z) \nu_t(dz)$ belongs to $L^1([0,T])$ thanks to \ref{hyp:lipschitz} and Hölder inequality.  
We can therefore reduce \eqref{eq:ito_1} to 
\begin{align} \label{eq: weak_step_1}
\mathbb{E} [\varPhi (\oeX_T) - \varPhi (X_T )] = & \int_0^T \EE \left[ \frac{\partial
  u}{\partial t} (t,\oeX_t) + \overline{\mathcal{L}}_{(\varepsilon, t, \oeX_{\eta(t^{{-}})})} u(t,\oeX_t) \right] dt
   = \int_0^T \EE \left[ \left( \overline{\mathcal{L}}_{(\varepsilon, t, \oeX_{\eta(t^{{-}})})} - \widetilde{X}_t - \widetilde{X}_{\eta(t)} \right) u(t,\oeX_t) \right] dt,
\end{align}
were we used the Kolmogorov PDE \eqref{eq:kolmogorov_pde} to rewrite the time derivative of $u$.
Making use of  pivot terms, we obtain the following decomposition 
\begin{align}
\left( \overline{\mathcal{L}}_{\varepsilon, t, \oeX_{\eta(t^{{-}})}} - \mathcal{L}_{(t)} \right) u(t,\oeX_t) & =  \mathcal{J}_{\oeX_{\eta(t^{{-}})}}^{(1)}(t,\oeX_t) + \mathcal{J}_{\oeX_{\eta(t^{{-}})}}^{(2)}(t,\oeX_{t^{{-}}} ) + \mathcal{J}_{\oeX_{\eta(t^{{-}})}}^{(3)}(t,\oeX_{t^{{-}}}) -\overline{\mathcal{H}}(t,\oeX_{t^{{-}}} ),
\end{align}
where we set  
\begin{align*}
  \mathcal{J}_{\oeX_{\eta(t^{{-}})}}^{(1)}(t,x) \coloneqq & (a({\eta(t^{-})}, {\oeX_{\eta(t^{-})}}) - a(t,x)) \frac{\partial u}{\partial x} (x, t) + \tfrac{1}{2} (b^2({\eta(t^{-})},{\oeX_{\eta(t^{-})}}) -  b^2(t,x)) \frac{\partial^2 u}{\partial x^2}(t,x), \\
  \mathcal{J}_{\oeX_{\eta(t^{{-}})}}^{(2)}(t,x) \coloneqq & \WEint \left\{ u(t,x+c(t,{\oeX_{\eta(t^{-})}},z)) - u(t,x+c(t,x,z))  - (c(t,{\oeX_{\eta(t^{-})}},z) - c(t,x,z))   \frac{\partial u}{\partial x}(t,x)\right\} \nu_t(dz),\\
   \mathcal{J}_{\oeX_{\eta(t^{{-}})}}^{(3)}(t,x) \coloneqq &  \tfrac12  \Eint \left( c^2 (t, \oeX_{\eta(t^{{-}})},z) -
    c^2 (t,x,z)\right) \nu_t(dz) \    \frac{\partial^2 u}{\partial x^2}(t,x),  
\end{align*}
and 
\begin{align} \label{eq:j_simplifie}
    \overline{\mathcal{H}}(t,\oeX_{t^{{-}}}) \coloneqq \Eint  \left(u(t,\oeX_{t^{{-}}}+c(t,\oeX_{t^{{-}}},z)) - u (t,\oeX_{t^{{-}}}) -  c (t, \oeX_{t^{{-}}}, z)  \frac{\partial u}{\partial x} (t,\oeX_{t^{{-}}}) - \tfrac{1}{2} c^2(t,\oeX_{t^{{-}}},z) \frac{\partial^2 u}{\partial x^2}(t,\oeX_{t^{{-}}})  \right)\nu_t (d z).
\end{align}

\paragraph{Upper-bound for $\EE[|\overline{\mathcal{H}}(t,\oeX_{t^{{-}}})|]$. } 
Using a Taylor expansion, one has:
\begin{align*}
    \overline{\mathcal{H}}(t,\oeX_{t^{{-}}}) = \frac{1}{2} \Eint c(t,\oeX_{t^{{-}}},z)^3 \left( \int_0^1 (1-\gamma)^2 \frac{\partial^3 u}{\partial x^3}(t,\oeX_{t^{{-}}} + \gamma c(t,\oeX_{t^{{-}}},z)) d\gamma \right) \nu_t(dz).
\end{align*}
Then we deduce from \eqref{eq:polynomial_growth_du} that
\begin{align} \label{eq:bound_of_j}
    \EE[|\overline{\mathcal{H}}(t,\oeX_{t^{{-}}})|] \leq \frac{1}{2} \EE \left[ \Eint c(t,\oeX_{t^{{-}}},z)^3 \left( \int_0^1 (1-\gamma)^2 (1 + |\oeX_{t^{{-}}}+\gamma  c(t,\oeX_{t^{{-}}},z)|^{q}) d\gamma  \right) \nu_t(dz) \right].
\end{align}
Using the sub-linearity of $c$ with respect to $z$ assumed in \ref{hyp:smooth}, and the fact that $|\gamma| \leq 1$, we can bound the right hand side of \eqref{eq:bound_of_j} as follows:
\begin{align*}
    c(t,\oeX_{t^{{-}}},z)^3 \left( \int_0^1 (1-\gamma)^2 (1 + |\oeX_{t^{{-}}}+\gamma  c(t,\oeX_{t^{{-}}},z)|^{q}) d\gamma  \right) \leqpt |z|^3 (1 + \sup_{t\in[0,T]} |\oeX_t|^{q+3}).
\end{align*}
Applying Lemma \ref{lem:uniform_bound_continuous_scheme}, {with  $q+3 \leq p$} one has:
\begin{align} \label{eq:bound_E_j}
    \EE[|\overline{\mathcal{H}}(t,\oeX_{t^{{-}}})|] \leqpt \Eint |z|^3 \nu_t(dz) \leq \Eint |z|^3 \mu(dz).
\end{align}
By definition of the Blumenthal--Geetor index $\beta$ of $\mu$, for any $0<\texttt{h}<3-\beta$:
\begin{align*}
    \EE[|\overline{\mathcal{H}}(t,\oeX_{t^{{-}}})|] \leqpt \varepsilon^{3-\beta-\texttt{h}   } \Eint |z|^{\beta+ \texttt{h}   } \mu(dz).
\end{align*}

\paragraph{Upper-bound of  $\sum_{i=1,2,3}\EE[\mathcal{J}_{\oeX_{\eta(t^{{-}})}}^{(i)}(t,\oeX_{t^{{-}}})]$.   }
This part is similar to standard weak convergence  proof arguments assuming smooth coefficients.  We will only detail the principle.  Using again  a Taylor expansion, the term  $\mathcal{J}_{\oeX_{\eta(t^{{-}})}}^{(2)}(t,\oeX_{t^{{-}}})$ gives a second order derivative of $u$ term: 
\begin{equation*}
\begin{aligned}
& \sum_{i=1,2,3}\mathcal{J}_{\oeX_{\eta(t^{{-}})}}^{(i)}(t,\oeX_{t^{{-}}}) \\
& = \left(a ({\eta(t)}, \oeX_{\eta(t^{{-}})}) -  a (t,\oeX_{t^{{-}}})\right)  \frac{\partial u}{\partial x}(t,\oeX_{t^{{-}}}) + \tfrac12 \left(b^2 ({\eta(t)}, \oeX_{\eta(t^{{-}})}) - b^2 (t,\oeX_{t^{{-}}})\right)  \frac{\partial^2 u}{\partial x^2}(t,\oeX_{t^{{-}}})  \\
&\quad  +  \tfrac12  \Eint \left( c^2 (t, \oeX_{\eta(t^{{-}})},z) -
    c^2 (t,\oeX_{t^{{-}}},z)\right) \nu_t(dz) \    \frac{\partial^2 u}{\partial x^2}(t,\oeX_{t^{{-}}}) \\
&\quad +   \WEint   ( c(t,\oeX_{\eta(t^{{-}})},z) - c(t,\oeX_{t^{{-}}},z))^2  \int_0^1 (1 -\gamma) \partial_{x^2}^2 u (\oeX_{t^{{-}}}+c(t,\oeX_{t^{{-}}},z) + \gamma  (c(t,\oeX_{\eta(t^{{-}})},z) -  c(t,\oeX_{t^{{-}}},z)) ) d\gamma  \  \nu_t(dz)\\
&\quad  + \WEint   ( c(t,\oeX_{\eta(t^{{-}})},z) - c(t,\oeX_{t^{{-}}},z)) \  c(t,\oeX_{t^{{-}}},z)  \int_0^1  \partial_{x^2}^2 u (\oeX_{t^{{-}}}  +  \gamma c(t,\oeX_{t^{{-}}},z)) d\gamma \ \nu_t(dz).
\end{aligned}
\end{equation*}
We observe that, by \ref{hyp:smooth}-\textit{(v)} for fixed $(t,x)$,  the map $z\mapsto c^2(t,x,z) + c^2(t,y,z)  + c(t,x,z) \  c(t,y,z)$ is in $L^1(\RR, \mu)$ (with a behaviour in $z^2$), as well as their derivatives with respect to $t$ and $x$.  
Each of the terms has a similar structure and we can summarise then, introducing the map,  for a fixed $t$, and $z$,   $(s\in(\eta(t),t) ,  \ x) \mapsto W(\eta(t),\oeX_{\eta(t^{{-}})}, z ; s,x)$ defined as 
\begin{equation*}
\begin{aligned}
& W(\eta(t),\oeX_{\eta(t^{{-}})}; s,x) \\
& = \left(a ({\eta(t)}, \oeX_{\eta(t^{{-}})}) -  a (s,x)\right)  \frac{\partial u}{\partial x}(s,x) + \tfrac12 \left(b^2 ({\eta(t)}, \oeX_{\eta(t^{{-}})}) - b^2 (s,x)\right)  \frac{\partial^2 u}{\partial x^2}(s,x)  \\
&\quad  +  \tfrac12  \Eint \left( c^2 (s, \oeX_{\eta(t^{{-}})},z) -
    c^2 (s,x,z)\right) \nu_s(dz) \    \frac{\partial^2 u}{\partial x^2}(s,x) \\
&\quad  +  \WEint   ( c(s,\oeX_{\eta(t^{{-}})},z) - c(s,x,z))^2 \left( \int_0^1 (1 -\gamma) \frac{\partial^2}{\partial x^2} u (x+c(s,x,z) + \gamma  (c(s,\oeX_{\eta(t^{{-}})},z) -  c(s,x,z)) ) d\gamma \right) \nu_s(dz)\\
&\quad  + \WEint   ( c(s,\oeX_{\eta(t^{{-}})},z) - c(s,x,z)) \  c(s,x,z) \left( \int_0^1 \frac{\partial^2}{\partial x^2} u ( x  +  \gamma c(s,x,z)) d\gamma \right)  \nu_s(dz).
\end{aligned}
\end{equation*} 
The map $W$ is differentiable in time,  and $\mathcal{C}^{2}$ in space, with controls on their derivatives.  In particular,
\begin{equation} \label{eq:bound_on_w_function}
\begin{aligned} 
\left(\left|\frac{\partial}{\partial s} W\right|  +   \left|\frac{\partial}{\partial x} W\right|  + \left|\frac{\partial^2}{\partial x^2} W\right| \right) (\theta, y; s, x) \leqpt  C_{a,b} ( 1 + |y|^2 + |x|^2) (1 + |x|^q) + C_c (1 + |y|^2 + |x|^2)( 1 + |x|^q + |y|^q), 
\end{aligned}
\end{equation}
Where $C_{a,b}$ is a positive constant that vanishes when $a$ and $b$ vanish, and $C_{c}$ is a positive constant that vanishes when $c$ vanishes. Hence, applying a second time the Itô formula
\begin{equation*}
\begin{aligned}
\left|\EE[\sum_{i=1,2,3}\mathcal{J}_{\oeX_{\eta(t^{{-}})}}^{(i)}(t,\oeX_{t^{{-}}})] \right|   = & 
\left|\int_{\eta(t)}^t  \EE[(\frac{\partial }{\partial s} + \overline{\mathcal{L}}_{\varepsilon, s, \oeX_{\eta(t^{{-}})}})   W (\eta(t),\oeX_{\eta(t^{{-}})};  s,\oeX_{s^{{-}}})]\ ds \right|. 
\end{aligned}
\end{equation*}
But with a Taylor expansion again, 
\begin{align*}
  \overline{\mathcal{L}}_{(\varepsilon, t, y)} f (x)  =  &   
  a(\eta (t), y)  \frac{\partial f}{\partial x} (x)  + \frac{1}{2} b^2 (\eta (t), y) 
  \frac{\partial^2 f}{\partial x^2} (x) \\
  &  + \WEint    c^2 (t, y, z) \left( \int_0^1 (1 - \gamma) \frac{\partial^2  f}{\partial x^2} ( x + \gamma c(t,y,z) ) d \gamma \right) \nu_t(dz) \\
  & + \frac{1}{2} \left( \Eint c^2(t,y,z) \nu_t(dz)   \right)^2 \frac{\partial^2 f}{\partial x^2} (x).
\end{align*}
All the integrals with respect to $\nu_t$ are uniformly bounded in time, {since  $c^2(t,y,z)$ behaves as $(1 + |y^2|) z^2$}. Thus we can drop to the conclusion,   counting the maximal power exponent $p\geq q + 4$ for the moments  $\sup_{ \varepsilon > 0} \EE[\sup_{\theta\in[0,T]} |\oeX_\theta|^p]$ we need to bound, 
\begin{equation*}
\begin{aligned}
\left|\EE[\sum_{i=1,2,3}\mathcal{J}_{\oeX_{\eta(t^{{-}})}}^{(i)}(t,\oeX_{t^{{-}}})] \right|  \leqpt (t - \eta(t)) (C_{a,b} + C_c). 
\end{aligned}
\end{equation*}
\paragraph{The weak error related to $\ueX$. }  The structure and arguments of the proof are all very similar for $\ueX$. We just need to adjust \eqref{eq: weak_step_1}, with 
\begin{align}\label{eq: weak_step_2}
\mathbb{E} [\varPhi (\ueX_T) - \varPhi (X_T )] 
   = \int_0^T \EE \left[ \left( \underline{\mathcal{L}}_{(\varepsilon, t, \ueX_{\eta(t^{{-}})})} - \widetilde{X}_t - \widetilde{X}_{\eta(t)} \right) u(t,\ueX_t) \right] dt,
\end{align}
and 
\begin{align*}
\left( \underline{\mathcal{L}}_{\varepsilon, t, \ueX_{\eta(t^{{-}})}} - \mathcal{L}_{(t)} \right) u(t,\ueX_t) & =  \mathcal{J}_{\ueX_{\eta(t^{{-}})}}^{(1)}(t,\ueX_t) + \mathcal{J}_{\ueX_{\eta(t^{{-}})}}^{(2)}(t,\ueX_{t^{{-}}}) + \mathcal{J}_{\ueX_{\eta(t^{{-}})}}^{(3)}(t,\ueX_{t^{{-}}}) -\underline{\mathcal{H}}(t,\ueX_{t^{{-}}}),
\end{align*}
where the $\mathcal{J}_{y}^{(i)}(t,x)$ are the same as before, whereas $\underline{\mathcal{H}}(t,x)$ looses its second order term: 
\begin{align*}
\underline{\mathcal{H}}(t,\ueX_{t^{{-}}}) & = \Eint  \left(u(t,\ueX_{t^{{-}}}+c(t,\ueX_{t^{{-}}},z)) - u (t,\ueX_{t^{{-}}}) -  c (t, \ueX_{t^{{-}}}, z)  \frac{\partial u}{\partial x} (\ueX_{t^{{-}}}, t) \right)\nu_t (d z)\\
& =  \Eint  c^2(t,\ueX_{t^{{-}}},z)  \int_0^1 (1 -\gamma) \frac{\partial^2 u}{\partial x^2} (t,\ueX_{t^{{-}}} + \gamma c(t,\ueX_{t^{{-}}},z)) d\gamma \ \nu_t (d z), 
\end{align*}
with 
\begin{equation*}
    c(t,\ueX_{t^{{-}}},z)^2   \frac{\partial^2 u}{\partial x^2} (t,\ueX_{t^{{-}}} + \gamma c(t,\ueX_{t^{{-}}},z))  \leqpt |z|^2 (1 + \sup_{t\in[0,T]} |\ueX_t|^{q+2}).
\end{equation*}
Then, for $ 0 < h < 2 - \beta$, 
\begin{align*} \label{eq:bound_E_j}
    \EE[|\underline{\mathcal{H}}(t,\oeX_{t^{{-}}})|] \leqpt \Eint |z|^2 \nu_t(dz) \leq \Eint |z|^2 \mu(dz)
\leqpt \varepsilon^{2-\beta-\texttt{h}   } \Eint |z|^{\beta+ \texttt{h}   } \mu(dz) \leqpt \varepsilon^{2 - \beta^+}. 
\end{align*}
\end{proof}

\section{Evaluation of the convergences rates through numerical experiments } \label{sec:numerics} 

\subsection{Strong convergence}
In this section, we investigate the optimality of the strong convergence rate of $\oeX$ with numerical simulations. We examine first the behaviour of the strong $L^p$--error in terms of its norm exponent $p$.  Then, we test the actual sensitivity of the error to the  integrability condition  driven by the parameter $\zeta$ in  \ref{hyp:peano}, based on Example~\ref{example:nu} for the choice of $\nu_t(dz)$. 
In these two behavioural studies, the numerical results appear to validate our theoretical result. However, there are two numerical limitations to this analysis that we would like to point out and  which are specific to strong error analysis.

\paragraph{Lack of exact trajectory solution. } Processing the computation  of the strong norm of the error always poses the problem of simulating a reference solution trajectory. In most of the cases where  the stochastic integral is not additive, it is impossible to have a solution that can be simulated exactly. 
It is therefore a common practice to compute an approximate reference solution, by pushing the approximation parameters to a limit value which serves as a bound for the experiments with coarse parameters.  However this approach requires to be able to construct coarser increments by summing up small ones. In our settings, this seems only possible when $c(s,x,z) = c(s,x)f(z)$, that limits our numerical investigation to such test cases for strong error simulation. 

\paragraph{Sampling the two-parameters increments. } Moreover, for the $\varepsilon$-EM algorithm under study, we have two  control-parameters, the time step $1/n$ and the small jumps cut $\varepsilon$.   Once the choice of $\varepsilon$ is fixed, the increments of the process $\int_0^\cdot \int_{\mathbb{R}/ B(\varepsilon)} f(z) \rpm$ can then be simulated on a very fine time grid,  and  next aggregated together to produce increments on a coarser time grid. But reversing the choices by fixing first a fine time grid before choosing $\varepsilon$ seems only feasible  in the  weak error framework. 
\medskip

Considering the above constraints, we leave aside the investigation of the scheme's behaviour in $\varepsilon$ to the analysis in weak convergence. For that reason, the comparison between $\oeX$ and $\ueX$ is also postponed to the weak error analysis.
In our strong error framework, we adopt the following strategy to sample the reference solution trajectories:  we fix the maximal step number to be $n_{\max} = 2^{23}$.  From Corollary \ref{cor:mainresult}-\eqref{eq:explicit_varepsilon}, selecting $\varepsilon_{\min} = n_{\max}^{ - \left( \frac{1}{2}+\left(\gamma \wedge \frac{1}{p}\right)\right)}$ ensures a minimal  $L^p$-distance between the reference trajectory $\overline{X}^{\varepsilon_{\min}, n_{\max}}$ and the exact trajectory 
\begin{align*}
    \left\| \sup_{0 \leq i \leq n}  \big| X_{t_i} - \overline{X}^{\varepsilon_{\min}, n_{\max}}_{t_i}\big| \right\|_{L^p(\Omega)} \quad = \quad \mathcal{O}\left(n_{\max}^{- \left\{ \gamma \wedge \left( \frac{2 \zeta }{p (1+\zeta)} \right) \right\} }\right), 
\end{align*}
regardless of the integrability exponent $\zeta \in (0,1]$, provided that the jump coefficient satisfies Condition \eqref{eq:bound_on_the_jump_coefficient}. As it is usually the case in strong error analysis, the numerical simulations allow then to investigate the error between the scheme  $\overline{X}^{\varepsilon_{\min}, n}$ for $n=2^k$, with $k< 23$, and the reference trajectory  $\overline{X}^{\varepsilon_{\min}, n_{\max}}$:
\begin{align*}
\left\| \sup_{0 \leq i \leq n}  \big| \overline{X}^{\varepsilon_{\min},n}_{t_i} - \overline{X}^{\varepsilon_{\min}, n_{\max}}_{t_i}\big| \right\|_{L^p(\Omega)}  \quad = \quad \mathcal{O}\left(n^{- \left\{ \gamma \wedge \left( \frac{2 \zeta }{p (1+\zeta)} \right) \right\} }\right).
\end{align*}
Since the reference trajectory is constructed on the basis of the small jumps cut at $\varepsilon_{\min}$, it turns out that this numerical evaluation is essentially focused on the Peano part of the error.

\subsubsection{Rate of convergence in terms of the norm exponent $p$}\label{sec:num_for_rem_optimal}

First, we investigate the behaviour of the $L^p$-strong error with respect to the choice of $p\geq 2$. For that purpose, we consider the following example, with $c(s,x,z) = \cos(x) z$ satisfying the condition \eqref{eq:bound_on_the_jump_coefficient}: 
\begin{align*}
    X_t = \int_{0}^t \cos(X_s) ds + \int_{0}^t \int_{-\infty}^{+\infty} \sin(X_{s^-})\  z \  \rpm,
\end{align*}
where the Poisson random measure $\widetilde{N}$ has time-homogeneous compensator 
\[\nu_s(dz) ds = \indic{|z|\leq 1} \frac{dz}{|z|^{3/2}} \ ds,\]
corresponding to a truncated $\frac{1}{2}$-stable noise-process. The coefficients are time-homogeneous and  Lipschitz in space, hence all the hypothesis of Theorem \ref{theoreme_maurer} are satisfied with $\gamma = \zeta = 1$. Moreover, Corollary \ref{cor:mainresult} holds, providing the theoretical rate of convergence $n^{-\frac{1}{p}}$.   
\medskip
\begin{figure}[H]
    \centering
    \includegraphics[width=0.65\paperwidth]{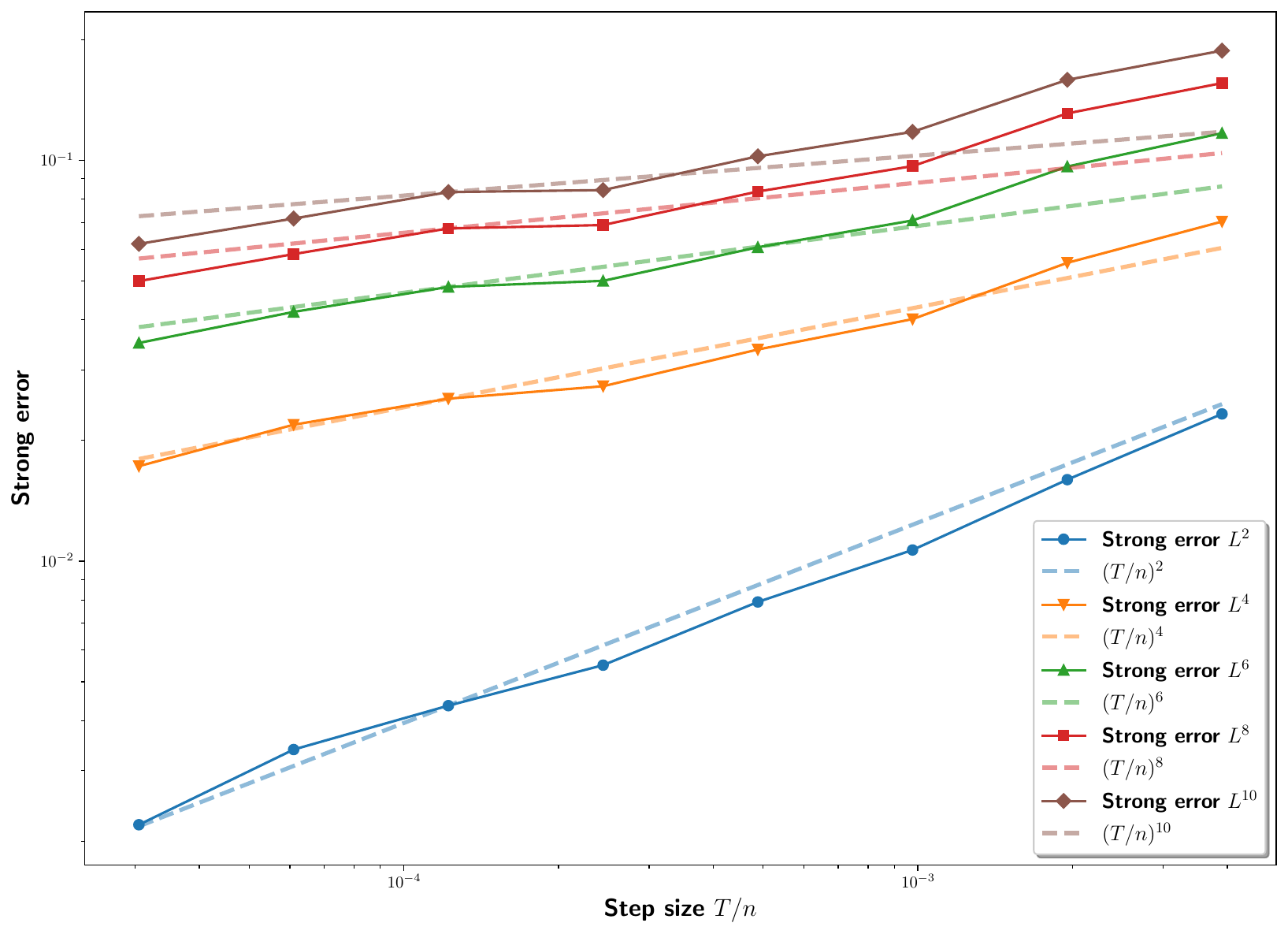}
    \caption{\em Behaviour of  $\left\| \sup_{0 \leq i \leq n}  \big| \overline{X}^{\varepsilon_{\min},n}_{t_i} - \overline{X}^{\varepsilon_{\min}, n_{\max}}_{t_i}\big| \right\|_{L^p(\Omega)}$ in terms of $T/n$, for various $L^p$-norms (lines with makers), and the corresponding theoretical (dash lines) rates.\label{fig:fournier_strong_lp}}
\end{figure}
The results shown in Figure \ref{fig:fournier_strong_lp} are performed with a $10^5$-paths Monte-Carlo simulation of the strong error, for $p=2\ell$, $\ell = 1,\ldots, 5$,  and $\sup$ in time is computed on the interval $[0,T\!=\!1]$.    Figure \ref{fig:fournier_strong_lp} reports the $L^p$-error in terms of the step size $T/n$, with $n=2^k$, $k=9,\ldots,17$,  varying over several decades,  and all computed with $\varepsilon_{\min} = n_{\max}^{-1}$, corresponding to the optimal value for the $L^2$-norm.

Even in this case, where \ref{hyp:peano} is trivially satisfies by both  time-homogeneous jump coefficient and  compensator,  the observations clearly indicate a dependence in $p$ in the rate of convergence with $n$. This behaviour stands out strongly from the EM-scheme for diffusions where (under adequate hypotheses) all the  $L^p$-errors converge with same rate $\frac{1}{2}$.    
In particular, for the $L^2$, $L^4$  and $L^6$ norms, the observed rates fit well the theoretical ones of $\frac{1}{2}$, $\frac{1}{4}$ and $\frac{1}{6}$, while the simulations seem to do  slightly better than expected for greater values of $p$. {This slight improvement may be attributed to our comparison with an approximate reference trajectory rather than an exact one.}

\subsubsection{Rate of convergence with low time-integrability} \label{ssec:numerical_low_integrability}

When the compensator of the Poisson random measure is time-dependent and presents a singularity on the interval $[0,T]$, according to Lemma \ref{lem_peano},  it results in a loss of speed of convergence for the Euler-Peano scheme which should theoretically also affect the $\varepsilon$-EM scheme.  More precisely, when $\psi_p$ is not better than $L^{1+\zeta}([0,T])$ with $\zeta \in (0,1)$, we may only recover the rate $n^{-\frac{2 \zeta}{p(1+\zeta)}}$ (assuming furthermore that all the coefficients are  Lipschitz in time, i.e. $\gamma = 1$).
 
\begin{figure}[H]
\centering
\includegraphics[width=0.65\paperwidth]{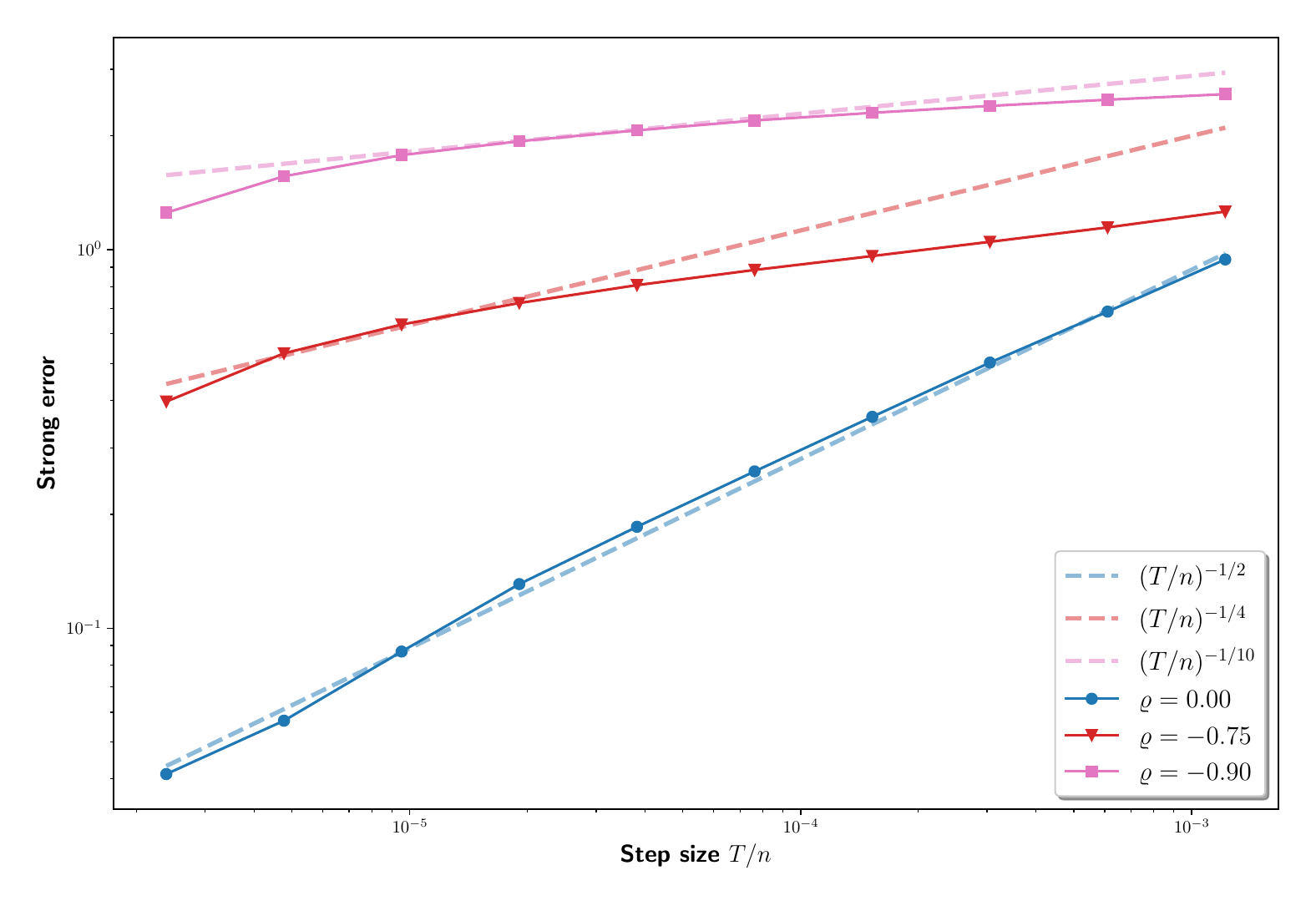}
\caption{\em Behaviour of  $\left\| \sup_{0 \leq i \leq n}  \big| \overline{X}^{\varepsilon_{\min},n}_{t_i} - \overline{X}^{\varepsilon_{\min}, n_{\max}}_{t_i}\big| \right\|_{L^2(\Omega)}$ in terms of $T/n$, for $\varrho\in\{0,-0.75,-0.9\}$  (lines with makers), and the corresponding theoretical  rates (dash lines). \label{fig:time_low_integrability}}
\end{figure}

To emphasise the influence of the time-integrability of the compensator on the convergence rate, we consider  the compensator
\[ \nu_s(dz) ds = \indic{|z|\leq10} \frac{dz}{|z|^{3/2}} \ s^\varrho  ds,\]
where $\varrho$ is negative and belongs to $(-1,0]$. We introduce this noise in the following SDE: 
\[X_t = 1+ \int_0^t \sin(X_s) ds + \int_0^t \int_{-\infty}^{+\infty} \cos(X_{s^-}) \ z \rpm,\]
for three different values of $\varrho$, namely $\varrho \in \{0, -0.75, -0.9\}$.  With the help of  Remark \ref{rem:low_integrability}, we compute the expected theoretical rate of convergence in each case: 
\begin{itemize}
    \item{$\varrho = 0$} corresponds to the standard time-homogeneous equation, leading to the rate of convergence $n^{-\frac{1}{2}}$; 
    
    \item{for $\varrho = -0.75$}, one has $\psi_p \in L^{1+\zeta}$ for every $\zeta \in (0,\frac{1}{3})$ but $\psi_p \notin L^{\frac{4}{3}} $, so the theoretical convergence rate in the $L^2$-norm should not be better than $n^{-\frac{1}{4}}$; 
    
    \item{for $\varrho = -0.9$}, one has $\psi_p \in L^{1+\zeta}$ for every $\zeta \in (0,\frac{1}{9})$ but $\psi_p \notin L^{\frac{10}{9}}$, so the theoretical convergence rate in the $L^2$-norm should not be better than $n^{-\frac{1}{10}}$.
\end{itemize}
We performed a $10^5$-paths Monte-Carlo simulation of the $L^2$-strong error on the time interval $[0,T\!=\!1]$. The computed errors are reported in Figure~\ref{fig:time_low_integrability} in terms of the step size $T/n$, with $n=2^k$, $k=11,\ldots,20$, varying over several decades, and computed with $\varepsilon_{\min} = n_{\max}^{-\frac{1}{2} - \frac{2(1-\varrho)}{p}}$, corresponding to the optimal value for the $L^2$-norm at each value of $\varrho$.

We observe in Figure \ref{fig:time_low_integrability} a clear influence of the parameter $\varrho$ in the error rates, which corresponds well with the expected theoretical ones. However, the algorithm seems to perform better than expected for the smallest value of $n$, especially in the cases of lower convergence ($\varrho = -0.9$ and $\varrho = -0.75$).  
This improvement can come from several sources, in first place, the comparison with a reference trajectory produced with the same algorithm favours error underestimation when $n$ is approaching the $n_{\max}$. And this is all the more evident as the convergence rate is low. 

\subsection{Weak convergence}

In this section, we provide numerical experiments for the theoretical weak rate of convergence. Making use of the Feynman-Kac formula, we get analytical expression for the reference values that can be compared to the numerical simulation in some specific cases. According to Theorem \ref{thm:weak_error_bound}, this allows us to obtain the behaviour of the weak error in terms of $\varepsilon$, while choosing the adequate number of time steps $n$. 

We consider two  test cases. In the first one,  we illustrate the behaviour of the convergence rate in response to the Blumenthal--Geetor index of a  singular Lévy measure.  While we restrict the SDE to the case of  multiplicative jump coefficient, we consider the case of a truncated $\alpha$-stable compensator for different values of $\alpha$. We observe that the theoretical rates of convergence are matched with a good precision.   

In the second case, we exploit the flexibility of the representation {\it in law} of the exact quantity to approximate, allowing  us to consider a much more complex jump coefficient involving  a  noise structure that cannot be represented as the increment of a stochastic process. We combine this jump coefficient with a truncated 0-stable compensator, and show that the convergence is still numerically achieving in this case, with the expected theoretical rate.
\medskip

The numerical tests presented in the following are based on two main methods.

\paragraph{Feynman-Kac formula with source term. }
We consider a $\mathcal{C}^4$-test function $\Phi$, together with a $\mathcal{C}^{1,4}$-function $G$,  and the flow $(X_{\theta}(t,x), \theta \geq t)$, solution of \eqref{eq:intro_SDEs} starting at $x$ at time $t$.  The Feynman-Kac formula proven in the first step proof of  Theorem \ref{thm:weak_error_bound} can be easily extended, with the same standard arguments to show that under \ref{hyp:smooth} the solution $u$ of the Kolmogorov PDE 
\begin{equation*}
\begin{aligned}
\left\{
    \begin{array}{ll}
         \frac{\partial u}{\partial t} +\mathcal{L}_{(t)} u  = 0, \quad (t,x)\in[0,T)\times\mathbb{R}, \\
         u (T, x)  = \varPhi (x), \quad x \in \mathbb{R},
    \end{array}
\right.
\end{aligned}
\end{equation*}
satisfies
\begin{equation*}
    u(t,x) = \EE\left[ \Phi(X_T(t,x)) - \int_t^T G(\theta, X_{\theta}(
    s,x)) ds \right].
\end{equation*}
Now, for a well chosen function $u(t,x)$, the function $G(t,x)$ can be explicitly identified. This method provides a way to compare a complex statistic of type  $\EE[\Phi(X_T(t,x))]$ to an exact reference value, up to the additional integral term $\EE[\int_t^T G(\theta, X_{\theta}(s,x)) ds]$.  As for the Feynman-Kac formula, the weak convergence rate  result in Theorem \ref{thm:weak_error_bound} can be  straightforwardly extended to that combination of test-functions, as soon as the time integral on $G$ is also discretised with a (at least) one-order numerical integration method.     

\paragraph{Sampling of the two-parameters noise. }   In this context of weak error analysis, it is now essential to vary $\varepsilon$ as much as $n$, while maintaining the Monte Carlo sample size of simulated paths sufficiently large (in view of the biases to be studied). We choose $\varepsilon = (1.5)^k$ for $k \in \{3,\ldots,8\}$ and for each value of $\varepsilon$, we consider the corresponding number $n_{\varepsilon}$ of time steps to match the rate of convergence given in Theorem \ref{thm:weak_error_bound}. From there, we sample the associated trajectory of the $\varepsilon$-EM scheme with Gaussian compensation ${(\oeX_{i \frac{T}{n_{\varepsilon}}}, i = 0,\ldots,n_{\varepsilon})}$ and the version without compensation ${(\ueX_{i \frac{T}{n_{\varepsilon}}}, i = 0,\ldots,n_{\varepsilon})}$. We can next sample approximations of $\EE\left[ \Phi(\oeX_T(t,x)) - \int_t^T G(\theta, \oeX_{\theta}(t,x)) ds \right]$ and $\EE\left[ \Phi(\ueX_T(t,x)) - \int_t^T G(\theta, \ueX_{\theta}(t,x)) ds \right]$ by the Monte-Carlo method, and using a $\tfrac{1}{3}$-Simpson rule to evaluate the time-integral of $G$. 

\subsection{Multiplicative jump coefficient with truncated $\alpha$-stable compensator}

In this test case, we consider the multiplicative jump coefficient $c(t,x,z) = \sin(x) \; z$,  associated with the  10-truncated $\alpha$-stable compensator $\nu_t(dz) dt= \tfrac{ 1}{|z|^{1+\alpha}} \indic{|z| \leq 10} dz \,dt$. We also add the linear drift $b(x) = -2 x$, so that the variance of the jump term does not increase the variance of the solution too quickly. This corresponds to the SDE 
\begin{align*}
    X_t = x - \int_0^t 2 X_s \ ds  + \int_0^t \int_{-\infty}^{+\infty} \sin(X_{s^-}) z \rpm, \quad t \in [0,T].
\end{align*}
We choose the Kolmogorov function $u$ and the test function $\Phi$ as
\begin{align*}
    u(t,x) & = \left(1-\tfrac{1}{2}e^{T-t}\right) x^2, \\
    \Phi(x) & = u(T,x) = \tfrac{1}{2} x^2,  \quad (t,x)\in[0,T]\times \RR.
\end{align*}
In this setting, we may compute the function $G(t,x)=\tfrac{\partial u}{\partial t}(t,x) + \mathcal{L}_{(t)} u(t,x)$ as 
\begin{equation*}
    G(t,x) = \tfrac{1}{2} x^2 e^{T-t}+ (1-\tfrac{1}{2} e^{T-t}) \ \sin^2(x) \  2 \frac{(10)^{2-\alpha}}{2-\alpha}.
\end{equation*}
Note that both the coefficients $b,c$ and the test function $\Phi$ (and $G$) satisfy \ref{hyp:smooth}. Using Theorem \ref{thm:weak_error_bound} and its specific remark concerning $\alpha$-stable type measures, this leads to the respective theoretical rates of convergence of $\varepsilon^{3-\alpha}$ for the $\varepsilon$-EM scheme $\oeX$ where the small jump are compensated with Gaussian variables and $\varepsilon^{2-\alpha}$ for the $\varepsilon$-EM scheme $\ueX$,  where they are simply ignored.
\begin{figure}[H]
    \includegraphics[width=0.7\paperwidth]{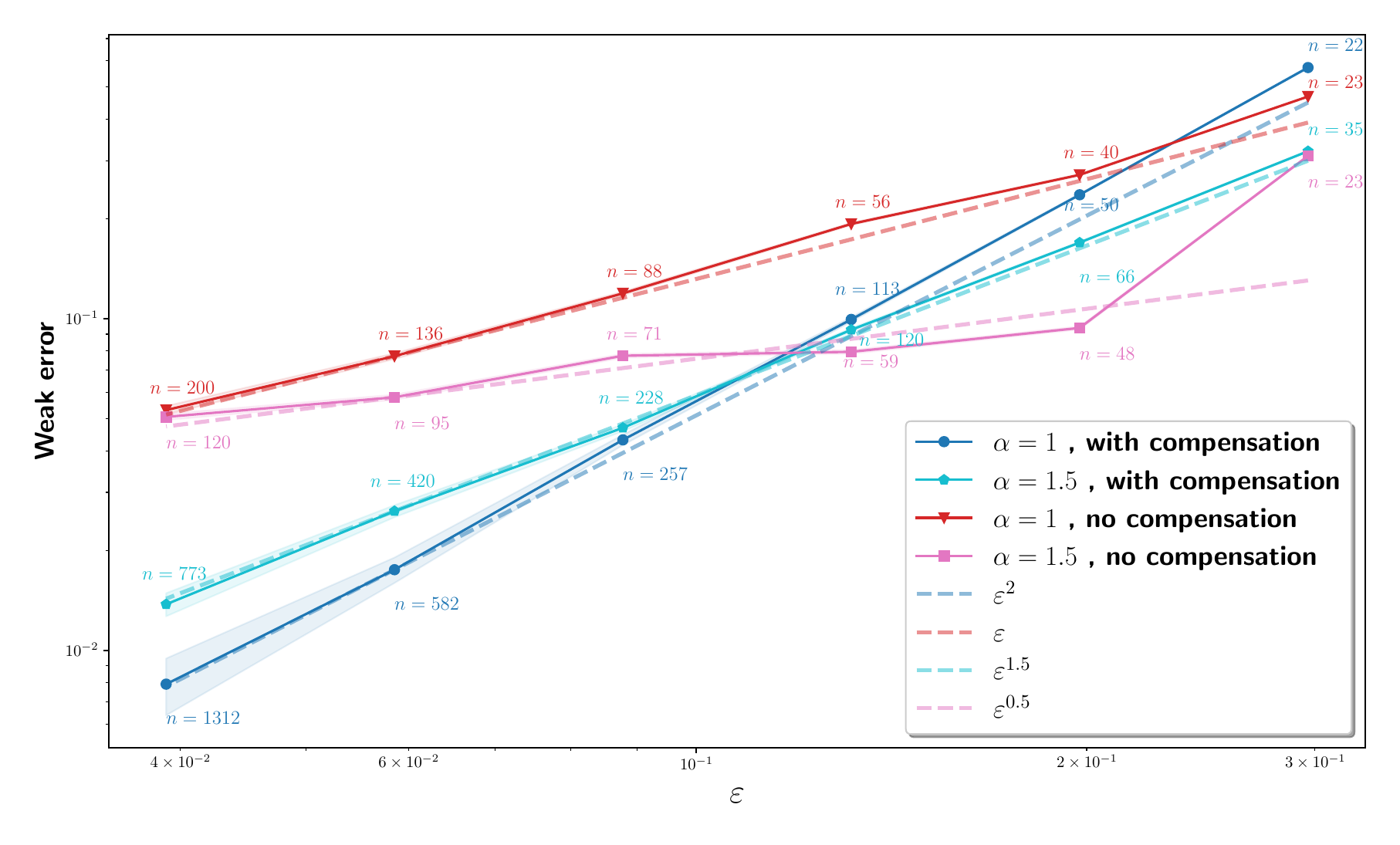}
    \caption{\em Behaviour of the weak error for the $\varepsilon$-EM schemes, for the values $1$ and $1.5$ of the noise parameter $\alpha$. \\
    Blue lines with $\circ$-markers report on the
    $\oeX$ version, with Gaussian compensation for the small jumps. \\
    Red lines with {\scriptsize{\text{$\square$}}}-{\scriptsize{\text{$\triangle$}}}-markers report on the $\ueX$ version, without compensation.\\
    The accompanying dash lines stand for the theoretical rates.  \label{fig:weak_error_multiplicative_case}}
\end{figure}
The results shown in Figure \ref{fig:weak_error_multiplicative_case} are  performed with  $10^8$-paths Monte-Carlo approximation of $\EE[\Phi(\oeX_T(t,x)) - \int_t^T G(\theta, \oeX_{\theta}(s,x)) ds]$ and $\EE[\Phi(\ueX_T(t,x)) - \int_t^T G(\theta, \ueX_{\theta}(s,x)) ds]$, that we compare in absolute value to the exact reference $u(t,x)$. We chose the parameters $(t,x)=(0,10)$, $T=1$, and $\varepsilon = (1.5)^{-k}$ for $k = 3, \ldots, 8$. The number $n_{\varepsilon}$ of time step is set to $\overline{n_{\varepsilon}} = 6 \times\lfloor \varepsilon^{-(3-\alpha)} \rfloor $ for the compensated scheme $\oeX$, and $\underline{n_{\varepsilon}} = 24 \times\lfloor \varepsilon^{-(2-\alpha)} \rfloor $ for the non-compensated one.  The observed errors reported in Figure \ref{fig:weak_error_multiplicative_case} match  very well the expected rates for both values of $\alpha$. 

\subsection{The case of a non-multiplicative jump coefficient} \label{sec:weak_error_arctan}

When the jump coefficient $c(x,z)$ cannot be written as $\bar{c}(x)g(z)$ for some functions $\bar{c}$ and $g$, then there is no explicit increment of a stochastic process behind the noise structure of the associated SDE. 
To the best of our knowledge,  the approximation of such case has not been numerically studied in the literature either so far. 
\medskip

There is two main difficulties that may arise when dealing with the simulation of non-increment based jump noise: the first lies in the simulation of the $\varepsilon$-EM scheme itself, because the variance of the small jump $\int_{t_{i-1}}^{t_i} \int_{B(\varepsilon)} c^2(s,\oeX_{t_{i - 1}},z) \rpc$ that appear in the Gaussian compensation term of the scheme might not have a closed form. The second, which is more specific to the weak error computation, is that the function $G(t,x)$ from the Kolmogorov equation might also not have a closed  form. 

Since it seemed difficult to find a non-trivial function $c(x,z)$ that allows to compute both the small jump integral and the function $G(t,x)$ exactly, we have chosen to focus on finding a coefficient $c$ that allows a closed form for $G$, to minimise the impact on the weak error estimation. 

For these reasons, we consider the jump coefficient $c(x,z) = \text{arctan}(x\, z)$ with the 1-truncated 0-stable compensator $\nu(dz) dt=  \indic{|z| \leq 1} \tfrac{1}{|z|} dz dt$. We also maintain the linear drift $b(x) = -2 x$. This corresponds to the SDE 
\begin{equation} \label{eq:sde_atan_drift}
    X_t = x - \int_0^t 2 X_s ds + \int_0^t \int_{-\infty}^{+\infty} \text{arctan}(z \ X_{s^-} ) \rpm, \quad t \in [0,T].
\end{equation}
In order to simplify the analytical form of the exact solution, we choose the Kolmogorov function $u$ to be homogeneous in time, equals to the test function $\Phi$: 
\begin{equation*}
\begin{aligned}
    u(t,x) & = \sin(x), \\
    \Phi(x) &  = \sin(x),  \quad (t,x)\in[0,T]\times \RR.
\end{aligned}
\end{equation*}
In this setting, we may obtain a closed form for the function $G(t,x)= \tfrac{\partial u}{\partial t}(t,x) + \mathcal{L}_{(t)} u(t,x)$, as
\begin{equation*}
    G(t,x) =  - 2 x \cos(x) + 2 \sin(x) (\log(2) - \log(\sqrt{x^2 + 1} +1 )).
\end{equation*}
We leave the details on the computation of $G$ and on the method used to approximate the small jump integral in Appendix~\ref{sec:details_arctan}.
In this case again,  the coefficients $b,c$ satisfy \ref{hyp:smooth}. The tests functions $\Phi$ and $G$ satisfy \eqref{eq:hypo_PHI}. Using Theorem \ref{thm:weak_error_bound}, the obtained rate of convergence behaves theoretically in $\mathcal{O}(\varepsilon^3)$ for $\oeX$ and in $\mathcal{O}(\varepsilon^2)$ for $\ueX$.

\begin{figure}[H]
    \centering
    \includegraphics[width=0.7\paperwidth]{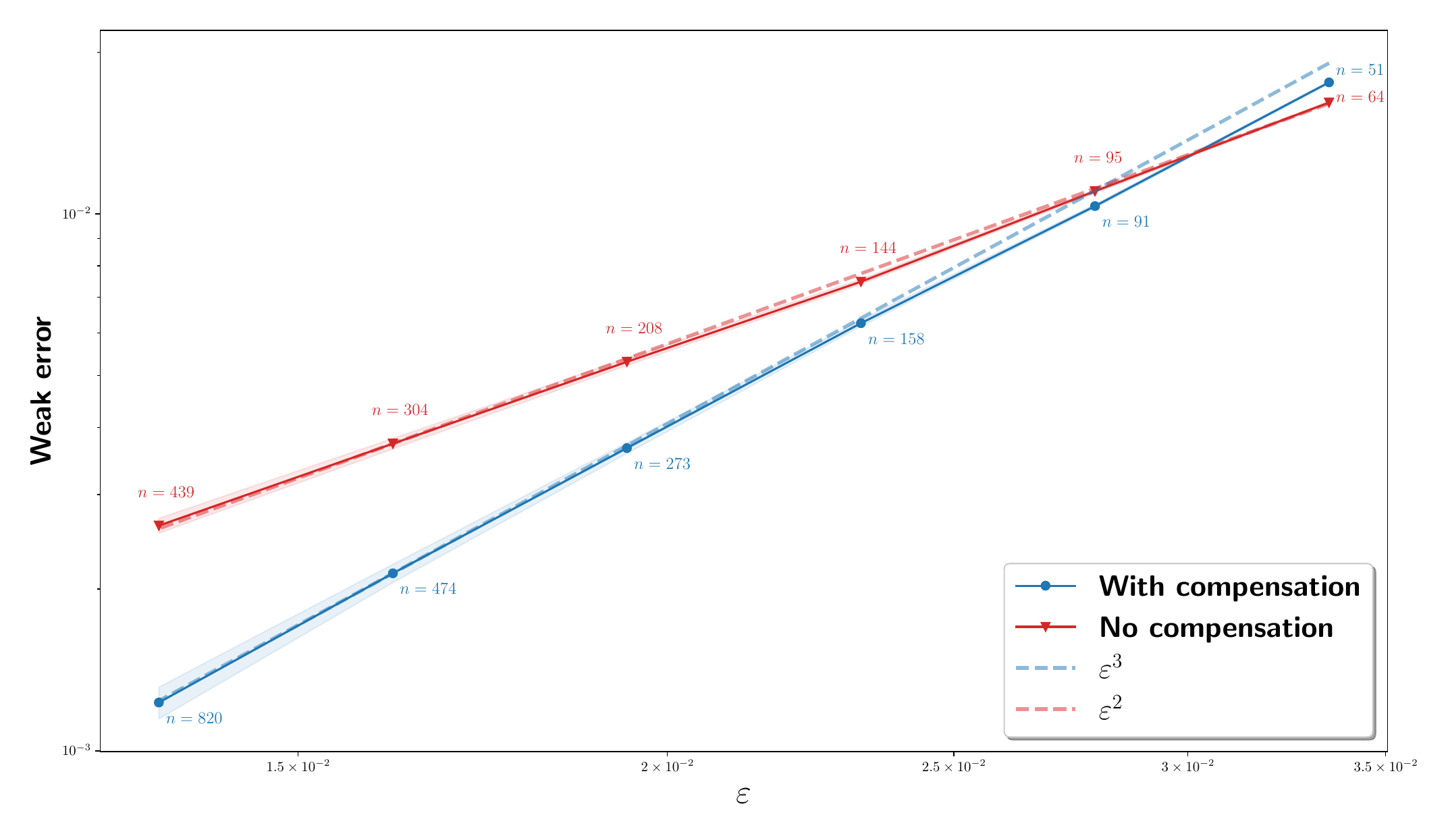}
    \caption{\em  Behaviour of the weak error for the $\varepsilon$-EM schemes in the non-multiplicative case \eqref{eq:sde_atan_drift}. \\
    Blue lines with $\circ$-markers report on the
    $\oeX$ version, with Gaussian compensation for the small jumps.      Red lines with {\scriptsize{\text{$\triangle$}}}-markers report on the $\ueX$ version, without compensation. 
    The accompanying dash lines stand for the theoretical rates.  \label{fig:weak_error_NONmultiplicative_case}}
\end{figure}
The results shown in Figure \ref{fig:weak_error_NONmultiplicative_case} are  performed with $10^9$-paths Monte-Carlo approximation of $\EE[\Phi(\oeX_T(t,x)) - \int_t^T G(\theta, \oeX_{\theta}(s,x)) ds]$ and $\EE[\Phi(\ueX_T(t,x)) - \int_t^T G(\theta, \ueX_{\theta}(s,x)) ds]$, that we compare in absolute value to the exact reference $u(t,x)=\sin(x)$. We chose the parameters $(t,x)=(0,10)$, $T=1$, and $\varepsilon = 10^{-1}\times (1.2)^{-k}$ for $k = 6, \ldots, 11$. The number $n_{\varepsilon}$ of time step is set to $\overline{n_{\varepsilon} }= 2 \times\lfloor 10\times \varepsilon^{-3} \rfloor $ for the compensated SDE and $\underline{n_{\varepsilon} }= 8\times \lfloor 10\times \varepsilon^{-2} \rfloor $ for the non-compensated one. 
The observed errors reported in Figure \ref{fig:weak_error_NONmultiplicative_case}  match very well the expected rates for the both versions of the scheme. 

\medskip

We would like to point out, however, that numerical tests in this case are particularly costly to carry out. This is due to the cost of the noise sampling and the computation of the compensations involved in the simulation of the scheme. For example, on a 160-cores 2400 Hz Intel Xeon CPU node, it took respectively around 48 hours 
and 19 hours 
to compute  the  $\oeX$ error value respectively for  $n=820$ and $n=474$,  with $10^9$ Monte-Carlo samples. 

\section{Application: the enhanced stochastic dynamics of rigid fibres in turbulence}\label{sec:application}

In this section, we would like to illustrate, at least on a qualitative level, the capability of a stochastic model driven by inhomogeneous jumps  measure to  replicate some  specific tail shapes transition, also called Lévy wings. Such Lévy wings are  observed in physics, on the probability density functions of  angular increments of small rod orientations in 2D turbulent flow.  We refer precisely to the Figure 8 in \cite{campana} where, in addition, a first stochastic Brownian SDE was introduced to model the dynamics of such  orientation angles in turbulence. 
\medskip

While stochastic models of the translational dynamics of (very) small objects transported by turbulent flow are well studied  (see e.g review in \cite{Pope1994} and the papers that have quoted it since),  stochastic approaches for the rotational dynamic of non-spherical objects in this same context remains little explored (see \cite{campana,CAMPANA_CFD} and references therein). In \cite{campana}, a  Brownian-noise based model for the  orientation (equivalently the angular velocity) is proposed and  compared  against  direct simulation of small rods in a 2D-turbulent flow.  Such direct numerical simulations (DNS) of the flow motion (solving the 2D Navier Stokes equation describing the velocity) are performed in order to resolve all the scales of turbulence. The result can then be considered as reference experiments for  phenomenological study  and  model validation.

The stochastic  model is derived from the analysis of the  instantaneous tangent dynamics of the turbulent flow $d {\bm r}(t) = d \mathbb{A}(t) r(t)$, where $\mathbb{A}$ is  the velocity gradient tensor of the flow at the (centre of mass) particle position. Based on empirical ergodicity arguments,  this gradient is further  decomposed in a given mean gradient $\av{\mathbb{A}}$ and a white noise part, represented with a given (measured) tensor  $\mathbb{B}$, multiplying  a $2\times{2}$ random matrix $\mathbb{W}$  made of independent Brownian motions. The model starts with the diffusion equation (in Stratonovich's sense) 
\begin{equation}\label{eq:separation_model}
	d \bm r(t) =  \av{\mathbb{A}} \, \bm r(t) \, d t+ (\mathbb{B}\, \circ d \mathbb{W}_t) \bm r(t). 
\end{equation}
When the flow is homogeneous in space with a shear,  $\av{\mathbb{A}}$  is reduced to the anisotropic  entries  $\av{\mathbb{A}_{ij}} = \sigma \delta_{i,1}\delta_{j,2}$, with the shear parameter $\sigma$ amplifying the anisotropic behaviour of the flow. 
The model \eqref{eq:separation_model} can then be further transformed in a diffusion equation for the (scalar) angular dynamics of the orientation angle $\theta_t = \arctan({r_y}/{r_x})$. With some measured coefficients $\sigma \mapsto \gamma_i(\sigma)$, this leads to the following SDE
\begin{equation}\label{eq:folded_theta}
d\theta_t = \frac{\sigma}{2}\left(\cos(2\theta_t)-1\right) dt + \sqrt{ \gamma_0 +\gamma_1\sin(2\theta_t) 
	+\gamma_2\sin(4\theta_t)
	+\gamma_3 \cos(2\theta_t) 
	+\gamma_4 \cos(4\theta_t)} \circ d W_t,
\end{equation}
where $W_t$ is now a one-dimensional Brownian motion. The evaluation of the model \eqref{eq:folded_theta} in \cite{campana}  can be summarised as follows. When the folded angle $\overline\theta_t = (\theta_t + \tfrac{\pi}{2}\, \text{mod }  \pi) - \tfrac{\pi}{2}$ is compared with the direct simulation of the flow, the model reproduces well the  measured PDF for the various values of $(\sigma, \gamma_i(\sigma))$  generated by the experiments.  

The situation is less straightforward when comparing the means and variances for the cumulative angle $\delta \theta_t = \theta_t - \theta_0$.   In particular the empirical statistics of the cumulative angle  presents  characteristics of an anomalous diffusion before to reach a diffusive regime. 

In the transitional regime, the behaviour of the variance $\EE[\delta \theta_t^2] - \EE[\delta \theta_t]^2$ is  in $\mathcal{O}(t^{q})$ with $q>1$, until a time $T^*$ and then turns to a linear behaviour with larger cumulative times. 
The values of $q$ and $T^*$ seem to be modulated with the shear parameter $\sigma$. This emphasises that the model \eqref{eq:folded_theta} is designed for a long-term perspective on the phenomena, without the ability to capture the transition regime characterised by superdiffusive variance and heavy-tailed distribution. 
Observing the vortex structures that emerge and vanish in 2D turbulence, it can be considered that a small transported rod has the potential to become ensnared within it, initiating a rapid swirling motion until the vortex dissipates. When observed from a longer-term perspective, this phenomenon can be interpreted as a sudden jump.

We then seek for a superdiffusive noise with jump,  that eventually converges to a Brownian motion in large times.  We introduce the time-inhomogeneous truncated $\alpha$-stable process
\[ L_t =
\int_0^t \int_{-\infty}^{+ \infty} z \ (N (d s, d z) - \nu_s (d z) d s)\] 
 with compensator 
 \[\nu_s (d z) d s = \kappa (s) | z |^{-
\alpha - 1}  \mathds{1}_{\{ -z_{\ast} \leq  z  \leq z^{\ast} \}} dz \ ds.\]
In the jump integral above,   $(z_{\ast}, z^{\ast})$ are non negative scalars driving the anticyclonic and cyclonic mean jump sizes, the time-dependent normalisation constant  $(\kappa (t))_{t \in \mathbb{R}_+}$ is defined by the continuous map $t \mapsto \kappa (t) = (t \wedge T^{\ast})^{q-1}$, where $T^{\ast}$ represents the average lifetime of the vortices in the flow. All these parameters, including the scale parameter  $\alpha \in (0,2)$,   are modulated with the shear value $\sigma$. 
This jump process admits the following time-dependent variance:
\[{\Var} (L_t)
= \int_0^t \int_{- \infty}^{+ \infty} | z |^2 \nu_s (d z) d s = \left(\tfrac{1}{q-1}  (t\wedge T^{
\ast})^{q} + (t\vee T^{\ast}) - T^{\ast} \right)  \ \left(\frac{z_{\ast}^{2 - \alpha}}{2 - \alpha} + \frac{{z^\ast}^{2 - \alpha}}{2 - \alpha}\right).\]
Denoting $\theta\mapsto a(\theta)$ and $\theta \mapsto b(\theta)$,  the  drift and diffusion coefficient in \eqref{eq:folded_theta}, the following cumulative angular dynamics 
\begin{equation}\label{eq:angular_with_jumps}
\theta_t = \theta_0 + \int_0^t \left( a(\theta_s)   +  \tfrac{1}{4} (b^2)''(\theta_s) \right)\,d s + \int_0^t b(\theta_{s^-}) \,d L_s
\end{equation}
inherits from the superdiffusivity   process $L$. 

We simulated the SDE \eqref{eq:angular_with_jumps} on the time interval $[0,1]$,  using the $\varepsilon$-EM scheme with Gaussian substitute \eqref{eq:numerical_scheme}, with parameters $n = 2^8$ and $\varepsilon = n^{-1}$, with $a$, $b$ as in \eqref{eq:folded_theta}, the $\gamma_i$ as reported in \cite[Table 1]{campana} for various values of $\sigma$. Sampling  $10^9$ Monte-Carlo trajectories, we estimated the time-evolution of the probability density function (PDF) of the renormalised modelled cumulative angle ${\delta\theta_t}/\sqrt{\Var{\theta_t}}$, aiming to replicate the PDF tails observed from measurements in \cite[Figure 8.]{campana}.

\begin{figure}[H]
\centering
\subfigure[isotropic case]{\includegraphics[width=0.37\paperwidth]{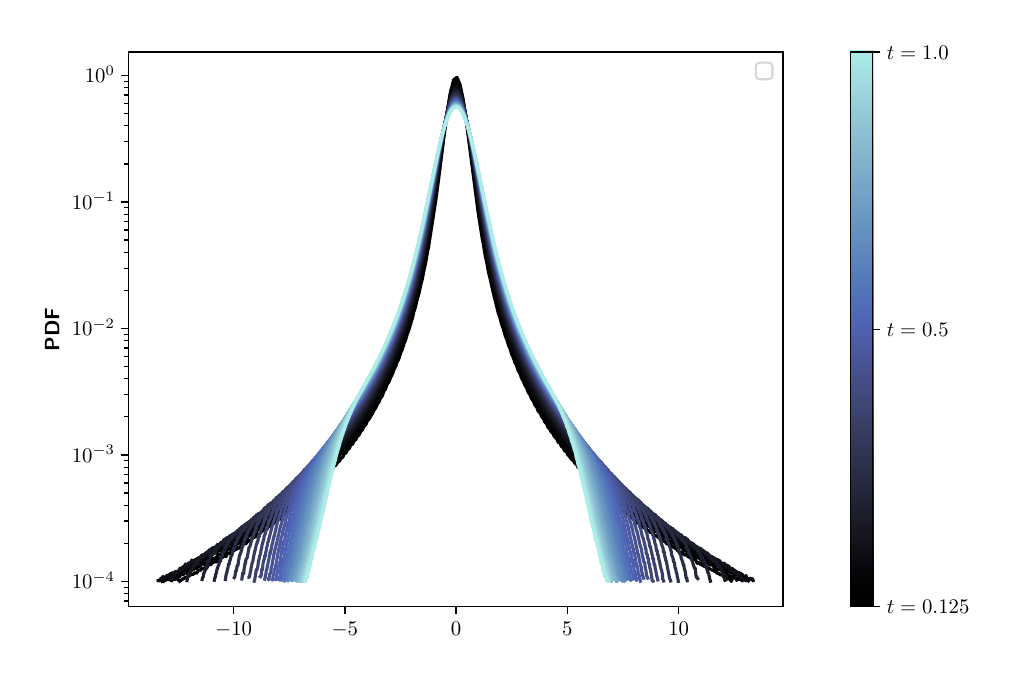}\label{subfig:pdf_shear0}}
\subfigure[with strong anisotropy]{\includegraphics[width=0.37\paperwidth]{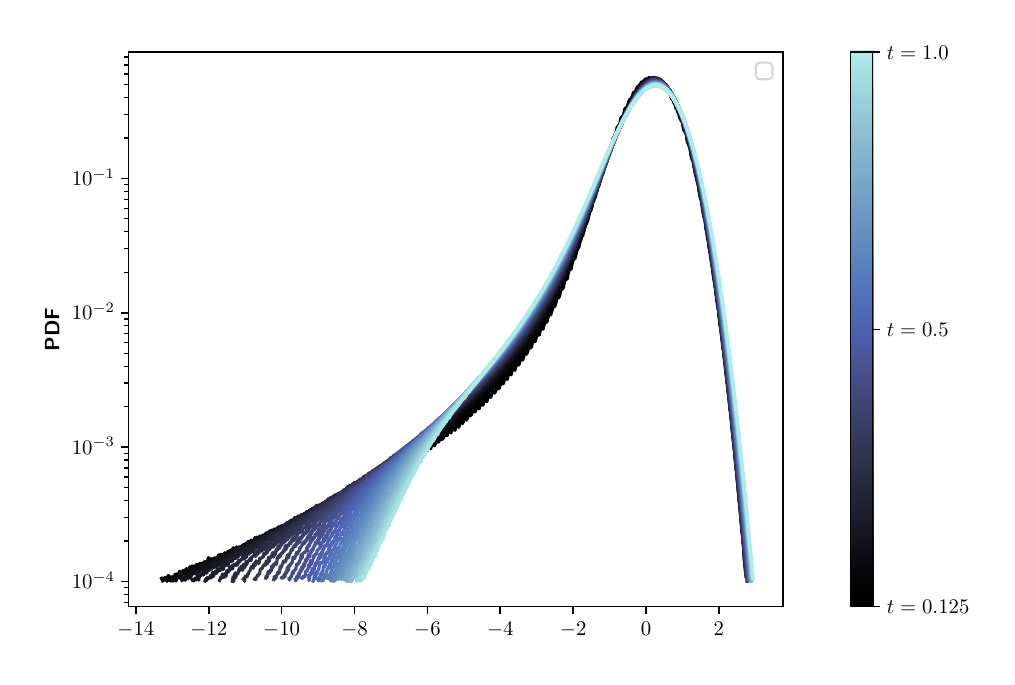}\label{subfig:pdf_shear8}}
\caption{\em 
Sampled PDF of the renormalised  cumulative angle ${\delta\theta_t}/\sqrt{\Var{\theta_t}}$  from  \eqref{eq:angular_with_jumps}, from $t=0.125$ to $t=1$. \\
{\rm (a)} $\sigma=0$,  $\alpha=3/2$ estimated on the DNS PDF (see \cite[Figure 8(a)]{campana}); $z_{\ast}=z^{\ast}=8$, $T^{\ast}=0.2$, $q=3/2$ estimated (and renormalised to fit final time 1) from the DNS variance map  \cite[Figure 7(a).]{campana}. \\[0.2cm]
{\rm (b)} $\sigma=2.8$,  $\alpha=1/2$ estimated on the DNS PDF (see \cite[Figure 8(b)]{campana}); $z_{\ast}=0$, $z^{\ast}=7$, $T^{\ast}=0.034$, $q=9/8$ estimated (and renormalised to fit final time 1) from the DNS variance map   \cite[Figure 7(a).]{campana}.
\label{fig:cluster_pdf}
}
\end{figure}

In Figure \ref{fig:cluster_pdf}, we present two typical cases: one without shear (isotropic case) and another with strong shear/anisotropy (with $\sigma=0$ and $\sigma=2.8$, respectively).  
For each $\sigma$, the specific values of the model parameters ($\alpha$, $(z_{\ast}, z^{\ast})$, $T^{\ast}$, $q$) are approximately estimated form the measurements in \cite{campana}. 
Although the estimators used are very crude, the results obtained with the model  \eqref{eq:angular_with_jumps} are qualitatively quite satisfactory. In Figure \ref{fig:cluster_pdf} (a and b), we observe a consistent pattern of wide tails that fold over time in a similar manner. 
There is still a great deal of physics to be understood, in particular to explain the modulation of the parameters as a function of shear. In this context, an even more refined stochastic approach involving continuous-time random walk processes would certainly be useful, although the analysis of this type of model is still incomplete and difficult. 

\section{Concluding remarks}

The proof of the strong convergence rate in Theorem \ref{theoreme_maurer} relies on the use of the Euler--Peano scheme as a pivot term. The comparison between the latter and the $\varepsilon$-EM scheme with substitute is only made possible by the use of the results of Bobkov, which are available in a discrete-time setting only, to the best of our knowledge. An eventual new result concerning bounds in $L^p$-Wasserstein distance between jumps SDEs could allow the strong convergence analysis of a continuous-time version of the $\varepsilon$-EM scheme, while making the use of the Euler--Peano scheme optional. Towards this direction, we point out again the work of Vlad Bally and Yifeng Qin in \cite{bally}, who showed a similar result in total variation distance.

The $\varepsilon$-EM scheme(s) may also be extended to a multidimensional form, and a strong convergence analysis could be done under similar proof structure (but with slightly more technicality), using for example Theorem 12 in \cite{bonis}, which extends the work of Bobkov to the multidimensional case. Regarding the weak error, we believe that the proof of Theorem \ref{thm:weak_error_bound} still holds, assuming the regularity of the flow for a multidimensional version of \eqref{eq:intro_SDEs}. Nevertheless, we point out that the proof of the latter from \cite{breton} relies on \cite{bichteler-etal-87} which is done in a one-dimensional setting only.

\medskip
In the recent literature on numerical schemes for jump SDEs, 
there is a growing interest in schemes with meshes adapted to jumps. For Lévy noise SDEs,  Dereich and Li \cite{dereich_li} give some $L^2(\Omega)$- and weak error asymptotic behaviours  in the case of  the $\varepsilon$-EM without substitute, intersecting a regular time-mesh, with the jump times of the cut noise.   Deligiannidis, Maurer and Tretyakov \cite{tretyakov} analyse the weak convergence for the $\varepsilon$-EM scheme with substitute with a restrictive jump-adapted time-stepping. In the same vein, for the case of finite activity jump diffusion,  Kelly, Lord and Sun \cite{kelly2023strong} analyse the strong error of a jump-adapted scheme. {In particular, this scheme allows them to handle the case of locally Lipschitz drift and diffusion coefficient.} We believe that performing a error analysis of the $\varepsilon$-EM scheme with substitute in the case of a jump-adapted scheme will be an interesting challenge.

\subsection*{Acknowledgements}
This work benefited from stimulating discussions with Jérémie Bec  and Lorenzo Campana who are warmly acknowledged. \smallskip

The authors are grateful to the OPAL infrastructure from Université Côte d'Azur and Inria for providing computational resources and support. \smallskip

The authors acknowledge the support of the French National Research Agency (ANR), under grant ANR-21-CE30-0040-01 (NETFLEX).

\appendix
\section{Appendix} \label{sec:appendix}
\subsection{Extended Gronwall Lemma}
\begin{lemma} \label{lem:generalized_gronwall}
  Let $f,\ell : [0, T] \rightarrow \mathbb{R}_+$ be
  non-negative non-decreasing functions and $g,h,k \in L^1([0,T],\mathbb{R_+})$. We assume that for some $p>1$, 
 \begin{equation}\label{eq:generalized_gronwall_hyp}
      \forall t \in [0, T], \text{ } f (t) \leq  \int_0^t  g(s) f (s) d s +  \left(
     \int_0^t h(s)  f^2 (s) d s \right)^{\frac{1}{2}} + \left( \int_0^t k (s) f^p
     (s) d s \right)^{\frac{1}{p}}  + \ell(t). 
\end{equation} 
Then
\begin{equation*}
\forall t \in [0, T], \text{ } f (t) \leq 2 p \ \ell(t)   \exp \left(
     \int_0^t (2 p \ g(s) + 2 p \ h(s) + k(s)) d s \right).
\end{equation*}
\end{lemma}
\label{sec:proof_gronwall}
\begin{proof}
Using the non-decreasing property of $f$ and then Young inequality, we have
 \begin{equation*}
\left( \int_0^t h (s) f^2 (s) d s \right)^{\frac{1}{2}}  \leq
f^{\frac{1}{2}} (t)  \  \left( \int_0^t h (s) f (s) d s\right)^{\frac{1}{2}} \leq \frac{1}{4} f (t) + \int_0^t h (s) f (s) d s
 \end{equation*}
 and
  \begin{equation*}
\left( \int_0^t k (s) f^p (s) d s \right)^{\frac{1}{p}}  \leq
f^{\frac{p - 1}{p}} (t)  \  \left( \int_0^t k (s) f (s) d s\right)^{\frac{1}{p}} \leq \frac{p - 1}{2 p} f (t) + \frac{2}{p}  \int_0^t k (s) f (s) d s
 \end{equation*}
 Combining this inequality with the hypothesis \eqref{eq:generalized_gronwall_hyp}, we obtain that for any $t\in[0,T]$,
  \[\frac{1}{2 p} f (t) \leq \int_0^t \left( g(s) + h(s) + \frac{2}{p} k(s) \right) f (s) d s + \ell (t). \]
Multiplying both sides by $2p$, we get
  \begin{align*}
    f (t) & \leq \int_0^t \left( 2 p \ g(s) + 2 p \ h(s) + k(s) \right) f(s) d s +2 p \ \ell (t), 
  \end{align*}
and we obtain the desired result by applying the standard Gronwall lemma.
\end{proof}

{
\subsection{Proof of Lemma \ref{lm:optimality_of_rate}} \label{sec:proof_lemma_rate}
Considering $X_t = 1 + \int_0^t X_{s^-} d L_s$, let $E_t = X_t - \wtX_t$, with $\wtX_t = 1 + \int_0^t \wtX_{\eta(s^-)} dL_s$ and $F_t = X_t - X_{\eta(t)}$, so that $E_t = \int_0^t (E_{\eta(s^-)} + F_{s^-}) d L_s .$ 
The Lévy--Itô decomposition yields $L_t = (\int_0^1 z \nu(dz)) t + \int_0^t \int_{-\infty}^{+\infty} z \rpm$ where $N$ is a Poisson random measure with intensity $\nu(dz) ds$. Applying \cite[Proposition 8.21]{tankov} we have $X_t = \prod_{0 \leq s \leq t} (1 + \Delta X_s)$. Since  $\nu$ has its support in $(0,+\infty)$,  $\Delta X_s \geq 0$ $\mathbb{P}$-a.s, which ensures that $X$ is positive and non-decreasing $\mathbb{P}$-a.s. In particular, for every $t\in[0,T]$, $F_t \geq 0$ $\mathbb{P}$-a.s. By induction,  $E_{t_i} \geq 0$,  using $E_{t_0} = 0$ and the formula $E_{t_i} = \sum_{j = 0}^{i-1} \int_{t_j}^{t_{j+1}} (E_{t_j} + F_{s^-}) d L_s$, for $i \in \{1,\ldots,n\}$.

Then, for every integer $p \geq 2$, the following inequality holds:
\begin{equation*}
    (E_t)^p = \bigg(\int_0^t (E_{\eta(s^-)} + F_{s^-}) d L_s\bigg)^p \geq \bigg(\int_0^t F_{s^-} d L_s\bigg)^p.
\end{equation*}
We apply the  Itô formula (e.g, \cite[Theorem 4.4.7]{apple}) to $\widetilde{E}= \int_0^t F_{s^-} dL_s$:
\begin{align*}
\widetilde{E}_t^{p} = & \bigg( \int_0^1 z \nu(dz) \bigg) p \int_0^t (\widetilde{E}_s)^{p-1} F_s \ ds \\ 
& +\int_0^t \int_{|z| < 1} ((\widetilde{E}_{s^-} 
+ F_{s^-}z)^{p} - \widetilde{E}_{s^-}^{p} - p z F_{s^-} \widetilde{E}_{s^-}^{p - 1}) \nu (d z) d s 
+ \int_0^t \int_{|z| < 1} ( (\widetilde{E}_{s^-} + F_{s^-}z)^{p} - \widetilde{E}_{s^-}^{p} ) \rpm.
\end{align*}
The first term of the right-hand side above  is bounded from below by zero, the second can be expressed by the mean of the binomial formula and the Fubini--Lebesgue Theorem, while the third term is a square integrable martingale. 
It follows that
\begin{align*}
\mathbb{E} [\widetilde{E}_t^{p}] \geq & \int_0^t \int_{|z| < 1} \left\{
\sum_{k = 2}^{p} \binom{p}{k} \mathbb{E} [\widetilde{E}_{s }^{p - k} F_{s }^k] z^k \right\}
\nu (d z) d s
= \int_0^t \left\{ \sum_{k = 2}^{p} \binom{p}{k} \left( \int_{0}^{1} z^k \nu (d z) \right)
\mathbb{E} [\widetilde{E}_s^{p - k} F_{s }^k] \right\} d s.
\end{align*}
Since every term in the above sum is positive, we retain only the critical one with $k=p$: 
\begin{equation*}
\mathbb{E} [\widetilde{E}_t^{p}] \geq \left( \int_{0}^{1} z^{p}
\nu (d z) \right)  \int_0^t \mathbb{E} [F_s^{p}] d s.
\end{equation*}
Then applying the It{\^o} formula once again on $F_t^{p} = \left( \int^t_{\eta (t)} X_{s^-} d L_s \right)^{p}$,
 and using previous arguments, 
\begin{equation*}
\mathbb{E} [F_t^{p}] 
\geq \left( \int_{0}^{1} z^{p}
\nu (d z) \right)  \int_{\eta (t)}^t \mathbb{E} [X_s^{p}] d s 
\geq \left( \int_{0}^{1} z^{p}\nu (d z) \right)^2  \int_0^t (s - \eta (s)) d s, 
\end{equation*}
since $\mathbb{E} [X_s^{p}] \geq  (\mathbb{E} [X_s])^{p} \geq  1$. 
In particular for $t = T$, 
\begin{align*}
\mathbb{E} [\widetilde{E}_{T}^{p}] \geq & 
\left( \int_{0}^{1} z^{p} \nu (d z) \right)^2 \ \sum_{i=0}^{n-1}\left(\int_{t_i}^{t_{i+1}}(s - t_i)ds\right) 
= \left( \int_{0}^{1} z^{p} \nu (d z) \right)^2 
\sum_{i=0}^{n-1} \frac{(t_{i+1}-t_i)^2}{2} 
= \frac{T^2}{2}  \left( \int_{0}^{1} z^{p} \nu (d z)
\right)^2 \ n^{- 1} .
\end{align*}
Letting $C_{p, T} = \left( \frac{T^2}{2}  \left( \int_{0}^{1}
    z^{p} \nu (d z) \right)^2 \right)^{\frac{1}{p}}$,  we obtain
\begin{equation*}
\Big\| \underset{t \in [0, T]}{\sup} |X_t - \wtX_t| \Big\|_{L^{ p} (\Omega)}  \geq \underset{t \in [0, T]}{\sup} \| X_t - \wtX_t \|_{L^{ p} (\Omega)}
\geq \| E_T \|_{L^{ p} (\Omega)} \geq \| \widetilde{E}_T \|_{L^{p} (\Omega)}
\geq \text{ } C_{p, T} n^{- \frac{1}{p}}.
\end{equation*}
}

\subsection{Technical details about the weak error simulation test case  \ref{sec:weak_error_arctan} } \label{sec:details_arctan}
In this section we detail the implementation of the $\varepsilon$-EM scheme in the case where $c(x,z)=\text{arctan}(x z)$. 

As mentioned in Section~\ref{sec:weak_error_arctan}, the numerical simulation of the weak error in this case requires first the computation of the function $G(t,x)$, and second  the computation of the variance  $\int_{t_{i-1}}^{t_i} \int_{B(\varepsilon)} c^2(s,\oeX_{t_{i - 1}},z) \rpc$ of the small jump integral.  Note that there is no need to compute the large jump martingale compensation $\int_{t_{i-1}}^{t_i} \int_{\RR \backslash B(\varepsilon)} c(s,\oeX_{t_{i - 1}},z) \rpc$ in our case, because this term is equal to zero due to the oddness of the   $\text{arctan}$ function.

\paragraph{Computation of $G(t,x)$.    } From the Kolmogorov equation, we get
\begin{equation*}
G(t,x)  = - 2x \cos(x) + \int_{-1}^{1} \Big\{\sin(x+\text{arctan}(x z)) - \sin(x) - \text{arctan}(x z) \cos(x)\Big\} \;\frac{dz}{|z|}.
\end{equation*}
With the remark that $\cos(\text{arctan}(y)) = \frac{1}{\sqrt{1+y^2}}$ and $\sin(\text{arctan}(y)) = \frac{y}{\sqrt{1+y^2}}$ for any $y \in \RR$ and using the sinus addition formula, we obtain that 
\begin{align*}
    \sin(x+\text{arctan}(x z)) - \sin(x) - \text{arctan}(x z) \cos(x) & = \sin(x) \left(\frac{1}{\sqrt{1+x^2 z^2}}-1\right) + \cos(x) \frac{x z}{\sqrt{1+x^2 z^2}}. 
\end{align*}
Then by imparity, one has
$\int_{-1}^{1} \frac{x z}{\sqrt{1+x^2 z^2}} \frac{1}{|z|} dz = 0$, 
and one may compute explicitly the integral
\begin{equation*}
    \int_{-1}^{1} \left(\frac{1}{\sqrt{1+x^2 z^2}}-1\right) \frac{1}{|z|} dz = 2\big(\log(2) - \log(\sqrt{x^2 + 1} +1 )\big),
\end{equation*}
leading to
\begin{equation*}
    G(t,x) = - 2 x \cos(x) + 2 \sin(x) \big(\log(2) - \log(\sqrt{x^2 + 1} +1 )\big).
\end{equation*}

\paragraph{Computation of $\int_{t_{i-1}}^{t_i} \int_{B(\varepsilon)} c^2(s,\oeX_{t_{i - 1}},z) \rpc$.    }  We now focus on the approximation of the variance of the small jump integral, which in our case writes
\begin{equation*}
    \frac{1}{n_{\varepsilon}} \int_{-\varepsilon}^{\varepsilon} \frac{\text{arctan}(\oeX_{t_{i - 1}},z)^2}{|z|} dz  =  \frac{2}{n_{\varepsilon}} \int_{0}^{\varepsilon |\oeX_{t_{i - 1}}|} \frac{\text{arctan}(z)^2}{z} dz. 
\end{equation*}
The latter integral does not have a simple closed form. However, when the upper-boundary  $\varepsilon |\oeX_{t_{i - 1}}|$ is small enough, we may approximate $\text{arctan}$ by its Taylor series near zero, which leads to the following approximation strategy: For each step $t_i = i \frac{T}{n_{\varepsilon}}$, given a threshold $\tau$, 
\begin{description}
    \item[\qquad $\bullet$ if $\varepsilon |\oeX_{t_{i - 1}}| < \tau$: ] we approximate $\text{arctan}(z)$ by its Taylor series at order 11.
    
    With $\tau=\tfrac12$, this leads  to an approximation error of $\frac{(0.5)^{2\times6+1}}{2\times6+1} \simeq 10^{-5}$ in the worse case, which is negligible with respect to our Monte-Carlo error.
    
    \item[\qquad $\bullet$ Otherwise: ]  separate the integral in two
    \begin{equation*}
        \int_{0}^{\varepsilon |\oeX_{t_{i - 1}}|} \frac{\text{arctan}(z)^2}{z} dz = \int_{0}^{\tau} \frac{\text{arctan}(z)^2}{z} dz + \int_{\tau}^{\varepsilon |\oeX_{t_{i - 1}}|} \frac{\text{arctan}(z)^2}{z} dz.
    \end{equation*}
    The first integral is approximated by the previous Taylor series method, while we use a $\tfrac{1}{3}$-Simpson rule to obtain the second. With discretisation step $m = 10^2$, this leads to an approximation error that behaves as $10^{-8} \times \sup_{q \in [\frac{1}{2},+\infty)}|\frac{d^4}{dq^4} \frac{\text{arctan}(q)^2}{q} |$, where the latter supremum can be shown numerically to be smaller than 10.
\end{description}
For $\varepsilon$ small enough, the second case will only appear in a few extreme trajectories of the Monte-Carlo simulations, which limits the impact of this additional numerical integration to the total computation time.

\end{document}